\documentclass[12pt,oneside,reqno]{amsart}

\usepackage{amssymb}
\usepackage{amsfonts}
\usepackage{bbm}
\usepackage{cases}
\usepackage{amsmath}
\usepackage{graphicx}
\usepackage{mathrsfs}
\usepackage{stmaryrd}
\usepackage{txfonts}
\usepackage{enumerate,amsmath,amssymb,amsthm,booktabs}

\numberwithin{equation}{section}

\newcommand{\be}{\begin{eqnarray}}
\newcommand{\ee}{\end{eqnarray}}
\newcommand{\ce}{\begin{eqnarray*}}
\newcommand{\de}{\end{eqnarray*}}
\newtheorem{theorem}{Theorem}[section]
\newtheorem{lemma}[theorem]{Lemma}
\newtheorem{remark}[theorem]{Remark}
\newtheorem{definition}[theorem]{Definition}
\newtheorem{proposition}[theorem]{Proposition}
\newtheorem{Examples}[theorem]{Example}
\newtheorem{corollary}[theorem]{Corollary}

\def\eps{\varepsilon}

\def\ga{{r}}
\def\e{\mathrm{e}}

\def\p{\partial}

\def\[{{\Big[}}
\def\]{{\Big]}}
\def\<{{\langle}}
\def\>{{\rangle}}
\def\({{\Big(}}
\def\){{\Big)}}

\def\bx{{\mathbf{x}}}
\def\tr{\mathrm {tr}}

\def\dif{{\mathord{{\rm d}}}}

\def\no{\nonumber}
\def\={&\!\!=\!\!&}

\def\cB{{\mathcal B}}

\def\cE{{\mathcal E}}
\def\cF{{\mathcal F}}

\def\cM{{\mathcal M}}

\def\cT{{\mathcal T}}

\def\mA{{\mathbb A}}
\def\mB{{\mathbb B}}

\def\mE{{\mathbb E}}

\def\mG{{\mathbb G}}
\def\mH{{\mathbb H}}
\def\mI{{\mathbb I}}

\def\mK{{\mathbb K}}
\def\mL{{\mathbb L}}
\def\mM{{\mathbb M}}
\def\mN{{\mathbb N}}

\def\mP{{\mathbb P}}
\def\mQ{{\mathbb Q}}
\def\mR{{\mathbb R}}

\def\mV{{\mathbb V}}
\def\mW{{\mathbb W}}

\def\1{{\mathbf{1}}}

\def\sB{{\mathscr B}}

\def\sF{{\mathscr F}}

\def\sL{{\mathscr L}}

\def\geq{\geqslant}
\def\leq{\leqslant}

\def\eps{\varepsilon}

\def\e{\mathrm{e}}

\def\p{\partial}

\def\[{{\Big[}}
\def\]{{\Big]}}
\def\<{{\langle}}
\def\>{{\rangle}}
\def\({{\Big(}}
\def\){{\Big)}}

\def\bx{{\mathbf{x}}}
\def\tr{\mathrm {tr}}

\def\dif{{\mathord{{\rm d}}}}

\def\no{\nonumber}
\def\={&\!\!=\!\!&}
\def\bt{\begin{theorem}}
\def\et{\end{theorem}}
\def\bl{\begin{lemma}}
\def\el{\end{lemma}}
\def\br{\begin{remark}}
\def\er{\end{remark}}
\def\bx{\begin{Examples}}
\def\ex{\end{Examples}}
\def\bd{\begin{definition}}
\def\ed{\end{definition}}
\def\bp{\begin{proposition}}
\def\ep{\end{proposition}}
\def\bc{\begin{corollary}}
\def\ec{\end{corollary}}

\def\geq{\geqslant}
\def\leq{\leqslant}

 \def\R{\mathbb R}
 \def\R{\mathbb R}

\def\<{\langle} \def\>{\rangle}

 \def\beq{\begin{equation}}  
 
\def\e{\text{\rm{e}}}

\allowdisplaybreaks

\begin{document}

\title{Ergodicity of stochastic differential equations with jumps and singular coefficients}
\date{}

\author{ LONGJIE XIE\ \ and \  XICHENG ZHANG}

\address{Longjie Xie:
School of Mathematics and Statistics, Jiangsu Normal University,
Xuzhou, Jiangsu 221000, P.R.China\\
Email: xlj.98@whu.edu.cn
 }
\address{Xicheng Zhang:
School of Mathematics and Statistics, Wuhan University,
Wuhan, Hubei 430072, P.R.China\\
Email: XichengZhang@gmail.com
 }

\thanks{Research of L. Xie is partially supported by the Project
Funded by the PAPD of Jiangsu Higher Education Institutions. X. Zhang is supported by NNSF of China (No. 11325105).}

\begin{abstract}
We show the strong well-posedness of SDEs driven by general multiplicative L\'evy noises with Sobolev diffusion
and jump coefficients and integrable drift. Moreover, we also study
the strong Feller property, irreducibility as well as the exponential ergodicity of the corresponding semigroup when
the coefficients are time-independent and singular dissipative.
In particular, the large jump is allowed in the equation.
To achieve our main results, we present a general approach for treating the SDEs with jumps 
and singular coefficients so that one just needs to focus on Krylov's {\it apriori} estimates for SDEs.

\bigskip


\end{abstract}

\maketitle \rm

\tableofcontents

\section{Introduction}

Let $(\Omega,\sF,(\sF_t)_{t\geq 0},\mP)$ be a filtered probability space, which satisfies the usual conditions.
On this probability space, let $(W_t)_{t\geq 0}$ be a $d$-dimensional standard $\sF_t$-Brownian motion and $N$ an $\sF_t$-Poisson random measure
with intensity measure $\dif t\nu(\dif z)$, where $\nu$ is a L\'evy measure on $\mR^d$, that is,
\begin{align}
\int_{\mR^d}\!\big(|z|^2\wedge1\big)\nu(\dif z)<+\infty,\ \nu(\{0\})=0. \label{nu}
\end{align}
The compensated Poisson random measure $\tilde N$ is defined as
$$
\tilde{N}(\dif t,\dif z):= N(\dif t,\dif z)-\dif t\nu(\dif z).
$$
Consider the following stochastic differential equation (SDE) in $\mR^d$ with jumps:
\begin{align}\label{sde1}
\dif X_t=\sigma_t(X_t)\dif W_t+b_t(X_t)\dif t+\!\!\int_{|z|<R}g_t(X_{t-},z)\tilde N(\dif t,\dif z)+\!\!\int_{|z|\geq R}g_t(X_{t-},z)N(\dif t,\dif z),
\end{align}
where $R>0$ is a fixed constant, and $\sigma: \mR_+\times\mR^d\rightarrow\mR^d\otimes\mR^d$, $b: \mR_+\times\mR^d\rightarrow\mR^d$
and $g: \mR_+\times\mR^d\times\mR^d\rightarrow\mR^d$ are Borel measurable functions,
which are called diffusion, jump and drift coefficients, respectively.
Recall that an $\sF_t$-adapted c\`adl\`ag (right continuous with left limit)
process $X$ is called a (strong) solution of SDE \eqref{sde1} if
for each $t>0$, the following random variables are finite $\mP$-almost surely,
$$
\int^t_0\!\!\|\sigma_s(X_s)\|^2\dif s,\int^t_0\!\!|b_s(X_s)|\dif s, \int^t_0\!\!\!\int_{|z|<R}\!\!|g_s(X_s,z)|^2\nu(\dif z)\dif s,
\int^t_0\!\!\!\int_{|z|\geq R}\!\!|g_s(X_s,z)|\nu(\dif z)\dif s,
$$
and
\begin{align*}
X_t&=X_0+\int^t_0\sigma_s(X_s)\dif W_s+\int^t_0b_s(X_s)\dif s+\int^t_0\!\!\!\int_{|z|<R}g_s(X_{s-},z)\tilde N(\dif s,\dif z)\\
&\qquad\quad+\int^t_0\!\!\!\int_{|z|\geq R}g_s(X_{s-},z)N(\dif s,\dif z),\ \mP-a.s.
\end{align*}
In this paper, one of our aims is to show the existence and uniqueness of a solution to the above SDE under some mild assumptions on the coefficients,
both in the non-degenerate diffusion case and in the multiplicative pure jump case.
Moreover, we also study the strong Feller property, irreducibility as well as the ergodicity of the semigroup associated with the above SDE
when the coefficients are time-independent and singular dissipative.

\vspace{2mm}
In the past decades, SDEs with singular drifts and driven by Brownian motion have been extensively studied.
In the case that $g\equiv0$ and $\sigma\equiv\mI_{d\times d}$, the identity matrix, a remarkable result due to Krylov and R\"ockner \cite{Kr-Ro}
says that SDE (\ref{sde1}) has a unique strong solution provided that
\begin{align*}
b\in L^q_{loc}\big(\mR_+;L^p(\mR^d)\big)\quad\text{with}\quad \tfrac{d}{p}+\tfrac{2}{q}<1.
\end{align*}
Latter, the second named author \cite{Zh3, Zh1} extended their result to the multiplicative noise case
under some non-degenerate and Sobolev conditions on the diffusion coefficient.
On the other hand, by studying the stochastic homeomorphism flow property of the SDEs with irregular drifts, Flandoli, Gubinelli and Priola \cite{F-G-P} obtained
a well-posedness result for a class of stochastic transport equations with irregular coefficients. After that,
there are many works devoted to the study of the regularities of the unique strong solution
to SDEs with rough coefficients, such as the Sobolev differentiability with respect to the initial value, stochastic homeomorphism flow
and the Malliavin differentiability with respect to the sample path.
The interested readers are referred to \cite{Fe-Fl-1,Fe-Fl-2,M-N-P-Z,Mo-Ni-Pr,W,Zh1,Zh4} and references therein.

\vspace{2mm}

In recent years, SDEs driven by pure jump L\'evy processes (i.e., $\sigma\equiv0$) and with irregular drifts have also attracted great interests
since it behaves quite differently.
In fact, when $d=1$ and $(L_t)_{t\geq 0}$ is a symmetric $\alpha$-stable process with $\alpha\in(0,1)$,
Tanaka, Tsuchiya and Watanabe \cite{Ta-Ts-Wa} showed that even if $b$ is time-independent, bounded and $\beta$-H\"older continuous
with $\beta<1-\alpha$, SDE
\begin{align}
\dif X_t=\dif L_t+b(X_t)\dif t,\quad X_0=x\in\mR^d \label{levy}
\end{align}
may not have a pathwise uniqueness strong solution, see also \cite{B-B-C} for related result. On the other hand, when $\alpha\in[1,2)$ and
$$
b\in C^\beta_b(\mR^d)\quad\text{with}\quad\beta>1-\tfrac{\alpha}{2},
$$
it was shown by Priola \cite{Pri} that there exists a unique strong solution $X_t(x)$ to SDE (\ref{levy}) for each $x\in\mR^d$,
which forms a stochastic $C^1$-diffeomorphism flow. Under the same condition, Haadem and Proske \cite{H-F} obtained the unique strong solution by using the Malliavin calculus. Recently, Zhang \cite{Zh00} obtained the pathwise uniqueness to SDE (\ref{levy}) when $\alpha\in(1,2)$, $b$ is bounded and in some fractional Sobolev space.
See also \cite{B-P-1,Ch-So-Zh,Pri2,Pri3} for related results. It is noticed that all the works mentioned above for SDE (\ref{sde1}) with $\sigma\equiv0$ are restricted to the additive noise case. We also mention that Bogachev and Pilipenko
\cite{B-P-2} treated the SDE with general L\'evy noise and discontinuous drift based on the heat kernel estimate.

\vspace{2mm}
In this paper, we shall first study the well-posedness of general SDEs with Sobolev diffusion and jump coefficients and integrable drifts.
In the mixing and non-degenerate diffusion case, we shall not make any assumptions on the pure jump L\'evy noise (or the L\'evy measure $\nu$ in (\ref{nu})).
Our result extends the existing results concerning singular SDEs driven by Brownian motion (see \cite{B-P-2,Ve,Kr-Ro,Zh1}), see Theorem \ref{main1}.
In the pure jump case, we shall assume that the L\'evy measure $\nu$ is symmetric and rotationally invariant $\alpha$-stable type,
which is because we need to use the heat kernel estimate established in \cite{Ch-Zh}, see Theorem \ref{main3}. Compared with \cite{H-F,Pri,Pri2,Pri3,Zh00}, we
are considering the multiplicative noise and drop the boundedness assumption on drift $b$.

\vspace{2mm}

Now we introduce the main argument adopted below: Zvonkin's transformation. Let $\sL_2^\sigma$ be the second order differential operator associated with the diffusion coefficient $\sigma$, that is,
$$
\sL_2^\sigma u(x):=\tfrac{1}{2}(\sigma^{ik}_t\sigma^{jk}_t\p_{i}\p_ju)(x).
$$
Here and below, we use Einstein's convention that the repeated indices in a product will be summed automatically.
Let $\sL_1^b$ be the first order differential operator associated with the drift coefficient $b$, that is,
$$
\sL_1^bu(x):=(b^i_t\p_iu)(x),
$$
and $\sL^g_{\nu}$ the nonlocal operator associated with the jump coefficient $g$,\! that is,
\begin{align}\label{UY7}
\begin{split}
\sL^g_{\nu} u(x)&:=\int_{|z|<R}\Big[u\big(x+g_t(x,z)\big)-u(x)-g_t(x,z)\cdot\nabla u(x)\Big]\nu(\dif z)\\
&\quad+\int_{|z|\geq R}\Big[u\big(x+g_t(x,z)\big)-u(x)\Big]\nu(\dif z)=:\sL^g_{\nu,R} u(x)+\bar\sL^g_{\nu,R} u(x).
\end{split}
\end{align}
Fix a time $T>0$. We consider the following Kolmogorov's backward  equation system:
\begin{align}\label{DE1}
\p_t \Phi+(\sL^\sigma_2+\sL^b_1+\sL^g_\nu)\Phi=0,\ \ \Phi_T(x)=x\in\mR^d.
\end{align}
Suppose that this equation has a regular enough solution $\Phi$ so that for each $t\in[0,T]$, the map $x\mapsto\Phi_t(x)$ forms a $C^2$-diffeomorphism on $\mR^d$.
By  It\^o's formula, one sees that
\begin{align*}
\Phi_t(X_t)&=\Phi_0(X_0)+\int^t_0\nabla\Phi_s(X_s)\sigma_s(X_s)\dif W_s\\
&+\int^t_0\!\!\!\int_{|z|<R}(\Phi_s(X_{s-}+g_s(X_{s-},z))-\Phi_s(X_{s-}))\tilde N(\dif s,\dif z)\\
&+\int^t_0\!\!\!\int_{|z|\geq R}(\Phi_s(X_{s-}+g_s(X_{s-},z))-\Phi_s(X_{s-}))N(\dif s,\dif z).
\end{align*}
Thus, if we let $Y_t:=\Phi_t(X_t)$ and
\begin{align}\label{HD7}
\tilde\sigma_t(y):=(\nabla\Phi_s\cdot\sigma_s)\circ\Phi^{-1}_t(y),\ \ \tilde g_t(y,z):=\Phi_t(\Phi^{-1}_t(y)+g_t(\Phi^{-1}_t(y),z))-y,
\end{align}
then $Y_t$ satisfies the following {\it new} SDE with {\it disappeared} drift:
\begin{align}\label{DE2}
\dif Y_t=\tilde\sigma(Y_t)\dif W_t+\int_{|z|<R}\tilde g_t(Y_{t-},z)\tilde N(\dif s,\dif z)+\int_{|z|<R}\tilde g_t(Y_{t-},z)\tilde N(\dif t,\dif z).
\end{align}
And vice versa, if $Y_t$ solves \eqref{DE2}, then $X_t:=\Phi^{-1}_t(Y_t)$ solves SDE \eqref{sde1}.
In the case $g$=0, if $\sigma$ is uniformly elliptic and Lipschitz continuous,
then $\tilde\sigma$ could be also Lipschitz continuous  due to the second order regularization effect of equation \eqref{DE1} even for H\"older $b$, see \cite{Zh1, Zh4}.
Thus, our main task is to solve equation \eqref{DE1} so that $\Phi$ has the desired properties.

\vspace{2mm}

However, for $b\in L^q_{loc}(\mR_+;L^p(\mR^d))$ being not H\"older continuous,
the transformed coefficients $\tilde\sigma$ and $\tilde g$ in \eqref{DE2}
are not expected to be Lipschitz continuous, but at most in the first order Sobolev space $\mW^{1,p}_{loc}(\mR^d)$.
In other words, we need to solve SDE \eqref{DE2} with coefficients being in $\mW^{1,p}_{loc}(\mR^d)$.
To this purpose, a key step is to show the following {\it apriori} Krylov's estimate:
for {\it any} solution $Y$, and any $T>0$ and $f\in L^q_{loc}\big(\mR_+;L^p(\mR^d)\big)$,
\begin{align}\label{KRY}
\mE\left(\int^T_0 f(t,Y_t)\dif t\right)\leq c\left(\int^T_0\left(\int_{\mR^d}|f(t,x)|^p\dif x\right)^{q/p}\dif t\right)^{1/q}.
\end{align}
For general continuous It\^o's process, such an estimate was established in \cite{Kr1} for $p=q\geq d+1$. For SDE \eqref{sde1} with $g=0$ and general $p,q$ with
$\frac{d}{p}+\frac{2}{q}<2$,  we refer to \cite{Zh1, Zh4, Zh5}.
However, for discontinuous semimartingales, there are few results. In Section 5, we shall devote to a detailed study about the above Krylov estimate
for any solution of SDE \eqref{sde1} under mild conditions.
In the non-degenerate diffusion case, we first use a Krylov's lemma (see Lemma \ref{pdi}) to show \eqref{KRY} for any solution of SDE \eqref{sde1} and
for any $p=q\geq d+1$. For general $p,q$ with $\frac{d}{p}+\frac{2}{q}<2$, we
need to first solve a non-homogeneous Kolmogorov's backward equation,
see Theorem \ref{pde}, and then use Girsanov's theorem,
see Theorem \ref{enou}.
While in the purely non-local $\alpha$-stable-like case, we shall use Duhamel's formula and
heat kernel estimates recently obtained in \cite{Ch-Zh, Ch-Zh0}
to solve a non-homogeneous nonlocal Kolmogorov's backward equation, see Theorem \ref{Th63}. Furthermore, by smoothing out the coefficient
we show \eqref{KRY} for any $p,q$ with $\frac{d}{p}+\frac{\alpha}{q}<\alpha$, see Theorem \ref{Th65}.

\vspace{2mm}
Although the well-posedness and regularity properties of strong solutions for SDEs with singular coefficients have been intensively studied,
it seems that there are less works devoted to studying the existence and uniqueness of invariant probability measures for time-independent SDEs with singular coefficients.
As we know, a general approach of proving the existence of invariant probability measures
is to verify the Lyapunov condition. More precisely, if there exists a positive function $\Phi_1\in C^2(\mR^d)$ and a positive compact function $\Phi_2$ such that
\begin{align}\label{Lya}
(\sL^\sigma_2+\sL^b_1+\sL^g_\nu)\Phi_1\leq C-\Phi_2
\end{align}
holds for some constant $C>0$,
then the associated semigroup of SDE \eqref{sde1} has an invariant probability measure $\mu$ with $\mu(\Phi_2)<\infty$, see for instance \cite{Has}.
Obviously, if $b\in L^p(\mR^d)$ for some $p>d$, then compact function $\Phi_2$ would not exist since $b$ can be singular at infinity.
In this direction, to the authors' knowledge, Wang \cite{Wa} obtained a first result about the ergodicity for SDEs with singular drifts by using perturbation argument
and his local dimension-free Harnack inequality. In particular, the main result in \cite{Wa} is applied to the following singular SDE so that
it admits a unique invariant probability measure:
$$
\dif X_t=(b(X_t)-\lambda_0 X_t)\dif t+\sqrt{2}\dif W_t,\ X_0=x\in\mR^d,
$$
where $\lambda_0>0$ and $b:\mR^d\to\mR^d$ satisfies
\begin{align}\label{DE3}
\int_{\mR^d}\e^{\lambda |b(x)|^2-\lambda_0|x|^2/2}\dif x<\infty \mbox{  for some $\lambda>\frac{1}{2\lambda_0}$}.
\end{align}
To prove the uniqueness of invariant measures, a usual way is to show the strong Feller property and irreducibility of the associated semigroup.
In \cite{XZ}, we have studied these two properties for SDEs driven by Brownian motion under some local conditions.

\vspace{2mm}

In the pure jump case, when the coefficients are locally Lipschitz continuous and satisfy a Lyapunov type dissipative condition such as \eqref{Lya},
it has been shown in \cite{Ku,Ma} that there is a unique invariant probability measure associated to the SDE, which is exponential ergodicity.
In \cite{Ku} and \cite{Ma}, they introduced some abstract conditions for the ergodicity.
In a recent work \cite{Ar-Bi-Ca}, the authors also introduced some Lyapunov stability condition
for the existence of invariant probability measures for general nonlocal operators. Clearly, the singular drift does not satisfy their Lyapunov conditions.

\vspace{2mm}

The second aim of this paper is to show the existence and uniqueness of invariant probability measures associated to SDE \eqref{sde1}
both in non-degenerate diffusion case and in pure jump case under suitable singular and dissipative assumptions. Our basic idea is as follows:
Suppose that $b$ can be decomposed into two parts:
$$
b=b_1+b_2,
$$
where $b_1$ is the singular part and $b_2$ is the dissipative part, see {\bf (H$^b$)} and {\bf ($\widetilde {\mathbf H}^b$)} below.
For the singular part $b_1$, we use Zvonkin's transformation to kill it and obtain a new SDE, which preserves the dissipativity fortunately.
Of course, to perform Zvonkin's transformation, here we need to solve a nonlocal elliptic equation rather than the parabolic equation \eqref{DE1}.
Moreover, we need to show the non-explosion and the apriori Krylov estimate for any solution of SDE \eqref{sde1} with singular and dissipative drift,
see Lemma \ref{Le74} and Lemma \ref{Le78}.
The following table figures out the methods of showing the existence and uniqueness of invariant probability measures for
time-independent SDE \eqref{sde1} with dissipative drift:
\begin{align*}
\begin{tabular}{ccc}
\toprule\midrule
Existence of IPMs &Strong Feller&Irreducibility\\ \midrule
&\textsf{ Diffusion with jump}&\\ \midrule
Lyapunov condition& Derivative formula & Coupling+GT\\ \midrule
&\textsf{  Pure jump SDE}&\\ \midrule
Lyapunov condition & Continuity of HK & Positivity of DHK \\ \midrule
\bottomrule
\end{tabular}
\end{align*}
Here, IPM, GT, HK and DHK stands for invariant probability measure, Girsanov's transform, heat kernel and Dirichlet heat kernel, respectively.

\vspace{2mm}

Below we provide two simple examples to illustrate the main results obtained in this paper.
\bx
{\rm Consider the following SDE of OU type:
$$
\dif X_t=\dif L_t-\lambda_0X_t\dif t+b(X_t)\dif t,\quad X_0=x\in\mR^d.
$$
When $L_t$ is a $d$-dimensional standard Brownian motion, we assume $b\in L^p(\mR^d)$ for some $p>d$.
When $L_t$ is a rotationally invariant symmetric $\alpha$-stable process with $\alpha\in(1,2)$, we assume
$b\in H^{\theta}_p(\mR^d)$ for some $\theta>1-\alpha/2$ and $p>2d/\alpha$, where $H^{\theta}_p(\mR^d)$ is the Bessel potential space. 
By Theorems \ref{main2} and \ref{main4},
the above SDE admits a unique strong solution
and there exists a unique invariant measure associated with it.
Note that  in both cases the classical Lyapunov condition (\ref{Lya}) is not satisfied, our result is new even in the existence of invariant measures.
Moreover, compared with Wang's global condition \eqref{DE3}, our global assumption $b\in L^p(\mR^d)$ is weaker locally, but the singularity is not comparable at infinity.
}
\ex

\bx
{\rm Consider the following mixing SDE with jumps:
$$
\dif X_t=\dif W_t+\lambda_1|X_{t-}|^\beta\dif L_t-\lambda_0X_t|X_t|^{\gamma-1}\dif t,\ \ X_0=x\in\mR^d,
$$
where $\beta\in(0,1)$, $\gamma\in(0,\infty)$ and $\lambda_0>0$, $\lambda_1\in\mR$,
$L_t$ is a $d$-dimensional pure jump L\'evy process. The main features of this SDE are that the jump coefficient $x\mapsto|x|^\beta$ is H\"older continuous
and the drift term can be polynomial growth.
By Theorem \ref{main22}, the above SDE has a unique strong solution. Moreover, there exists a unique invariant probability
measure which is $V$-ergodicity (see Definition \ref{de}) in the case $\gamma\in(0,1]$ and exponential ergodicity in the case $\gamma>1$.
}
\ex

Finally, recall that a probability measure $\mu$ on $\mR^d$ is called an invariant probability measure
of operator $\sL:=\sL_2^\sigma+\sL_1^b+\sL_\nu^g$ if it satisfies the following Fokker--Planck--Kolmogorov equation
\begin{align}
\sL^*\mu=0\Leftrightarrow\mu(\sL \varphi)=0,\ \ \varphi\in C_0^{\infty}(\mR^d), \label{meas}
\end{align}
where the asterisk stands for the formal adjoint operator. Obviously, any invariant probability measure of the semigroup associated with SDE \eqref{sde1}
satisfies \eqref{meas}. When $g\equiv0$, the existence of solutions to (\ref{meas}) was obtained in \cite{B-P-R}
by analytic methods under a Lyapunov-type condition, which is much weaker than those needed for
the existence of a solution to SDE (\ref{sde1}). Moreover, under some quite weak conditions,
the uniqueness and regularities of the solutions for \eqref{meas} are also studied in \cite{B-K-R-1, B-K-R-2,B-R-S}, see also \cite{Wa}.
To our knowledge, these results cannot cover our results stated above.

\vspace{2mm}

This paper is organized as follows: In Section 2, the main results including the existence-uniqueness and ergodicity for SDE \eqref{sde1} are stated.
In Section 3, after introducing the notion of Krylov's estimate,
we present two general results: Stability and Zvonkin's transformation for SDE \eqref{sde1}. Moreover, we also prove a useful stochastic Gronwall's inequality,
which extends Scheutzow's result \cite{Sc}. In Section 4, we study the regularities of parabolic integral-differential equations. In Section 5,
applying the results obtained in the previous section, we show various Krylov's estimates for the solution of SDE \eqref{sde1}.
By the general results in Section 3, the strong well-posedness results are proved in Section 6.
Finally, the strong Feller property and irreducibility as well as the ergodicity for SDE \eqref{sde1} are proven in Section 7.

\vspace{2mm}

Throughout this paper, we use the following convention: $c$ with or without subscripts will denote a positive constant, whose value may change in different places.
Moreover, we use $A\lesssim B$ to denote $A\leq cB$ for some constant $c>0$.

\section{Statement of main results}

\subsection{Strong well-posedness of singular SDEs with jumps}
To state our main results, we first introduce some spaces and notations.
For $p,q\in[1,\infty]$ and $0\leq S<T<\infty$,
let $\mL^q_p(S,T)$  be the space of all Borel functions on $[S,T]\times\mR^d$ with norm
$$
\|f\|_{\mL^q_p(S,T)}:=\Bigg(\int_S^T\!\Bigg(\int_{\mR^d}|f(t,x)|^p\dif x\Bigg)^{q/p}\dif t\Bigg)^{1/q}<\infty.
$$
For $p=\infty$ or $q=\infty$, the above norm is understood as the usual $L^\infty$-norm. We shall simply write
$$
\mL^q_p(T):=\mL^q_p(0,T),\ \ \mL^p(T):=\mL^p_p(T).
$$
Given a $R>0$, we shall write $B_R:=\{x\in\mR^d: |x|<R\}$. For a measurable function $g_t(x,z):\mR_+\times\mR^d\times\mR^d\to\mR^d$ and $0\leq\eps<R\leq\infty$, we introduce
the following functions, which will be used frequently below:
for $j=0,1$ and $\alpha\geq 1$,
\begin{align}\label{LK1}
 \Gamma^{j,\alpha}_{\eps, R}(g)(t,x):= \Gamma^{j,\alpha}_{\eps, R}(g_t)(x):=\|\nabla^j_xg_t(x,\cdot)\|^\alpha_{L^\alpha(B_R\setminus B_\eps;\nu)}
:=\! \! \int_{\eps\leq |z|<R}\!\!\!|\nabla^j_x g_t(x,z)|^\alpha \nu(\dif z).
\end{align}
Here and below, $\nabla_x$ denotes the generalized gradient with respect to $x$.

We make the following assumptions on  the diffusion coefficient $\sigma$:
\begin{enumerate}
\item [{\bf (H$^\sigma$)}] There are constants $c_0\geq 1$ and $\beta\in(0,1)$ such that for all $(t,x)\in\mR_+\times\mR^d$,
$$
c_0^{-1}|\xi|^2\leq |\sigma^*_t(x)\xi|^2\leq c_0|\xi|^2,\quad\forall \xi\in\mR^d,
$$
where $\sigma^*$ stands for the transpose of $\sigma$, and
$$
\|\sigma_t(x)-\sigma_t(x')\|\leq c_0|x-x'|^\beta.
$$
Here and below, $\|\cdot\|$ denotes the Hilbert-Schmidt norm of a matrix.
\end{enumerate}

Our first main result of this paper is:

\bt[Non-degenerate diffusion with jumps]\label{main1}
Let $\Gamma^{j,2}_{0,R}(g)$ be defined as in \eqref{LK1}. Suppose that {\bf (H$^\sigma$)} holds
and for any $T>0$,
$$
\Gamma^{0,2}_{0,R}(g)\in\mL^\infty(T),\ \lim_{\eps\to 0}\|\Gamma^{0,2}_{0,\eps}(g)\|_{\mL^\infty(T)}=0,
$$
and for some $p,q\in(2,\infty)$ with $\frac{d}{p}+\frac{2}{q}<1$,
$$
|\nabla\sigma|,\ b,\ \big(\Gamma^{1,2}_{0,R}(g)\big)^{1/2}\in \mL^{q}_p(T).
$$
Then for any initial value $X_0=x\in\mR^d$, SDE \eqref{sde1} admits a unique strong solution $X_t(x)$.
Moreover, for any $T>0$, there is a constant $c_T>0$ such that
for all $t\in(0,T]$, $x,y\in\mR^d$ and bounded measurable $\varphi$,
\begin{align}
\big|\mE \varphi(X_t(x))-\mE \varphi(X_t(y))\big|\leq \frac{c_T}{\sqrt{t}}\|\varphi\|_{\infty}|x-y|. \label{feller}
\end{align}
\et

Let us make some comments on the above result.

\br
If $g_t(x,z)=\bar\sigma_t(x)z$ with $\bar\sigma_t(x)\in\mL^\infty(T)$ and $\nabla \bar\sigma_t(x)\in\mL^q_p(T)$
in the above theorem, then the assumptions on $\Gamma^{j,2}_{0,R}(g)$ automatically hold.
In particular, if $\bar\sigma_t(x)=\bar\sigma(x)=|x|^\beta\mI$ for some $\beta\in(0,1)$, then one can check
$\nabla\bar\sigma\in L^p_{loc}(\mR^d)$ for any $p<d/(1-\beta)$.
\er
\br
It is noticed that in the estimate \eqref{feller}, we do not make any assumption about the large jump coefficient
since the large jump part is independent with small jump part
and
%it does not make any contribution about the  regularization effect.
has only finitely many jumps in any finite time interval.
%since \eqref{feller} holds before the first large jump occurs (see \cite{W-X-Zh}).
\er

In the above mixing case, the non-degenerate diffusion part plays a dominant role.
In the pure jump case, we need to use the regularization effect of the jump noise. For this, we assume
$\nu(\dif z)=|z|^{-d-\alpha}\dif z$ for some $\alpha\in(1,2)$, and $g$ satisfies that
\begin{enumerate}[{\bf (H$^g$)}]
\item $g_t(x,0)=0$ and there is a constant $c_1>0$ such that for all $t\geq 0,x, x',z,z'\in\mR^d$,
\begin{align}\label{CN1}
c^{-1}_1|z-z'|\leq |g_t(x,z)-g_t(x,z')|\leq c_1|z-z'|,
\end{align}
and for some $\beta\in(0,1)$ and $j=0,1$,
\begin{align*}
|\nabla^{j}_z g_t(x,z)-\nabla^{j}_z g_t(x',z)|\leq c_1|x-x'|^\beta(|z|+|z|^{1-j}),
\end{align*}
\begin{align*}
|\nabla_z g_t(x,z)-\nabla_z g_t(x,z')|\leq c_1|z-z'|^\beta.
\end{align*}
\end{enumerate}
Under \eqref{CN1}, the map $z\mapsto g_t(x,z)$ admits an inverse denoted by $g^{-1}_t(x,z)$.
The main point for us   is that the assumption {\bf (H$^g$)} is invariant under Zvonkin's transformation,
that is, $\tilde g$ defined by \eqref{HD7} still satisfies {\bf (H$^g$)}, see Proposition \ref{Pr63} below.
More importantly, by the change of variables, it allows us to write
\begin{align}
\sL^g_\nu u(x)&=\int_{\mR^d}\Big[u\big(x+z\big)-u(x)-z\cdot\nabla u(x)\Big]\frac{\kappa(t,x,z)}{|z|^{d+\alpha}}
\dif z+\bar b^g_t\cdot\nabla u(x),\label{HF3}
\end{align}
where $\bar b^g_t(x):=\int_{|z|\geq 1}g_t(x,z)\nu(\dif z)$ and
\begin{align}\label{ka}
\kappa(t,x,z):=\frac{|z|^{d+\alpha}\det(\nabla_z g^{-1}_t(x,z))}{|g^{-1}_t(x,z)|^{d+\alpha}}.
\end{align}
Our second well-posedness result is:

\bt[Multiplicative pure jump noise]\label{main3}
Suppose that $\sigma\equiv0$, $\nu(\dif z)=|z|^{-d-\alpha}\dif z$ for some $\alpha\in(1,2)$, and {\bf (H$^g$)} holds.
Moreover, we also suppose that
for some $\theta\in(1-\frac{\alpha}{2},1)$,
$p\in(\frac{2d}{\alpha}\vee 2,\infty)$ and $q\in(\frac{2\alpha}{\alpha+2(\theta-1)},\infty)$ with $\frac{d}{p}+\frac{\alpha}{q}<\frac{\alpha}{2}$,
$$
\big(\Gamma^{1,2}_{0,R}(g)\big)^{1/2},\ (\mI-\Delta)^{\theta/2}b\in \cap_{T>0}\mL^{q}_p(T).
$$
Then for each $X_0=x\in\mR^d$,  SDE \eqref{sde1} admits a unique strong solution $X_t(x)$. Moreover, $X_t(x)$ has a density $\rho(t,x,y)$,
 which enjoys the following estimates:
 \begin{enumerate}[(i)]
 \item (Two sides estimate) For any $T>0$, there are two constants $c_1,c_2>0$ such that for all $t\in(0,T)$ and $x,y\in\mR^d$,
 \begin{align}\label{GD5}
c_1t(t^{1/\alpha}+|x-y|)^{-d-\alpha}\leq \rho(t,x,y)\leq c_2t(t^{1/\alpha}+|x-y|)^{-d-\alpha}.
 \end{align}
 \item (Gradient estimate) For any $T>0$, there is a constant $c_3>0$ such that for all $t\in(0,T)$ and $x,y\in\mR^d$,
 \begin{align}\label{GD4}
|\nabla_x\rho(t,x,y)|\leq c_3t^{1-1/\alpha}(t^{1/\alpha}+|x-y|)^{-d-\alpha}.
 \end{align}
 \end{enumerate}
\et

%We would like to make the following important comments.
We would like to make the following comment.

\br
%{\it(1)}
If $g_t(x,z)=\bar\sigma_t(x)z$ with $\bar\sigma$ satisfying {\bf (H$^\sigma$)}, then  {\bf (H$^g$)} holds.
In this case, $\big(\Gamma^{1,2}_{0,R}(g)\big)^{1/2}=c|\nabla\bar\sigma|$ for some $c>0$.
Compared with the additive noise case considered in \cite{H-F,Pri,Pri2,Zh00}, we drop the boundness condition on the drift $b$, which is essentially used in their proof. For interesting examples of the drift $b$ which can be discontinuous, see \cite{Zh00}.
\er

Let $\chi:\mR^d\to[0,1]$ be a smooth function with $\chi(x)=0$ for $|x|\geq 2$ and $\chi(x)=1$ for $|x|\leq 1$.
For $m\in\mN$, define the cutoff function  $\chi_m$ by
\begin{align}\label{Cut}
\chi_m(x):=\chi(m^{-1} x).
\end{align}
Using suitable localization technique, we have: %the following local well-posedness result.

\bc[Local well-posedness]\label{local}
Suppose that for each $m\in\mN$,
\begin{align*}
\sigma^m_t(x):=\sigma_t(x\chi_m(x)),\ b^m_t(x):=b_t(x)\chi_m(x),\
g^m_t(x,z):=g_t(x\chi_m(x),z)
\end{align*}
satisfy the same assumptions as in Theorem \ref{main1} or Theorem \ref{main3}. Then SDE \eqref{sde1}
admits a unique strong solution $X_t$ up to the explosion time $\zeta$, that is,
$\lim_{t\uparrow\zeta}X_t=\infty$. \ec
\begin{proof}
By Theorem \ref{main1} or Theorem \ref{main3}, for each $m\in\mN$, there exist a unique global strong solution $X_t^m$ to SDE (\ref{sde1})
with coefficients $\sigma^m, g^m$ and $b^m$. For $m\geq k$, define
$$
\zeta_{m,k}:=\inf\{t\geq0: |X_t^m|\geq k\}\wedge m.
$$
By the uniqueness of the solution, we have
$$
\mP\Big(X_t^m=X_t^k,\,\forall t\in[0,\zeta_{m,k})\Big)=1,
$$
which implies that for $m\geq k$,
$$
\zeta_{k,k}\leq \zeta_{m,k}\leq \zeta_{m,m},\quad a.s..
$$
Hence, if we let $\zeta_k:=\zeta_{k,k}$, then $\zeta_k$ is an increasing sequence of $(\sF_t)$-stopping times and for $m\geq k$,
$$
\mP\Big(X_t^m=X_t^k,\,\forall t\in[0,\zeta_{k})\Big)=1.
$$
Now, for each $k\in\mN$, we can define $X_t:=X^k_t$ for $t<\zeta_k$
and $\zeta:=\lim_{k\rightarrow\infty}\zeta_k$. It is easy to see that $X_t$ is the unique solution of SDE (\ref{sde1})
up to the explosion time $\zeta$ and $\lim_{t\uparrow\zeta}X_t=\infty$ a.s.
\end{proof}

As for the non-explosion, under some Lyapunov conditions, we may show the existence of global solutions (for instance, see Lemma \ref{Le71} below).

\subsection{Ergodicity of SDEs with singular dissipative coefficients}
Below we turn to the study of the ergodicity of SDE \eqref{sde1}. We first recall some basic notions about the ergodicity.
Let $(P_t)_{t\geq 0}$ be a semigroup of bounded linear operators on Banach space $\cB_b(\mR^d)$, where $\cB_b(\mR^d)$ denotes
the space of all bounded Borel measurable functions. Let $\mu$ be a probability measure on Borel space $(\mR^d,\sB(\mR^d))$.
We use the following standard notation:
$$
\<\mu, \varphi\>:=\int_{\mR^d}\varphi(x)\mu(\dif x).
$$
\begin{enumerate} [$\bullet$]
\item $\mu$ is said to be an invariant probability measure (or stationary distribution) of $P_t$ if
$$
\<\mu, P_t\varphi\>=\<\mu,\varphi\>,\quad\forall t>0,\,\,\forall \varphi\in\cB_b(\mR^d).
$$
\item One says that $P_t$ is ergodic if $P_t$ admits a unique invariant probability measure $\mu$, which amounts to say that
\begin{align}\label{Erg}
\lim_{t\to\infty}\frac{1}{t}\int^t_0 P_sf(x)\dif s=\<\mu,f\>,\ \ f\in \cB_b(\mR^d).
\end{align}
\item  One says that $P_t$ has the strong Feller property if for all $\varphi\in\cB_b(\mR^d)$, $P_t\varphi\in C_b(\mR^d)$.
\item $P_t$ is said to be irreducible if for each open ball $B$ and $x\in\mR^d$, $P_t1_B(x)>0$.
\end{enumerate}

About the ergodicity, we have the following more precise classification (cf. \cite{Me-Tw} and \cite{Ha}).

\bd\label{de}
Let $V:\mR^d\rightarrow[1,\infty)$ be a measurable function and $\mu$ an invariant probability measure of $P_t$.
We say $P_t$ to be $V$-uniformly exponential ergodic
if there exist $c_0,\gamma>0$ such that for all $t\geq 0$ and $x\in\mR^d$,
$$
\sup_{\|\varphi\|_V\leq 1}\big|P_t\varphi(x)-\<\mu,\varphi\>\big|\leq c_0V(x)\e^{-\gamma t},
$$
where $\|\varphi\|_V:=\sup_{x\in\mR^d}\frac{|\varphi(x)|}{V(x)}<+\infty$.
If $V\equiv1$, then $P_t$ is said to be uniformly exponential ergodic, which is equivalent to
$$
\|P_t(x,\cdot)-\mu\|_{Var}\leq c_0\e^{-\gamma t},\quad \forall x\in\mR^d,
$$
where $\|\cdot\|_{Var}$ is the total variation of a signed measure, $P_t(x,\cdot)$ is the kernel of bounded linear operator $P_t$.
\ed

It is useful to observe that the above notions are invariant under homeomorphism transformation of the phase space.
More precisely, let $\Phi:\mR^d\to\mR^d$ be a homeomorphism. Define a new semigroup of bounded linear operators on $\cB_b(\mR^d)$ by
$$
P^\Phi_t\varphi(y):=[P_t(\varphi\circ\Phi)](\Phi^{-1}(y)),
$$
where $\Phi^{-1}$ is the inverse of $\Phi$. We have the following simple observations, which are direct by definition.
\bp\label{Pr4}
\begin{enumerate}[(i)]
\item $\mu$ is an invariant probability measure of $P_t$ if and only if $\mu\circ\Phi^{-1}$ is an invariant probability measure of $P^\Phi_t$.
\item $P_t$ has the strong Feller property if and only if $P^\Phi_t$ has the strong Feller property.
\item $P_t$ is irreducible  if and only if $P^\Phi_t$ is  irreducible.
\item $P_t$ is $V$-uniformly exponential ergodic  if and only if $P^\Phi_t$ is $V\circ\Phi^{-1}$-uniformly exponential ergodic.
\end{enumerate}
\ep

To study the ergodicity of SDE \eqref{sde1}, we shall assume that the coefficients are time-independent, i.e.,
\begin{align}
\dif X_t&=\sigma(X_t)\dif W_t+b(X_t)\dif t+\int_{|z|<R}g(X_{t-},z)\tilde N(\dif t,\dif z)+\int_{|z|\geq R}g(X_{t-},z)N(\dif t,\dif z). \label{sde4}
\end{align}
We have the following ergodicity result when the drift $b$ is  dissipative.
\bt\label{main22}
Suppose that the local condition in Corollary \ref{local} holds for $(\sigma,b,g)$, and
for some $r>-1$ and $\kappa_1,\kappa_2,\kappa_3>0$,
$$
2\<x,b(x)\>+\|\sigma(x)\|^2\leq-\kappa_1|x|^{2+r}+\kappa_2,\ \ |b(x)|\leq \kappa_3(1+|x|^{1+r}),
$$
and for any $\eps>0$ and $\lambda\geq R$, there is a $c_{\eps,\lambda}>0$ such that
$$
\Gamma^{0,2}_{0,\lambda}(g)(x)+\Gamma^{0,1}_{\lambda,\infty}(g)(x)\leq\eps|x|^{1+r}+c_{\eps,\lambda}.
$$
Then for each $X_0=x\in\mR^d$, SDE \eqref{sde4} has a unique global strong solution $X_t(x)$.
If we let $P_t\varphi(x):=\mE\varphi(X_t(x))$, then $P_t$ admits a unique invariant probability measure $\mu$. Moreover,
if $r=0$, then $P_t$ is $V$-uniformly exponential ergodic with $V(x)=1+|x|$; if $r>0$, then $P_t$ is uniformly exponential ergodic.
\et

In the above ergodicity result, the drift $b$ is assumed to be locally bounded. Next we show two new ergodicity results, which allow the drift
to be singular at infinity. To this end, we shall assume that
\begin{enumerate}[{\bf (H$^b$)}]
\item $b=b_1+b_2$, where $b_1$ is the singular part and for some $p>d$,
\begin{align*}
b_1\in L^p(\mR^d),
\end{align*}
and $b_2$ is the dissipative part which satisfies for some $\kappa_1,\kappa_2,\kappa_3>0$ and $r>-1$,
\begin{align}
\<x, b_2(x)\>\leq -\kappa_1|x|^{2+r}+\kappa_2\quad\text{and}\quad|b_2(x)|\leq \kappa_3(1+|x|^{1+r}).    \label{diss}
\end{align}
\end{enumerate}
We have the following ergodicity result.

\bt[Ergodicity for diffusion with jumps]\label{main2}
Suppose that {\bf (H$^\sigma$)}  and  {\bf (H$^b$)} hold and for the same $p$ in  {\bf (H$^b$)},
$$
|\nabla\sigma|,\ (\Gamma^{1,2}_{0,R}(g))^{1/2}\in L^p(\mR^d),
$$
and for any $\lambda\geq R$,
\begin{align}\label{JF}
\Gamma^{0,2}_{0,\lambda}(g),\Gamma^{0,1}_{\lambda,\infty}(g)\in L^\infty(\mR^d),\ \lim_{\eps\to 0}\|\Gamma^{0,2}_{0,\eps}(g)\|_{\infty}=0.
\end{align}
Then, the conclusions of Theorem \ref{main22} still hold and the invariant probability measure $\mu$ has a density $\rho\in L^q(\mR^d)$ with $q<d/(d-1)$.
\et

\br
If $b\in L^p_{loc}(\mR^d)$ for some $p>d$, and for some $m>0$, $b$ satisfies \eqref{diss} for all $|x|>m$, then {\bf (H$^b$)} holds.
In fact, it suffices to take $b_1=\chi_mb$ and $b_2=(1-\chi_m)b$. The typical function satisfying \eqref{diss} is given
by $b_2(x)=-x|x|^rc(x)$ with $0<c_0\leq c(x)\leq c_1$.
Moreover, if $g(x,z)=\bar\sigma(x)z$ with $\bar\sigma(x)\in L^\infty(\mR^d)$ and $\nabla \bar\sigma(x)\in L^p(\mR^d)$ for some $p>d$
and $\int_{|z|>1}|z|\nu(\dif z)<\infty$, then all the assumptions on $g$ in the above theorem hold.
\er
\br
Under some minimal assumptions, Bogachev, R\"ockner and Shaposhnikov \cite{B-R-S}
have already shown the absolute continuity of $\mu$ with respect to the Lebesgue measure. However, it seems that their results can not be used
to our singular case since it is not known whether $b\in L^1_{loc}(\mu)$.
\er
In the pure jump case, we assume $\nu(\dif z)=|z|^{-d-\alpha}\dif z$ for some $\alpha\in(1,2)$ and
\begin{enumerate}[{\bf ($\widetilde {\mathbf H}^b$)}]
\item $b=b_1+b_2$, where $b_2$ satisfies \eqref{diss}, and $b_1$ satisfies that for some $\theta\in(1-\alpha/2,1)$ and $p>2d/\alpha$,
\begin{align*}
(\mI-\Delta)^{\theta/2}b_1\in L^p(\mR^d).
\end{align*}
\end{enumerate}
We have:

\bt[Ergodicity for pure jump SDE]\label{main4}
Suppose that {\bf (H$^g$)}  and  {\bf ($\widetilde {\mathbf H}^b$)} hold and for the same $p$ in  {\bf ($\widetilde {\mathbf H}^b$)} and any $\lambda\geq R$,
$$
(\Gamma^{1,2}_{0,R}(g))^{1/2}\in L^p(\mR^d),\ \Gamma^{0,2}_{0,\lambda}(g),\Gamma^{0,1}_{\lambda,\infty}(g)\in L^\infty(\mR^d).
$$
Then the same conclusions of Theorem \ref{main22} hold and the invariant probability measure $\mu$ has a density $\rho\in L^q(\mR^d)$ with $q<d/(d-\alpha+1)$.
\et

\section{General stability and Zvonkin's transformation}

In this section we introduce two basic results: stability and Zvonkin's transformation for SDE \eqref{sde1}
under quite general assumptions.
First of all, we introduce the following important notion about Krylov's estimate.

\bd
Let $X=(X_t)_{t\geq 0}$ be an $\sF_t$-adapted process and $p,q\in[1,\infty)$. We say that Krylov's estimate holds for $X$ with index $p,q$, if for
all $T>0$, there is a constant $c_0>0$ such that for all $0\leq t_0\leq t_1\leq T$ and $f\in \mL^q_p(t_0,t_1)$,
\begin{align}\label{Kry}
\mE\left(\int^{t_1}_{t_0}f(s,X_s)\dif s\Big|\sF_{t_0}\right)\leq c_0\|f\|_{\mL^q_p(t_0,t_1)},
\end{align}
where $c_0$ will be called Krylov's constant of $X$.
\ed
\br\label{Rem}
Krylov's estimate \eqref{Kry} implies that for Lebesgue almost all $s$, the distribution of the random variable $X_s$ admits a density $\rho_s(y)$ with respect to
the Lebesgue measure so that
$$
\|\rho\|_{\mL^{q'}_{p'}(T)}\leq c_0,\ \ \tfrac{1}{p'}+\tfrac{1}{p}=1,\ \ \tfrac{1}{q'}+\tfrac{1}{q}=1,
$$
where $c_0$ is the Krylov constant of $X$.
\er
\br\label{Rem22}
Suppose that for some $p,q\in[1,\infty)$, Krylov's estimate holds for $X$ with index $p,q$. Then the Krylov estimate for $X$ also holds for any $p'\in[p,\infty)$ and
$q'$ with $p'-\frac{p'}{q'}=p-\frac{p}{q}$. In fact, by Remark \ref{Rem}, it  automatically holds that
$$
\mE\left(\int^{t_1}_{t_0}f(s,X_s)\dif s\Big|\sF_{t_0}\right)\leq \|f\|_{\mL^1_\infty(t_0,t_1)}.
$$
Notice that by the interpolation theorem (see \cite[Theorem 5.1.2]{Be-Lo}), we have
$$
(\mL^1_\infty(T), \mL^q_p(T))_{[\theta]}=\mL^{q'}_{p'}(T),
$$
where $\theta\in(0,1)$, $\frac{1}{q'}=1-\theta+\frac{\theta}{q}$ and $\frac{1}{p'}=\frac{\theta}{p}$, $(\cdot,\cdot)_{[\theta]}$ stands for the complex interpolation.
The desired Krylov estimate for $p'\in[p,\infty)$ and $q'$ with $p'-\frac{p'}{q'}=p-\frac{p}{q}$ follows by the interpolation theorem (see \cite{Zh5}).
\er

\br\label{Rem23}
Let $\{X^{(n)}, n\in\mN\}$ be a sequence of $\sF_t$-adapted processes. Suppose that $X^{(n)}$ satisfies
Krylov's estimate with the same index $p,q$ and Krylov's constant $c_0$. If $X^{(n)}_t$ converges to $X_t$ in probability for each $t$ as $n\to\infty$,
then $X$ still satisfies the Krylov estimate with the same index $p,q\in[1,\infty)$ and Krylov's constant $c_0$.
\er
The above definition about Krylov's estimate has the following useful consequence.

\bl[Khasminskii's type estimate]\label{Cor1}
Let $X=(X_t)_{t\geq 0}$ be an $\sF_t$-adapted process. Suppose that $X$ satisfies Krylov's estimate for some $p,q\in[1,\infty)$. Then for any $\lambda, T>0$,
$0\leq t_0\leq t_1\leq T$ and $f\in\mL^q_p(T)$,
\begin{align*}
\mE^{\sF_{t_0}}\exp\left(\lambda\int^{t_1}_{t_0}f(s,X_s)\dif s\right)\leq 2^{n},
\end{align*}
where $\mE^{\sF_{t_0}}(\cdot):=\mE(\cdot|{\sF_{t_0}})$,
and $n$ is chosen so that $\|f\|_{\mL^q_p((j-1)T/n, jT/n)}\leq\frac{1}{2\lambda c_0}$ for all $j=1,\cdots, n$, and $c_0$ is the Krylov constant of $X$.
\el
\begin{proof}
Without loss of generality, we assume $t_0=0$ and $t_1=T$.
For $\lambda>0$, let us choose $n$ large enough so that for $t_j=\frac{jT}{n}$,
\begin{align}\label{UT1}
\lambda c_0\|f\|_{\mL^p(t_j, t_{j+1})}\leq1/2, \ \ j=0,\cdots, n-1.
\end{align}
For $m\in\mN$, noticing that
$$
\left(\int^{t_{j+1}}_{t_j}g(s)\dif s\right)^m=m!\int\!\!\!\cdots\!\!\!\int_{\Delta^m}g(s_1)\cdots g(s_m)\dif s_1\cdots\dif s_m,
$$
where
$$
\Delta^m:=\Big\{(s_1,\cdots, s_m):  t_j\leq s_1\leq s_2\leq\cdots\leq s_m\leq t_{j+1}\Big\},
$$
by \eqref{Kry}, we have
\begin{align*}
\mE^{\sF_{t_j}}\left(\int^{t_{j+1}}_{t_j}f(s, X_s)\dif s\right)^m
&=m!\mE^{\sF_{t_j}}\left(\int\!\!\!\cdots\!\!\!\int_{\Delta^m}f(s_1,X_{s_1})\cdots f(s_m, X_{s_m})\dif s_1\cdots\dif s_m\right)\\
&=m!\mE^{\sF_{t_j}}\Bigg(\int\!\!\!\cdots\!\!\!\int_{\Delta^{m-1}}f(s_1, X_{s_1})\cdots f(s_{m-1}, X_{s_{m-1}})\\
&\quad\times\mE^{\sF_{s_{m-1}}}\left( \int^{t_{j+1}}_{s_{m-1}}f(s_m,X_{s_m})\dif s_m\right)\dif s_1\cdots\dif s_{m-1}\Bigg)\\
&\leq m!\mE^{\sF_{t_j}}\Bigg(\int\!\!\!\cdots\!\!\!\int_{\Delta^{m-1}}f(s_1, X_{s_1})\cdots f(s_{m-1}, X_{s_{m-1}})\\
&\quad \times c_0\|f\|_{\mL^q_p(t_j,t_{j+1})}\dif s_1\cdots\dif s_{m-1}\Bigg)\leq\cdots\leq m!(c_0\|f\|_{\mL^q_p(t_j, t_{j+1})})^m,
\end{align*}
which implies by \eqref{UT1} that
\begin{align*}
\mE^{\sF_{t_j}}\exp\left(\lambda\int^{t_{j+1}}_{t_j}f(s,X_s)\dif s\right)
&=\sum_{m}\frac{1}{m!}\mE^{\sF_{t_j}}\left(\lambda\int^{t_{j+1}}_{t_j}f(s,X_s)\dif s\right)^m\leq 2.
\end{align*}
Hence,
\begin{align*}
&\mE^{\sF_0}\exp\left(\lambda\int^T_0f(s,X_s)\dif s\right)=\mE^{\sF_0}\left(\prod_{j=0}^{n-1}\exp\left(\lambda\int^{t_{j+1}}_{t_j}f(s,X_s)\dif s\right)\right)\\
&=\mE^{\sF_0}\left(\prod_{j=0}^{n-2}\exp\left(\lambda\int^{t_{j+1}}_{t_j}f(s,X_s)\dif s\right)\mE^{\sF_{t_{n-1}}}
\exp\left(\lambda\int^{t_{n}}_{t_{n-1}}f(s,X_s)\dif s\right)\right)\\
&\leq 2\mE^{\sF_0}\left(\prod_{j=0}^{n-2}\exp\left(\lambda\int^{t_{j+1}}_{t_j}f(s,X_s)\dif s\right)\right)\leq\cdots\leq 2^{n}.
\end{align*}
The proof is complete.
\end{proof}
Let $f$ be a locally integrable function on $\mR^d$. The Hardy-Littlewood maximal function of $f$ is defined by
$$
\cM f(x):=\sup_{s>0}\fint_{B_s}|f(x+y)|\dif y,
$$
where $\fint_{B_s}:=\frac{1}{|B_s|}\int_{B_s}$ and
$|B_s|$ denotes the Lebesgue measure of ball $B_s:=\{x: |x|<s\}$.
The following well known results can be found in  \cite[Lemma 5.4]{Zh1} and \cite[p. 5, Theorem 1]{St}.

\bl\label{Le23}
(i) Let $\mB$ be a Banach space and $f:\mR^d\to\mB$ a  locally integrable function with $\nabla f\in L^1_{loc}(\mR^d;\mB^d)$. Then there is a Lebesgue zero set $E$ such that for all $x,y\notin E$,
\begin{align}
\|f(x)-f(y)\|_\mB\leq 2^d\int^{|x-y|}_0\!\!\!\fint_{B_s}\Big[\|\nabla f\|_\mB(x+w)+\|\nabla f\|_\mB(y+w)\Big]\dif w\dif s.   \label{w11}
\end{align}
In particular, if $\nabla f\in L^p(\mR^d;\mB^d)$ for some $p>d$, then
\begin{align}\label{w12}
\|f(x)-f(y)\|_\mB\leq c_{d,p}|x-y|^{1-d/p}\|\nabla f\|_{p}.
\end{align}
(ii) For $p\in(1,\infty]$, there is a constant $c_{d,p}>0$ such that for all $f\in L^p(\mR^d)$,
\begin{align}
\|\cM f\|_p\leq c_{d,p}\|f\|_p.   \label{mf}
\end{align}
\el
\begin{proof}
We prove \eqref{w12}. For $\alpha\in(\frac{d}{p}-\frac{1}{p},1-\frac{1}{p})$, by H\"older's inequality, we have
\begin{align*}
\int^{|x-y|}_0\!\!\!\fint_{B_s}\|\nabla f\|_\mB(x+w)\dif w\dif s
&\leq\left(\int^{|x-y|}_0\!\!\! s^{-\alpha p^*}\dif s\right)^{\frac{1}{p^*}}\left(\int^{|x-y|}_0\!\!\!s^{\alpha p}\fint_{B_s}\|\nabla f\|^p_\mB(x+w)\dif w\dif s\right)^{\frac{1}{p}}\\
&\lesssim\left(\int^{|x-y|}_0 s^{-\alpha p^*}\dif s\right)^{\frac{1}{p^*}}
\left(\int^{|x-y|}_0s^{\alpha p-d}\dif s\right)^{\frac{1}{p}}\|\nabla f\|_p\\
&\lesssim |x-y|^{1-d/p}\|\nabla f\|_p.
\end{align*}
Substituting this into \eqref{w11}, we obtain \eqref{w12}.
\end{proof}

The following lemma is an easy consequence of Lemmas \ref{Cor1} and \ref{Le23}.
\bl\label{Le25}
Let $X, Y$ be two $\sF_t$-adapted processes, which satisfy Krylov's estimate with the same index $p,q\in(1,\infty)$ and Krylov's constant $c_0$.
Let $f_t(x):\mR_+\times\mR^d\to\mR$ and $g_t(x,z):\mR_+\times\mR^d\times\mR^d\to\mR$ be two Borel functions.
Let $T,R>0$. Suppose that for some $r\geq 1$,
$$
h(t,x):=|\nabla f_t(x)|^r+\Gamma^{1,r}_{0,R}(g_t)(x)\in\mL^{q}_p(T).
$$
Then there exists an $\sF_t$-adapted process $\ell_t$ with the property
$$
\mE \e^{\lambda \int^T_0 \ell_s\dif s}\leq c\big(\lambda,d,r,c_0,\|h\|_{\mL^q_p(T)}\big)<\infty, \ \forall\lambda>0,
$$
such that for Lebesgue almost all $t\in[0,T]$,
\begin{align}\label{LL1}
|f_t( X_t)-f_t( Y_t)|^r+\int_{|z|<R}|g_t(X_t,z)-g_t(Y_t,z)|^{r}\nu(\dif z)\leq \ell_t|X_t-Y_t|^r\ \ a.s.
\end{align}
\el
\begin{proof}
First of all, by \eqref{w11} with $\mB=\mR$, we have
\begin{align}\label{EG1}
|f_t(x)-f_t(y)|^r\leq 2^{dr}|x-y|^r(\cM|\nabla f_t|(x)+\cM|\nabla f_t|(y))^r,
\end{align}
and by \eqref{w11} with $\mB=L^r(B_R;\nu)$,
\begin{align}
&\int_{|z|<R}|g_t(x,z)-g_t(y,z)|^{r}\nu(\dif z)\leq2^{dr}|x-y|^r\Big(\cM\|\nabla_x g_t\|_{\mB}(x)+\cM\|\nabla_x g_t\|_{\mB}(y)\Big)^r\no\\
&\qquad\qquad\qquad\leq2^{dr+r}|x-y|^r\Big(\cM(\Gamma^{1,r}_{0,R}(g_t))(x)+\cM(\Gamma^{1,r}_{0,R}(g_t))(y)\Big).\label{EG2}
\end{align}
Now let us define
\begin{align*}
\ell_t=2^{dr+r}\Big[\cM h(t,\cdot)(X_t)+\cM h(t,\cdot)(Y_t)\Big].
\end{align*}
It follows by \eqref{mf} and Lemma \ref{Cor1}
that  $\ell_t$ has the desired property. The desired estimate \eqref{LL1} follows by \eqref{EG1}, \eqref{EG2} and Remark \ref{Rem}.
\end{proof}

Next we show a stochastic Gronwall's inequality, which has independent interest.

\bl[Stochastic Gronwall's inequality]\label{im}
Let $\xi(t)$ and $\eta(t)$ be two nonnegative measurable $\sF_t$-adapted processes,
$A_t$ a continuous nondecreasing $\sF_t$-adapted process with $A_0=0$, $M_t$ a local martingale with $M_0=0$. Suppose that
$t\mapsto \xi_t$ is c\`adl\`ag and
\begin{align}\label{Gron}
\xi(t)\leq\xi(0)+\int^t_0\eta(s)\dif s+\int^t_0\xi(s)\dif A_s+M_t,\ \forall t\geq 0.
\end{align}
Then for any $0<q<p<1$ and $t\geq 0$, we have
\begin{align}
\big[\mE(\xi(t)^*)^{q}\big]^{1/q}\leq \Big(\tfrac{p}{p-q}\Big)^{1/q}\Big(\mE \e^{pA_t/(1-p)}\Big)^{(1-p)/p}\mE\Bigg(\xi(0)+\!\int_0^t\!\eta(s)\dif s\Bigg),  \label{gron}
\end{align}
where $\xi(t)^*:=\sup_{s\in[0,t]}\xi(s)$.
\el
\begin{proof}
Without loss of generality, we assume that the right hand side of \eqref{gron} is finite. Let
$\bar\xi(t)$ be the right hand side of \eqref{Gron} and $\bar A_t:=\int^t_0\xi(s)/\bar\xi(s)\dif A_s$. Then
$$
\xi(t)\leq {\bar\xi}(t)=\xi(0)+\int^t_0\eta(s)\dif s+\int^t_0\bar\xi(s)\dif \bar A_s+M_t.
$$
By It\^o's formula, one has
$$
\e^{-\bar A_t}\bar \xi(t)=\xi(0)+\int_0^t\!\e^{-\bar A_s}\eta(s)\dif s+\int_0^t\!\e^{-\bar A_s}\dif M_s.
$$
Using $\e^{-\bar A_s}\leq 1$ and suitable stopping time technique, we get for any stopping time $\tau$,
$$
\mE\Big(\e^{-\bar A_{t\wedge\tau}}\bar\xi(t\wedge\tau)\Big)\leq \mE\Bigg(\xi(0)+\!\int_0^{t\wedge\tau}\!\eta(s)\dif s\Bigg),
$$
which yields by H\"older's inequality, $\xi(t)\leq\bar\xi(t)$ and $\bar A_t\leq A_t$ that for any $p\in(0,1)$,
$$
\mE\xi(t\wedge\tau)^p\leq
\mE \bar\xi(t\wedge\tau)^p\leq \Big(\mE \e^{pA_{t}/(1-p)}\Big)^{1-p}\left[\mE\Bigg(\xi(0)+\!\int_0^t\!\eta(s)\dif s\Bigg)\right]^p.
$$
Now, for any $\lambda>0$, define a stopping time
$$
\tau_\lambda:=\inf\{s\geq 0: \xi(s)\geq\lambda\}.
$$
Since $\xi$ is c\`adl\`ag, we have $\xi_{\tau_\lambda}\geq\lambda$ and
\begin{align*}
\lambda^p\mP\left(\xi(t)^*>\lambda\right)&\leq\lambda^p\mP\left(\tau_\lambda\leq t\right)
\leq\mE \xi(t\wedge\tau_\lambda)^p\\
&\leq \Big(\mE \e^{p A_{t}/(1-p)}\Big)^{1-p}\left[\mE\Bigg(\xi(0)+\!\int_0^t\!\eta(s)\dif s\Bigg)\right]^p=:\delta,
\end{align*}
and for any $q\in(0,p)$,
\begin{align*}
\mE|\xi(t)^*|^{q}&=q\int^\infty_0\lambda^{q-1}\mP(\xi(t)^*>\lambda)\dif\lambda\\
&\leq q\int^\infty_0\lambda^{q-1}((\lambda^{-p}\delta)\wedge 1)\dif\lambda=p\delta^{q/p}/(p-q).
\end{align*}
The proof is complete.
\end{proof}
\br
In \cite{Sc}, Scheutzow proved \eqref{gron} for continuous martingales. His proof depends on
a martingale inequality of Burkholder, which does not hold for discontinuous martingale as pointed out by him.
Compared with the proof provided in \cite{Sc}, our proof is more elementary.
\er

The following stability result and Zvonkin's transformation will be our cornerstone.
\bt[Stability]\label{Sta}
For $i=1,2$, let $X^{(i)}_t$ satisfy the following SDE
$$
X^{(i)}_t=X^{(i)}_0+\!\!\int_0^t\sigma^{(i)}_s(X^{(i)}_s)\dif W_s+\!\!\int_0^tb^{(i)}_s(X^{(i)}_s)\dif s
+\!\!\int_0^t\!\!\!\int_{|z|<R}g^{(i)}_s(X^{(i)}_{s-},z)\tilde{N}(\dif s, \dif z),
$$
where $(\sigma^{(i)}, b^{(i)}, g^{(i)})$ are two families of measurable coefficients. Let $\ga\geq 1$.
Suppose that $X^{(1)}$ satisfies Krylov's estimate with index $p,q\in(1,\!\infty)$, and for all $T\!>0$,
there are $p_i\in[p,\infty]$ and $q_i=1/(1-(p-p/q)/p_i)$, $i=1,2,3,4$ such that
$$
\hbar:=\|\nabla \sigma^{(1)}\|^2_{\mL^{2q_1}_{2p_1}(T)}+\|\nabla b^{(1)}\|_{\mL^{q_2}_{p_2}(T)}
+\big\|\Gamma^{1,2}_{0,R}(g^{(1)})\big\|_{\mL^{q_3}_{p_3}(T)}+\big\|\Gamma^{1,2\ga}_{0,R}(g^{(1)})\big\|_{\mL^{q_4}_{p_4}(T)}<\infty,
$$
where $\Gamma^{j,\alpha}_{0,R}(g^{(1)})$ is defined by \eqref{LK1}.
Then for any $\theta\in(0,1)$ and $T>0$,
\begin{align*}
\left[\mE\left(\sup_{t\in[0,T]}|X^{(1)}_t-X^{(2)}_t|^{2\ga\theta}\right)\right]^{1/\theta}
\!\!\!\leq c_1\Bigg[\mE|X^{(1)}_0-X^{(2)}_0|^{2\ga}+\mE\left(\int^T_0\delta_s(X^{(2)}_s)\dif s\right)\Bigg],
\end{align*}
where $c_1$ only depends on $T,r,\theta,p,q,d,\hbar$ and the Krylov constant of $X^{(1)}$, and
\begin{align}\label{Del}
\begin{split}
\delta_s(x)&:=\|\sigma^{(1)}_s(x)-\sigma^{(2)}_s(x)\|^{2\ga}+|b^{(1)}_s(x)-b^{(2)}_s(x)|^{2\ga}\\
&\quad+\Gamma^{0,2\ga}_{0,R}\big(g^{(1)}_s-g^{(2)}_s\big)(x)+\Big(\Gamma^{0,2}_{0,R}\big(g^{(1)}_s-g^{(2)}_s\big)(x)\Big)^\ga.
\end{split}
\end{align}
\et
\begin{proof}
For simplicity of notations, we write $Z_t:=X^{(1)}_t-X^{(2)}_t$ and
\begin{align*}
&\Sigma_t:=\sigma^{(1)}_t(X^{(1)}_t)-\sigma^{(2)}_t(X^{(2)}_t),\ \ B_t:=b^{(1)}_t(X^{(1)}_t)-b^{(2)}_t(X^{(2)}_t),\\
&\qquad\qquad G_t(z):=g^{(1)}_t(X^{(1)}_{t},z)-g^{(2)}_t(X^{(2)}_{t},z).
\end{align*}
Since $X^{(1)}$ satisfies Krylov's estimate with index $p,q\in(1,\infty)$,
by the assumption, Remark \ref{Rem22} and Lemma \ref{Le25}, there exist $\sF_t$-adapted processes $\ell^{(j)}_t$ with
\begin{align}\label{YT1}
\mE\e^{\lambda\int^T_0\ell^{(j)}_s\dif s}\leq c(\lambda, \hbar)<\infty,\ \ \lambda>0,\ \ j=1,2,3,4,
\end{align}
such that
\begin{align}\label{UT3}\begin{split}
&\qquad|\Sigma_t|^2\leq \ell^{(1)}_t|Z_t|^2+2\|\sigma^{(1)}_t-\sigma^{(2)}_t\|^2(X^{(2)}_t),\\
&\qquad |B_t|\leq \ell^{(2)}_t|Z_t|+|b^{(1)}_t-b^{(2)}_t|(X^{(2)}_t),\\
&\int_{|z|<R}|G_t(z)|^{2}\nu(\dif z)\leq \ell^{(3)}_t|Z_t|^2+2\Gamma^{0,2}_{0,R}\big(g^{(1)}_t-g^{(2)}_t\big)(X^{(2)}_t),\\
&\int_{|z|<R}|G_t(z)|^{2\ga}\nu(\dif z)\leq \ell^{(3)}_t|Z_t|^{2r}+2^r\Gamma^{0,2\ga}_{0,R}\big(g^{(1)}_t-g^{(2)}_t\big)(X^{(2)}_t),
\end{split}
\end{align}
where $\ga\geq 1$. Now, by It\^o's formula, we have
\begin{align*}
\dif |Z_t|^{2\ga}&=\Big(\ga\|\Sigma_t\|^2|X_t|^{2(\ga-1)}+2\ga(\ga-1)|\Sigma_t X_t|^2|X_t|^{2(\ga-2)}
+2\ga\<B_t,Z_t\>|X_t|^{2(\ga-1)}\Big)\dif t\\
&\quad+\Bigg[\int_{|z|<R}\Big(|Z_t+G_t(z)|^{2\ga}-|Z_t|^{2\ga}-2\ga\<G_t(z), Z_t\>|Z_t|^{2(\ga-1)}\Big)\nu(\dif z)\Bigg]\dif t+\dif M_t,
\end{align*}
where $M_t$ is a local martingale.
Noticing that
$$
|x+y|^{2\ga}-|x|^{2\ga}-2\ga\<y,x\>|x|^{2(\ga-1)}\lesssim |y|^{2\ga}+|y|^2|x|^{2(\ga-1)},
$$
by \eqref{UT3} and Young's inequality, we get
\begin{align*}
\dif |Z_t|^{2\ga}&\lesssim |Z_t|^{2\ga}(\ell^{(1)}_t+\ell^{(2)}_t+\ell^{(3)}_t+\ell^{(4)}_t+1)\dif t+\delta_t(X^{(2)}_t)\dif t+\dif M_t,
\end{align*}
where $\delta_t(x)$ is defined by \eqref{Del}.
By Lemma \ref{im} and \eqref{YT1}, we obtain the desired estimate.
\end{proof}

The following proposition provides a way of transforming SDE \eqref{sde1} into a new SDE, which is called Zvonkin's transformation in the literature.
\bp[Zvonkin's transformation]\label{Pr313}
For each $t\geq 0$, let $\Phi_t(x)$ be a homeomorphism over $\mR^d$. Let $p,q\in(1,\infty)$.
Suppose that there exist a sequence of smooth functions $\Phi^n$ and a function $\bar b\in L^q_{loc}(\mR_+; L^p_{loc}(\mR^d))$ such that
for each $T>0$ and $(t,x)\in\mR_+\times\mR^d$,\ $m\in\mN$,
\begin{align}\label{UU9}
\sup_{n\in\mN}\|\nabla\Phi^n\|_{\mL^\infty(T)}<\infty,\ \ \lim_{n\to\infty}\Phi^n_t(x)=\Phi_t(x),\
\lim_{n\to\infty}\big\|\nabla(\Phi^n-\Phi)\chi_m\big\|_{\mL^q_p(T)}=0,
\end{align}
and
$$
\lim_{n\to\infty}\big\|((\p_s+\sL^\sigma_2+\sL^b_1+\sL^g_{\nu,R})\Phi^n-\bar b)\chi_m\big\|_{\mL^q_p(T)}=0,
$$
where $\chi_m$ is the cutoff function defined by \eqref{Cut}.
If $X$ solves SDE \eqref{sde1} and satisfies Krylov's estimate with the above index $p,q$,
then $Y_t\!:=\!\Phi_t(X_t)$ solves the following SDE:
\begin{align}\label{SDE4}
\dif Y_t=\tilde\sigma_t(Y_t)\dif W_t+\tilde b_t(Y_t)\dif t+\!\!\int_{|z|<R}\tilde g_t(Y_{t-},z)\tilde{N}(\dif t, \dif z)
+\!\!\int_{|z|\geq R}\!\!\tilde g_t(Y_{t-},z){N}(\dif t, \dif z),
\end{align}
where
\begin{align*}
\begin{split}
&\tilde\sigma_t(y):=\big(\nabla\Phi_t\cdot\sigma_t\big)\circ\Phi_t^{-1}(y),\quad \tilde b_t(y):=\bar b_t\big(\Phi_t^{-1}(y)\big),  \\
&\qquad\tilde g_t(y,z):=\Phi_t\Big(\Phi_t^{-1}(y)+g_t\big(\Phi_t^{-1}(y),z\big)\Big)-y.
\end{split}
\end{align*}
\ep
\begin{proof}
By It\^o's formula, we have
\begin{align*}
\Phi^n_t(X_t)&=\Phi^n_0(X_0)+\int_0^t\big(\nabla\Phi^n_s\cdot\sigma_s\big)(X_s)\dif W_s\\
&\quad+\!\!\int_0^t\!\!\!\int_{|z|<R}\Big[\Phi^n_s\big(X_{s-}+g_s\big(X_{s-},z\big)\big)-\Phi^n_s\big(X_{s-}\big)\Big]\tilde{N}(\dif s, \dif z)\\
&\quad+\!\!\int_0^t\!\!\!\int_{|z|\geq R}\Big[\Phi^n_s\big(X_{s-}+g_s\big(X_{s-},z\big)\big)-\Phi^n_s\big(X_{s-}\big)\Big]{N}(\dif s, \dif z)\\
&\quad+\int^t_0\Big((\p_s+\sL^\sigma_2+\sL^b_1+\sL^g_{\nu,R})\Phi^n_s\Big)(X_s)\dif s.
\end{align*}
Since $X$ satisfies Krylov's estimate with index $p,q$, by the assumptions and taking limits $n\to\infty$, we obtain SDE \eqref{SDE4}.
For example, for each $m\in\mN$, define
$$
\tau_m:=\inf\left\{t>0: |X_t|+\int^t_0\!\!\!\int_{|z|<R}|g_s(X_s,z)|^2\nu(\dif z)\dif s>m\right\}.
$$
By (\ref{UU9}) and the dominated convergence theorem, we have
\begin{align*}
&\mE\bigg|\int_0^{t\wedge\tau_m}\!\!\!\int_{|z|<R}\Big[\Phi^n_s\big(X_{s-}+g_s(X_{s-},z)\big)-\Phi^n_s(X_{s-})\\
&\qquad\quad\quad-\Phi_s\big(X_{s-}+g_s(X_{s-},z)\big)+\Phi_s(X_{s-})\Big]\tilde{N}(\dif s, \dif z)\bigg|^2\\
&=\mE\int_0^{t\wedge\tau_m}\!\!\!\int_{|z|<R}\Big|\Phi^n_s\big(X_s+g_s(X_{s},z)\big)-\Phi^n_s(X_s)\\
&\qquad\quad\qquad-\Phi_s\big(X_s+g_s(X_{s},z)\big)+\Phi_s(X_s)\Big|^2\nu(\dif z)\dif s\rightarrow0,\quad n\rightarrow\infty.
\end{align*}
Moreover, by Krylov's estimate for $X$ and the assumption, we also have
\begin{align*}
&\mE\left(\int^{t\wedge\tau_m}_0\big|(\p_s+\sL^\sigma_2+\sL^b_1+\sL^g_{\nu,R})\Phi^n_s-\bar b_s\big|(X_s)\dif s\right)\\
&\leq \mE\left(\int^{t}_0\big|((\p_s+\sL^\sigma_2+\sL^b_1+\sL^g_{\nu,R})\Phi^n_s-\bar b_s)\chi_m\big|(X_s)\dif s\right)\\
&\leq c\big\|((\p_s+\sL^\sigma_2+\sL^b_1+\sL^g_{\nu,R})\Phi^n-\bar b)\chi_m\big\|_{\mL^q_p(t)}\to 0,\ \ n\to\infty.
\end{align*}
The proof is complete since $\tau_m\to\infty$ as $m\to\infty$.
\end{proof}

To end this section, we recall a characterization about the Sobolev differentiability of random fields  in \cite[Theorem 1.1]{XZ},
which will be used to prove the Sobolev differentiability of the strong solution of the SDE with respect to the initial value.
For $p,q,r\in[1,\infty]$ and $T>0$, letting $L^r(T):=L^r([0,T])$, we define
$$
H^{1}_q\big(\mR^d;L^p(\Omega; L^r(T))\big):=\Big\{f(x,\omega,t): f,\,\nabla f\in L^{q}\big(\mR^d;L^p(\Omega; L^r(T))\big)\Big\},
$$
and
$$
\|f\|_{H^{1}_q(\mR^d;L^p(\Omega; L^r(T)))}:=\|f\|_{L^{q}(\mR^d;L^p(\Omega; L^r(T)))}+\|\nabla f\|_{L^{q}(\mR^d;L^p(\Omega;L^r(T)))}.
$$
\bt\label{th1}
Let $f\in L^{q}\big(\mR^d;L^p(\Omega; L^r(T))\big)$ for some $p\in(1,\infty)$ and $q,r\in(1,\infty]$.
Then $f\in H^1_q\big(\mR^d;L^p(\Omega; L^r(T))\big)$
if and only if there exists a nonnegative function $g\in L^q(\mR^d)$ such that for Lebesgue-almost all $x,y\in \mR^d$,
\begin{align}
\|f(x,\cdot)-f(y,\cdot)\|_{L^p(\Omega; L^r(T))}\leq |x-y|\big(g(x)+g(y)\big).\label{Mo23}
\end{align}
Moreover, if (\ref{Mo23}) holds, then for Lebesgue-almost all $x\in \mR^d$,
\begin{align*}
\|\p_if(x,\cdot)\|_{L^p(\Omega; L^r(T))}\leq 2g(x),\ \ i=1,\cdots,d,
\end{align*}
where $\p_i f$ is the weak partial derivative of $f$ with respect to the  $i$-th spacial variable.
\et

\section{A study of parabolic integral-differential equations}

This section is devoted to a careful study of the Kolmogorov backward equation associated to SDE \eqref{sde1}.

\subsection{Second order integral-differential equations}

First of all, we  introduce some Sobolev spaces and
notations for later use. For $(p,\alpha)\in[1,\infty]\times(0,2]\setminus\{\infty\}\times\{1\}$,  let $H^{\alpha}_p:=(\mI-\Delta)^{-\alpha/2}\big(L^p(\mR^d)\big)$
be the usual Bessel potential space with the norm
\begin{align}
\|f\|_{\alpha,p}:=\|(\mI-\Delta)^{\alpha/2}f\|_p\asymp\|f\|_p+\|\Delta^{\alpha/2}f\|_p,  \label{norm}
\end{align}
where $\|\cdot\|_p$ is the usual $L^p$-norm in $\mR^d$, and $(\mI-\Delta)^{\alpha/2}f$
and $\Delta^{\alpha/2}f$ are defined through the Fourier transformation
$$
(\mI-\Delta)^{\alpha/2}f:=\cF^{-1}\big((1+|\cdot|^2)^{\alpha/2}\cF f\big),\ \ \Delta^{\alpha/2}f:=\cF^{-1}\big(|\cdot|^{\alpha}\cF f\big).
$$
For $(p,\alpha)=(\infty,1)$, we define $H^1_\infty$ as the space  of Lipschitz functions with finite norm
$$
\|f\|_{1,\infty}:=\|f\|_\infty+\|\nabla f\|_\infty<\infty.
$$
Notice that for $n=1,2$ and $p\in(1,\infty)$, an equivalent norm in $H^n_p$ is given by% (\cite[p. 135, Theorem 3]{St})
$$
\|f\|_{n,p}=\|f\|_p+\|\nabla^nf\|_{p},
$$
and for $\alpha\in(0,2)$, up to a multiplying constant, an alternative expression of $\Delta^{\alpha/2}$ is given by
\begin{align}\label{Fra}
\Delta^{\alpha/2}f(x):={\rm p.v.}\int_{\mR^d}\frac{f(x+y)-f(x)}{|y|^{d+\alpha}}\dif y,
\end{align}
where p.v. stands for the Cauchy principal value.
We need the following Sobolev embedding: for $p\in[1,\infty]$ and $\alpha\in[0,2]$,
\begin{align}\label{Sob}
\left\{
\begin{aligned}
&H^\alpha_p\subset L^q,\ \ q\in[p,\tfrac{dp}{d-\alpha p}],\ \ &\alpha p<d;\\
&H^\alpha_p\subset H^{\alpha-d/p}_\infty\subset C_b^{\alpha-d/p},\ \ &\alpha p>d,
\end{aligned}
\right.
\end{align}
where $C^\beta_b$ is the usual H\"older space.
Moreover, for $\alpha\in[0,1]$ and $p\in(1,\infty]$, there is a constant $c=c(p,d,\alpha)>0$ such that for all $f\in H^\alpha_p$,
\begin{align}\label{UY1}
\|f(\cdot+y)-f(\cdot)\|_p\leq c(|y|^\alpha\wedge 1)\|f\|_{\alpha,p}.
\end{align}
The above facts are standard and can be found in \cite[Chapter 6]{Be-Lo} or \cite{Tri}.

The following lemma due to  \cite[Lemma 5]{M-P} strengthens the estimate \eqref{UY1}, which will play an important role in the following.
\bl\label{impo}
For $\alpha\in(0,2]$, write $y^{(\alpha)}:=y 1_{\alpha\in[1,2]}$. For any $p\in(\frac{d}{\alpha}\vee 1,\infty]$,
there is a constant $c=c(p,d,\alpha)>0$ such that for all $f\in H^\alpha_p$,
$$
\Big\|\sup_{y\not=0}|y|^{-\alpha}|f(x+y)-f(x)-y^{(\alpha)}\cdot\nabla f(x)|\Big\|_p\leq c\|f\|_{\alpha,p}.
$$
\el

Let $a(t,x):\mR_+\times\mR^d\to\mM^d_{sym}$ be a Borel measurable function, where $\mM^d_{sym}$ denotes the space of
all symmetric $d\times d$-matrices. Suppose that
\begin{enumerate}[{\bf (H$^a$)}]
\item There are constants $c_0\geq 1$ and $\beta\in(0,1)$ such that for all $(t,x)\in\mR_+\times\mR^d$,
$$
c_0^{-1}|\xi|^2\leq a^{ij}(t,x)\xi_i\xi_j\leq c_0|\xi|^2,\ \ \xi\in\mR^d,
$$
and
$$
\|a(t,x)-a(t,y)\|\leq c_0|x-y|^\beta.
$$
\end{enumerate}
We introduce the following second order partial differential operator:
\begin{align*}
\sL_2^au:=a^{ij}\p_{i}\p_ju,
\end{align*}
and for $0\leq S\leq T<\infty$, $\alpha\in(0,2]$ and $q,p\in(1,\infty]$, we introduce
$$
\mH^{\alpha,q}_p(S,T):=L^q\big([S,T];H^{\alpha}_p\big),\quad \mH^{\alpha,q}_p(T):=\mH^{\alpha,q}_p(0,T).
$$
Under {\bf (H$^a$)}, it is well known that $\sL_2^a$ admits a fundamental solution $\rho(s,x;t,y)$ so that (c.f. \cite{CHXZ})
$$
\p_s\rho(s,x;t,y)+\sL_2^a\rho(s,\cdot;t,y)(x)=0,\ \lim_{s\uparrow t}\rho(s,x;t,y)=\delta_{x-y}.
$$
Moreover, $\rho(s,x;t,y)$ enjoys the following upper and gradient estimates:
\begin{align}\label{GD2}
|\nabla^j_x\rho(s,\cdot;t,y)|(x)\leq c^{-1}_1(t-s)^{-(d+j)/2}\e^{-c_1|x-y|^2/(t-s)},\ j=0,1,2,
\end{align}
and fractional  derivative estimate: for any $\alpha\in(0,2)$,
\begin{align}\label{GD1}
|\Delta^{\alpha/2}\rho(s,\cdot;t,y)|(x)\leq c_2\big(|x-y|+(t-s)^{1/2}\big)^{-d-\alpha}.
\end{align}
The following $\mL^q_p$-estimate of the second order derivative was proven by Kim \cite{Ki} for $1<p\leq q<\infty$.
By duality, one in fact can drop the restriction $p\leq q$.

\bl
Let $\lambda, T\geq 0$ and  $p,q\in(1,\infty)$. Under {\bf (H$^a$)}, for any $f\in \mL^q_p(T)$, there exists a unique solution $u\in \mH^{2,q}_p(T)$ to the following backward PDE:
\begin{align*}
\p_tu+(\sL^a_2-\lambda) u=f,\quad u(T)=0.
\end{align*}
Moreover, for any $\vartheta\in(0,2)$ and $p'\in[p,\infty]$, $q'\in[q,\infty]$ satisfying
\begin{align}\label{Con1}
\frac{d}{p}+\frac{2}{q}<2-\vartheta+\frac{d}{p'}+\frac{2}{q'},
\end{align}
there exists a constant $c_1=c_1(d,p,q,\vartheta,p',q',T, c_0)>0$ such that for all $\lambda\geq 0$ and $S\in(0,T)$,
\begin{align}
(\lambda\vee 1)^{\frac{1}{2}(2-\vartheta+\frac{d}{p'}+\frac{2}{q'}-\frac{d}{p}-\frac{2}{q})}\|u\|_{\mH^{\vartheta,q'}_{p'}(S,T)}
+\|\nabla^2 u\|_{\mL^q_p(S,T)}\leq c_1\|f\|_{\mL^q_p(S,T)}.   \label{es1}
\end{align}
\el
\begin{proof}
It suffices to show the estimate \eqref{es1}. Without loss of generality, we assume $f\in C^\infty_c([0,T]\times\mR^d)$. First of  all, by \cite{Ki}, we have
$$
\|\nabla^2 u\|_{\mL^q_p(S,T)}\lesssim\|f\|_{\mL^q_p(S,T)}.
$$
On the other hand, by Duhamel's formula, we can write
$$
u(s,x)=\int^T_s\e^{-\lambda(t-s)}\left(\int_{\mR^d}\rho(s,x;t,y) f(t,y)\dif y\right)\dif t.
$$
Let $r=1/(1-1/p+1/p')$ and $\varrho_\vartheta(t,x):=\big(|x|+t^{1/2}\big)^{-d-\vartheta}$.
Suppose $(p',\vartheta)\not =(\infty,1)$. By \eqref{GD1} and Young's convolution inequality, we have
\begin{align*}
\|\Delta^{\vartheta/2}u(s)\|_{p'}&\leq\int^T_s\e^{-\lambda(t-s)}\left\|\int_{\mR^d}\Delta^{\vartheta/2}_x\rho(s,\cdot;t,y) f(t,y)\dif y\right\|_{p'}\dif t\\
&\lesssim\int^T_s\e^{-\lambda(t-s)}\left\|\int_{\mR^d}\varrho_\vartheta(t-s,\cdot-y) |f(t,y)|\dif y\right\|_{p'}\dif t\\
&\leq\int^T_s\e^{-\lambda(t-s)}\|\varrho_\vartheta(t-s,\cdot)\|_{r}\|f(t)\|_{p}\dif t\\
&\lesssim\int^T_s\e^{-\lambda(t-s)} (t-s)^{(d/r-\vartheta-d)/2}\|f(t)\|_{p}\dif t=(h_\lambda*\|f(\cdot)\|_{p})(s),
\end{align*}
where $h_\lambda(t):=\e^{-\lambda t} t^{(d/r-\vartheta-d)/2}1_{t>0}$. Hence, by Young's convolution inequality again,
$$
\|\Delta^{\vartheta/2}u\|_{\mL^{q'}_{p'}(S,T)}\lesssim\|h_\lambda\|_{L^{1/(1+1/q'-1/q)}(0,T-S)}\|f\|_{\mL^q_p(S,T)}
\lesssim(\lambda\vee 1)^{\frac{1}{\alpha}(\vartheta-2-\frac{d}{p'}-\frac{2}{q'}+\frac{d}{p}+\frac{2}{q})}\|f\|_{\mL^q_p(S,T)}.
$$
For $(p',\vartheta)=(\infty,1)$, by the gradient estimate \eqref{GD2}, we still have
$$
\|\nabla u\|_{\mL^{q'}_{\infty}(S,T)}(\lambda\vee 1)^{\frac{1}{2}(-1-\frac{2}{q'}+\frac{d}{p}+\frac{2}{q})}\|f\|_{\mL^q_p(S,T)}.
$$
Moreover, using the upper bound estimate of the heat kernel, we also have
$$
\|u\|_{\mL^{q'}_{p'}(S,T)}\lesssim (\lambda\vee 1)^{-1-\frac{1}{2}(\frac{d}{p'}+\frac{\alpha}{q'}-\frac{d}{p}-\frac{\alpha}{q})}\|f\|_{\mL^q_p(S,T)}.
$$
Combining the above calculations, we get \eqref{es1}.
\end{proof}
For $\lambda\geq 0$ and $R,T>0$, we consider the following backward second order parabolic integral-differential equation:
\begin{align}
\p_tu+(\sL^a_2-\lambda)u+\sL^b_1u+\sL^g_{\nu, R} u=f,\quad u(T)=0,\label{pide55}
\end{align}
where $\sL^g_{\nu,R}$ is defined in \eqref{UY7}.
We now prove the following solvability to the above equation.
\bt\label{pde}
Let $p\in(d/2\vee 1,\infty), q\in(1,\infty)$ and $T>0$.
Let $\Gamma^{0,2}_{0,R}(g)$ be defined by \eqref{LK1}. Assume that {\bf (H$^a$)} holds and 
\begin{enumerate}[(i)]
\item for some $p_1\in[p,\infty]$ and $q_1\in[q,\infty]$ with $\frac{d}{p_1}+\frac{2}{q_1}<1$, $b\in\mL^{q_1}_{p_1}(T)$;
\item $\Gamma^{0,2}_{0,R}(g)\in\mL^\infty(T)$ and $\lim_{\eps\to 0}\|\Gamma^{0,2}_{0,\eps}(g)\|_{\mL^\infty(T)}=0$.
\end{enumerate}
Then for some $\lambda_0>0$ depending on $\|b\|_{\mL^{q_1}_{p_1}(T)}$ and $\|\Gamma^{0,2}_{0,R}(g)\|_{\mL^\infty(T)}$,
and for all $\lambda\geq\lambda_0$ and $f\in \mL^q_p(T)$, 
there exists a unique solution $u\in \mH^{2,q}_p(T)$ to the equation (\ref{pide55}). Moreover,
 in this case the estimate \eqref{es1} still holds and $\p_t u\in\mL^q_p(T)$.
\et

\begin{proof}
By the standard continuity method, it suffices to show
the apriori estimate \eqref{es1} for equation \eqref{pide55} under the assumptions in the theorem. First of all, 
for any $\vartheta\in(0,2)$ and $p'\in[p,\infty],q'\in[q,\infty]$ satisfying \eqref{Con1}, by \eqref{es1} we have
\begin{align}\label{NM1}
\begin{split}
\lambda^{\frac{1}{2}(2-\vartheta+\frac{d}{p'}+\frac{2}{q'}-\frac{d}{p}-\frac{2}{q})}\|u\|_{\mH^{\vartheta,q'}_{p'}(S,T)}
+\|\nabla^2 u\|_{\mL^q_p(S,T)}\leq c_1\big\|f+\sL_1^bu+\sL_{\nu,R}^gu\big\|_{\mL^q_p(S,T)}.
\end{split}
\end{align}
Below, for simplicity of notation, we drop the time variable $t$.
Recalling the definitions of $\sL^g_{\nu,R}u$ and $\Gamma^{0,\alpha}_{\eps, R}(g)$ (see \eqref{LK1}), we have for any $\eps\in(0,R)$,
\begin{align*}
|\sL^g_{\nu,\eps}u(x)|&\leq\int_{|z|\leq\eps}\Big|u\big(x+g(x,z)\big)-u(x)-g(x,z)\cdot\nabla u(x)\Big|\nu(\dif z)\no\\
&\leq \sup_{y\not =0}|y|^{-2}|u(x+y)-u(y)-y\cdot\nabla u(x)|\,|\Gamma^{0,2}_{0,\eps}(g)(x)|,
\end{align*}
and for $\alpha\in(d/p\vee 1, 2)$,
\begin{align*}
&|\sL^g_{\nu,R}u(x)-\sL^g_{\nu,\eps}u(x)|\leq\int_{\eps<|z|<R}\Big|u\big(x+g(x,z)\big)-u(x)-g(x,z)\cdot\nabla u(x)\Big|\nu(\dif z)\no\\
&\qquad\qquad\leq \sup_{y\not =0}|y|^{-\alpha}|u(x+y)-u(y)-y\cdot\nabla u(x)\,|\Gamma^{0,\alpha}_{\eps,R}(g)(x)|.
\end{align*}
Thus, thanks to $p>d/\alpha \vee1$, by Lemma \ref{impo}, we obtain that for any $\eps\in(0,R)$,
\begin{align}\label{BB2}
\begin{split}
\|\sL_{\nu,R}^gu\big\|_{\mL^q_p(S,T)}&\leq \|\sL_{\nu,\eps}^gu\big\|_{\mL^q_p(S,T)}+\|\sL_{\nu,R}^gu-\sL_{\nu,\eps}^gu\big\|_{\mL^q_p(S,T)}\\
&\lesssim\|\Gamma^{0,2}_{0,\eps}(g)\|_{\mL^\infty(T)}\|u\|_{\mH^{2,q}_{p}(S,T)}+ \|\Gamma^{0,\alpha}_{\eps,R}(g)\|_{\mL^\infty(T)}\|u\|_{\mH^{\alpha,q}_{p}(S,T)}.
\end{split}
\end{align}
On the other hand, letting $q_2:=qq_1/(q_1-q)$ and $p_2:=pp_1/(p_1-p)$, by H\"older's inequality, we have
\begin{align}\label{BB1}
\|\sL_1^bu\|_{\mL^q_p(S,T)}\leq \|b\|_{\mL^{q_1}_{p_1}(S,T)}\|u\|_{\mH^{1,q_2}_{p_2}(S,T)}.
\end{align}
Now by \eqref{NM1}, \eqref{BB2} and \eqref{BB1}, there are $c_2, c_3>0$ such that for all $\eps\in(0,R)$,
\begin{align*}
&\lambda^{\frac{1}{2}(1+\frac{d}{p_2}+\frac{2}{q_2}-\frac{d}{p}-\frac{2}{q})}\|u\|_{\mH^{1,q_2}_{p_2}(S,T)}
+\lambda^{1-\frac{\alpha}{2}}\|u\|_{\mH^{\alpha,q}_{p}(S,T)}+\|\nabla^2 u\|_{\mL^q_p(S,T)}\\
&\leq c_2\Big(\|\Gamma^{0,2}_{0,\eps}(g)\|_{\mL^\infty(T)}\|u\|_{\mH^{2,q}_{p}(S,T)}
+\|\Gamma^{0,\alpha}_{\eps,R}(g)\|_{\mL^\infty(T)}\|u\|_{\mH^{\alpha,q}_{p}(S,T)}\Big)\\
&\qquad+c_3\Big(\|b\|_{\mL^{q_1}_{p_1}(T)}\|u\|_{\mH^{1,q_2}_{p_2}(S,T)}+\|f\|_{\mL^q_p(S,T)}\Big),
\end{align*}
which implies that for $\eps$ small enough and some $\lambda_0$ large enough and all $\lambda\geq\lambda_0$,
$$
\|u\|_{\mH^{1,q_2}_{p_2}(S,T)}+\|u\|_{\mH^{\alpha,q}_{p}(S,T)}+\|\nabla^2 u\|_{\mL^q_p(S,T)}\lesssim\|f\|_{\mL^q_p(S,T)}.
$$
Here we have used that $\lim_{\eps\to 0}\|\Gamma^{0,2}_{0,\eps}(g)\|_{\mL^\infty(T)}=0$ and
$$
\|\Gamma^{0,\alpha}_{\eps,R}\|_{\mL^\infty(T)}\leq \|\Gamma^{0,2}_{\eps,R}\|^{\alpha/2}_{\mL^\infty(T)}\nu(\{z: \eps<|z|<R\})^{1-\frac{\alpha}{2}}.
$$
Substituting this estimate into \eqref{NM1}, \eqref{BB2} and \eqref{BB1}, we get the estimate \eqref{es1}.
\end{proof}
\subsection{Non-local parabolic equations}
In this subsection we assume $\alpha\in(1,2)$ and  consider the following nonlocal operator
\begin{align}\label{49}
\sL^\kappa_\alpha f(x):=\int_{\mR^d} [f(x+z)-f(x)-z\cdot\nabla f(x)]\frac{\kappa(t,x,z)\dif z}{|z|^{d+\alpha}},
\end{align}
where $\kappa(t,x,z):\mR_+\times\mR^d\times\mR^d\to\mR$ satisfies that for some $\kappa_0>1$ and $\beta,\beta'\in(0,1]$,
\begin{align}\label{Con2}
\kappa^{-1}_0\leq \kappa(t,x,z)\leq\kappa_0,\ \ |\kappa(t,x,z)-\kappa(t,y,z)|\leq \kappa_0|x-y|^\beta(|z|^{\beta'}+1).
\end{align}
For $\beta\in[0,1]$ and $\gamma\in\mR$, let
$$
\varrho^\beta_\gamma(t,x):=t^{\gamma/\alpha}(|x|^\beta\wedge 1)(|x|+t^{1/\alpha})^{-d-\alpha}.
$$
It is easy to see that for any $p\geq 1$, there is a $c>0$ such that
\begin{align}\label{YT2}
\|\varrho^\beta_\gamma(t)\|_p\leq c t^{(\beta+\gamma-d-\alpha)/\alpha+d/(\alpha p)},\ \ t>0.
\end{align}

First of all, we introduce the following
Kato's class as in \cite{Wa-Zh}.
\bd
A Borel measurable function $f: [0,\infty)\times\mR^d\to\mR$ is said to be in Kato's class $\mK^\alpha_d$ if
$$
\lim_{\eps\downarrow 0}K^\alpha_f(\eps)=0,
$$
where
\begin{align}\label{Kf}
K^\alpha_f(\eps):=\eps  \sup_{(t,x)\in [0,\infty)\times\mR^d} \int^\eps_0\!\!\!\int_{\mR^d}
 \frac{\varrho^0_\alpha(s,x-y)|f(t\pm s, y)|}{s^{1/\alpha}(\eps-s)^{1/\alpha}}\dif y\dif s.
\end{align}
Here we have extended $f$ to $\mR$ by setting $f(t,\cdot)=0$ for $t<0$.
\ed
\br
For any $p,q\in[1,\infty]$, by H\"older's inequality, one sees that $\mL^q_p(T)\subset\mK^\alpha_d$ if $\frac{d}{p}+\frac{\alpha}{q}<\alpha-1$.
\er
We need the following result proved in \cite{Ch-Zh} and \cite{Ch-Zh0}.
\bt\label{heat}
Let $\alpha\in(1,2)$ and $b\in\mK^\alpha_d$. Under \eqref{Con2},
there is a unique continuous function $\rho_{\kappa,b}(s,x;t,y)$  satisfying
\begin{enumerate}[(i)]
\item (Two-sides estimate) There is a constant $c_0>1$ such that for all $s<t$ and $x,y\in\mR^d$,
\begin{align*}
c^{-1}_0\varrho^0_\alpha(t-s,x-y)\leq \rho_{\kappa,b}(s,x;t,y)\leq c_0 \varrho^0_\alpha(t-s,x-y).
\end{align*}
\item (Conservativeness)
For all $0\leq s<t<\infty$ and $x,y\in\mR^d$, it holds that
\begin{align*}
\int_{\mR^d}\rho_{\kappa,b} (s,x;t,y)\dif y=1.
\end{align*}
\item (C-K equation) For all $0\leq s<r<t<\infty$ and $x,y\in\mR^d$,
the following Chapman-Kolmogorov equation holds:
\begin{align*}
\int_{\mR^d}\rho_{\kappa,b}(s,x;r,z)\rho_{\kappa,b}(r,z;t,y)\dif z=\rho_{\kappa,b}(s,x;t,y).
\end{align*}
 \item (Generator) For any $f\in C_b^2(\mR^d)$, we have
      \begin{align*}
        P^{\kappa,b}_{s,t}f(x)-f(x)=\int^t_s\!P^{\kappa,b}_{s,r}(\sL^\kappa_\alpha+\sL^b_1)f(x)\dif r=\int^t_s\!(\sL^\kappa_\alpha+\sL^b_1)P^{\kappa,b}_{r,t}f(x)\dif r,
      \end{align*}
    where  $P^{\kappa,b}_{s,t}f(x):=\int_{\mR^d}\rho_{\kappa,b}(s,x;t,y)f(y)\dif y$.
  \item (Continuity) For any bounded and uniformly continuous function $f(x)$, we have
      \begin{align*}
        \lim_{|t-s|\to 0}\|P^{\kappa,b}_{s,t}f-f\|_\infty=0.
      \end{align*}
      \item (Gradient estimate) There is a constant $c_1>0$ such that  for all $0\leq s<t\leq T$ and $x,y\in\mR^d$,
\begin{align}\label{Grad1}
|\nabla_x \rho_{\kappa,b}(s,x;t,y)|\leq c_1 (s-t)^{-1/\alpha}\rho_{\kappa,b}(s,x;t,y).
\end{align}
\item (Fractional derivative estimate) For any $\theta\in[0,\alpha)$, there is a constant $c_3>0$ such that  for all $0\leq s<t\leq T$ and $x,y\in\mR^d$,
\begin{align}\label{Grad}
|\Delta^{\theta/2}_x\rho_{\kappa,b}(s,x;t,y)|\leq c_3\varrho^0_{\alpha-\theta}(t-s,x-y),
\end{align}
and for some $\eps\in(0,2-\alpha)$ and any $\theta\in[0,\alpha+\eps)$,
\begin{align}\label{Grad0}
|\Delta^{\theta/2}_x\rho_{\kappa,0}(s,x;t,y)|\leq c_4\varrho^0_{\alpha-\theta}(t-s,x-y).
\end{align}
\end{enumerate}
\et

As in the proof of \eqref{es1}, by Young's convolution inequality, the following result is an easy consequence of \eqref{YT2}, \eqref{Grad} and \eqref{Grad1}.
\bt\label{Th62}
Assume \eqref{Con2} and $b\in\mK^\alpha_d$.
Let $p,q\in(1,\infty)$ and  $p'\in[p,\infty]$, $q'\in[q,\infty]$ and $\vartheta\in[0, \alpha)$ with
$$
\frac{d}{p}+\frac{\alpha}{q}<\alpha-\vartheta+\frac{d}{p'}+\frac{\alpha}{q'}.
$$
For any $T>0$,
there is a constant $c>0$ only depending on $T, d,\alpha,\kappa_0,\beta,\gamma, p,q,p',q',\vartheta$ and the function $K^\alpha_b(\eps)$ defined by \eqref{Kf}
such that for all $\lambda\geq 0$ and $S\in(0,T)$,
\begin{align}\label{NR1}
(\lambda\vee 1)^{\frac{1}{\alpha}(\alpha-\vartheta+\frac{d}{p'}+\frac{\alpha}{q'}-\frac{d}{p}-\frac{\alpha}{q})}\|u\|_{\mH^{\vartheta,q'}_{p'}(S,T)}
\leq c\|f\|_{\mL^q_p(S,T)},
\end{align}
where $u(s,x):=\int^T_s\e^{-\lambda(t-s)}P^{\kappa,b}_{s,t}f(t,x)\dif t$.
\et

In \eqref{NR1}, the regularity of $u$ is at most $\vartheta$-order with $\vartheta<\alpha$.
In order to obtain higher regularity, we have to make further regularity assumption on $b$.
We first prepare two lemmas.
\bl
Assume \eqref{Con2}. Let $\theta\in[0,1)$ and $\gamma\in(0,\alpha)$ with $\gamma+\theta<\alpha+\eps$, where $\eps$ is the same as in (vii) of Theorem \ref{heat}.
For any $p\in(1,\infty]$ and $T>0$,
there is a constant $c>0$ such that for all $f\in H^\theta_p$ and $0\leq s<t\leq T$,
\begin{align}\label{UT4}
\|P^{\kappa,0}_{s,t}f\|_{\gamma+\theta,p}\leq c(t-s)^{-\gamma/\alpha}\|f\|_{\theta,p}.
\end{align}
\el
\begin{proof}
Since
$$
\int_{\mR^d}\Delta^{(\gamma+\theta)/2}_x\rho_{\kappa,0}(s,x;t,y) \dif y=\Delta^{(\gamma+\theta)/2}_x 1=0,
$$
by definition and \eqref{Grad0}, we have
\begin{align*}
|\Delta^{(\gamma+\theta)/2}_xP^{\kappa,0}_{s,t}f(x)|&=\left|\int_{\mR^d}\Delta^{(\gamma+\theta)/2}_x\rho_{\kappa,0}(s,x;t,y) (f(y)-f(x))\dif y\right|\\
&\lesssim\int_{\mR^d}\varrho^0_{\alpha-\gamma-\theta}(t-s,x-y)|f(y)-f(x)|\dif y\\
&=\int_{\mR^d}\varrho^0_{\alpha-\gamma-\theta}(t-s,y)|f(x-y)-f(x)|\dif y.
\end{align*}
By Minkovskii's inequality, \eqref{UY1} and \eqref{YT2}, we get
\begin{align*}
\|\Delta^{(\gamma+\theta)/2}P^{\kappa,0}_{s,t}f\|_p&\lesssim\int_{\mR^d}\varrho^0_{\alpha-\gamma-\theta}(t-s,y)\|f(\cdot-y)-f(\cdot)\|_p\dif y\\
&\lesssim\|f\|_{\theta,p}\int_{\mR^d}\varrho^\theta_{\alpha-\gamma-\theta}(t-s,y)\dif y
\lesssim\|f\|_{\theta,p} (t-s)^{-\gamma/\alpha}.
\end{align*}
The proof is complete.
\end{proof}

The following lemma can be regarded as an extension of H\"older's inequality to $H^\alpha_p$.
\bl\label{Le61}
For any $\alpha,\gamma_1,\gamma_2\in[0,1)$ and $p,p_1,p_2\in[1,\infty]$ with
$$
\frac{1}{p_i}<\frac{1}{p}+\frac{\gamma_i}{d},\quad \frac{\gamma_i}{d}\leq\frac{1}{p_1}+\frac{1}{p_2}-\frac{1}{p}<\frac{\gamma_1+\gamma_2+\alpha}{d}, \ \ i=1,2,
$$
there is a constant $c=c(p_i,\gamma_i,p,\alpha,d)>0$ such that
$$
\|fg\|_{\alpha,p}\leq c\|f\|_{\alpha+\gamma_1,p_1}\|g\|_{\alpha+\gamma_2,p_2}.
$$
\el
\begin{proof}
Let $p_0:=dp_1/(d-p_1\gamma_1)\geq p$. First of all, by H\"older's inequality and Sobolev's embedding \eqref{Sob}, we have
\begin{align}\label{GD8}
\|fg\|_p\leq\|f\|_{p_0}\|g\|_{pp_0/(p_0-p)}\lesssim\|f\|_{\gamma_1,p_1}\|g\|_{\alpha+\gamma_2,p_2}.
\end{align}
Notice that by \eqref{Fra},
$$
\Delta^{\alpha/2}(fg)=\int_{\mR^d}\frac{(f(\cdot+y)-f(\cdot))(g(\cdot+y)-g(\cdot))}{|y|^{d+\alpha}}\dif y+(\Delta^{\alpha/2}f)g+f\Delta^{\alpha/2}g.
$$
Hence,
\begin{align}\label{U0}
\begin{split}
\|\Delta^{\alpha/2}(fg)\|_p&\leq\int_{\mR^d}\frac{\|(f(\cdot+y)-f(\cdot))(g(\cdot+y)-g(\cdot))\|_p}{|y|^{d+\alpha}}\dif y\\
&\qquad+\|(\Delta^{\alpha/2}f)g\|_p+\|f\Delta^{\alpha/2}g\|_p.
\end{split}
\end{align}
As above, by H\"older's inequality and Sobolev's embedding, we have
\begin{align}\label{U1}
\|(\Delta^{\alpha/2}f)g\|_p\leq\|f\|_{\alpha,p_0}\|g\|_{p_0p/(p_0-p)}\lesssim\|f\|_{\alpha+\gamma_1,p_1}\|g\|_{\alpha+\gamma_2,p_2},
\end{align}
and by symmetry,
\begin{align*}
\|f(\Delta^{\alpha/2}g)\|_p\lesssim\|f\|_{\alpha+\gamma_1,p_1}\|g\|_{\alpha+\gamma_2,p_2}.
\end{align*}
Moreover, for $\eps\in(0,\gamma+\alpha-\frac{d}{p_1}-\frac{d}{p_2}+\frac{d}{p})$, by H\"older's inequality, Sobolev's embedding
and \eqref{UY1}, we have
\begin{align}\label{U3}
\begin{split}
&\|(f(\cdot+y)-f(\cdot))(g(\cdot+y)-g(\cdot))\|_p\\
&\leq \|f(\cdot+y)-f(\cdot)\|_{p_0}\|g(\cdot+y)-g(\cdot)\|_{pp_0/(p_0-p)}\\
&\lesssim \|f(\cdot+y)-f(\cdot)\|_{\gamma_1,p_1}\|g(\cdot+y)-g(\cdot)\|_{\alpha+\gamma_2-\eps,p_2}\\
&\lesssim \Big((|y|^\alpha\|f\|_{\alpha+\gamma_1,p_1})\wedge (2\|f\|_{\gamma_1,p_1})\Big)
\Big((|y|^\eps\|g\|_{\alpha+\gamma_2,p_2})\wedge(2\|g\|_{\alpha+\gamma-\eps,p_2})\Big)\\
&\lesssim (|y|^{\alpha+\eps}\wedge 1)\|f\|_{\alpha+\gamma_1,p_1}\|g\|_{\alpha+\gamma_2,p_2}.
\end{split}
\end{align}
Substituting \eqref{U1}-\eqref{U3} into \eqref{U0}, we obtain
$$
\|\Delta^{\alpha/2}(fg)\|_p\lesssim\|f\|_{\alpha+\gamma_1,p_1}\|g\|_{\alpha+\gamma_2,p_2},
$$
which together with \eqref{GD8} yields the desired estimate.
\end{proof}
We now study the following nonlocal parabolic equation:
\begin{align}\label{UY2}
\p_t u+\sL^\kappa_\alpha u-\lambda u+\sL^b_1 u=f,\ \ u(T)=0.
\end{align}
By Duhamel's formula, we shall consider the following mild form:
\begin{align}\label{In}
u(s,x)=\int^T_s\e^{-\lambda(t-s)}P^{\kappa,0}_{s,t}(\sL^b_1 u+f)(t,x)\dif t.
\end{align}
\br\label{Rem44}
It should be noticed that if $u$ is regular enough, saying, $u\in \mH^{\gamma,\infty}_{\infty}(T)$ for some $\gamma>\alpha$, then $u$ solves \eqref{In}
if and only if $u$ solves \eqref{UY2} because $u$ is in the domain of $\sL^\kappa_\alpha$ and $\sL^b_1$.
\er

We have the following regularity estimate.
\bt\label{Th63}
Let $\theta\in[0,1)$ and $p,q\in[1,\infty]$ with
$\frac{d}{p}+\frac{\alpha}{q}<\alpha+\theta-1$ and  $b\in\mH^{\theta,q}_{p}(T)$.
Under \eqref{Con2}, for any $\gamma\in(0, \alpha-\frac{\alpha}{q})$ with $\gamma+\theta<\alpha+\eps$, where $\eps$ is the same as in (vii) of Theorem \ref{heat},
there exists a constant
$c>0$ such that for all $\lambda\geq 0$ and $f\in\mH^{\theta,q}_p(S,T)$,
there is a unique solution $u\in\mH^{\gamma+\theta,q}_p(S,T)$ to \eqref{In} so that
\begin{align}\label{GR3}
(\lambda\vee 1)^{1-\frac{1}{q}-\frac{\gamma}{\alpha}}\|u\|_{\mH^{\gamma+\theta, \infty}_p(S,T)}\leq c\|f\|_{\mH^{\theta,q}_p(S, T)}.
\end{align}
Moreover, the estimate \eqref{NR1} still holds.
\et
\begin{proof}
By the standard Picard's iteration, it suffices to show the apriori estimate \eqref{GR3}.
Let $q^*=q/(q-1)$. Since $\gamma\in(0, \alpha-\frac{\alpha}{q})$ and $\gamma+\theta<\alpha+\eps$,
by \eqref{In}, \eqref{UT4} and Lemma \ref{Le61}, we have
\begin{align*}
\|u(s)\|_{\gamma+\theta, p}&\leq \int^T_s\e^{-\lambda(t-s)}\|P^{\kappa,0}_{s,t}(\sL^b_1 u+f)\|_{\gamma+\theta,p}\dif t\\
&\lesssim \int^T_s\e^{-\lambda(t-s)}(t-s)^{-\gamma/\alpha}(\|b\cdot\nabla u\|_{\theta,p}+\|f\|_{\theta,p})\dif t\\
&\lesssim \int^T_s\e^{-\lambda(t-s)}(t-s)^{-\gamma/\alpha}(\|b\|_{\theta,p}\|u\|_{\gamma+\theta,p}+\|f\|_{\theta,p})\dif t,
\end{align*}
which implies by Volterra-Gronwall's inequality (cf. \cite[Lemma 4.1]{Wa-Zh1}) and H\"older's inequality that
\begin{align*}
\|u(s)\|_{\gamma+\theta, p}
&\lesssim\int^T_s\e^{-\lambda(t-s)}(t-s)^{-\gamma/\alpha}\|f\|_{\theta,p}\dif t\\
&\leq \left(\int^T_s\e^{-q^*\lambda(t-s)}(t-s)^{-q^*\gamma/\alpha}\dif t\right)^{1/q^*}\|f\|_{\mH^{\theta,q}_p(s,T)}\\
&\lesssim(\lambda\vee 1)^{\frac{1}{q}+\frac{\gamma}{\alpha}-1}\|f\|_{\mH^{\theta,q}_p(s,T)}.
\end{align*}
Thus, we get the desired estimate \eqref{GR3}.
\end{proof}

\section{Krylov's estimate for semimartingales}

This section is devoted to the study of Krylov's estimates for discontinuous semimartingales,
which can be regarded as apriori estimates for the solution of SDE \eqref{sde1}.

\subsection{General discontinuous semimartingales}
The classical Krylov estimate on the distribution of continuous martingales is well known,
see \cite{Kry} or \cite[Lemma 3.1]{Gy-Ma}.
%, and has been widely used (see \cite{M-N-P-Z,Zh4}). 
Below, we generalize it to discontinuous semimartingales.

The following important result on the existence of a solution for a partial differential inequality comes from Krylov \cite[Chapter III, Theorem 2.4]{Kry}.

\bl\label{pdi}
Given a nonnegative smooth function $f$ on $\mR_{+}\times\mR^d$ and $\lambda>0$, there exists a nonnegative smooth function $u(t,x)$ such that
for all nonnegative definite symmetric matrices $a=(a^{ij})_{d\times d}$ and $\beta\geq 0$,
\begin{align}
&\beta\p_t u+a^{ij}\p_i\p_{j}u-\lambda(\beta+\tr a) u+(\beta\det a)^{1/(d+1)}f\leq 0,\label{k4}
\end{align}
and
\begin{align}
&|\nabla u|\leq \sqrt{\lambda}u,\quad \ u\leq K_{d}\lambda^{-d/(2(d+1))}\|f\|_{\mL^{d+1}(T)},    \label{k2}
\end{align}
where $K_{d}>0$ depends only on the dimension $d$.
\el

Using this lemma, we can show the following Krylov estimate for general discontinuous semimartingales.

\bl\label{kk}
Let $m=m(t)$ be an  $\mR^d$-valued continuous local martingale, $V=V(t)$ an $\mR^d$-valued
continuous adapted process with finite variation on finite time intervals,
$N(\dif t,\dif z)$  a Poisson random measure with compensator $\dif t\nu(\dif z)$, where $\nu$ is a L\'evy measure,
and $G:\mR_+\times\Omega\times\mR^d\rightarrow\mR^d$ a predictable process with
$$
\int^t_0\!\!\!\int_{|z|\leq R}|G(s,z)|^2\dif s\nu(\dif z)<\infty,\ \ a.s.,
$$
where $R>0$.
Suppose that
$$
m(0) = V(0) = 0,\ \ \dif \<m^i, m^j\>_t\ll\dif t.
$$
Let $a^{ij}(t):=\tfrac{\dif\<m^i,m^j\>_t}{2\dif t}$ and
$$
X(t):=m(t)+ V(t)+\int_0^t\!\!\!\int_{|z|\leq  R}G(s,z)\tilde N(\dif s,\dif z)
+\int_0^t\!\!\!\int_{|z|> R}G(s,z)N(\dif s,\dif z).
$$
Then for any $T>0$, $p\geq d+1$ and  $\alpha\in[1,2]$, there is a constant $c=c(T,p,d,\alpha)>0$
such that for any stopping time $\tau$ and $f\in\mL^p(T)$,
\begin{align}\label{Kr1}
\mE\left(\int_0^{T\wedge\tau}\!\!\big(\det a(t)\big)^{\frac{1}{p}}f(t,X_t)\dif t\right)\leq
c\Big(\mV^2+\mA+\mG_\alpha^{\frac{2}{\alpha}}\Big)^{\frac{d}{2p}}\|f\|_{\mL^{p}(T)},
\end{align}
where
\begin{equation*}
\begin{aligned}
&\mV:=\mE\left(\int_0^{T\wedge\tau}\!|\dif V(t)|\right),\,\, \mA:=\mE\left(\int_0^{T\wedge\tau}\!\!\tr \,a(t)\dif t\right),\\
&\mG_\alpha:=\mE\left(\int_0^{T\wedge\tau}\!\!\!\int_{|z|<R}|G(t,z)|^{\alpha}\dif t\nu(\dif z)\right).
\end{aligned}
\end{equation*}
\el
\begin{proof}
By standard approximation, we may assume that $f\in C^{\infty}_0(\mR_+\times\mR^d)$ and $\mV,\mA,\mG_\alpha$ are finite.
For a given constant $\lambda>0$ whose precise value will be decided latter, let $u$ be the nonnegative smooth function given by Lemma \ref{pdi}
corresponding to $\lambda$ and $f$. By It\^o's formula, we have
\begin{align*}
Z_t&:=u(t,X_t)-\int_0^t\!\Big(\p_s u+a^{ij}\p_{ij}u+\sL^G_{\nu}u\Big)(s,X_s)\dif s
-\int_0^t\p_i u(s,X_s)\dif V^i_s\\
&=\int^t_0\p_i u(s,X_s)\dif m^i_s+\int^t_0\!\!\!\int_{\mR^d}\big(u(s, X_{s-}+G(s,z))-u(s,X_{s-})\big)\tilde N(\dif s,\dif z),
\end{align*}
is a local martingale, where
$$
\sL^G_{\nu}u(t,x):=\int_{\mR^d}\big[u\big(t,x+G(t,z)\big)-u(t,x)-1_{|z|\leq  R}G^i(t,z)\p_i u(t,x)\big]\nu(\dif z).
$$
Observing that for $|z|<R$,
\begin{align*}
\Sigma_t(x,z)&:=u\big(t,x+G(t,z)\big)-u(t,x)-G^i(t,z)\p_i u(t,x)\\
&=G^i(t,z)G^j(t,z)\int^1_0\!\!\int^1_0\p_i\p_ju\big(t,x+s_1s_2G(t,z)\big)\dif s_1\dif s_2,
\end{align*}
by (\ref{k4}) with $\beta=0$, we have
$$
\Sigma_t(x,z)\leq \lambda|G(t,z)|^2\|u\|_{\mL^\infty(T)}.
$$
Moreover, by \eqref{k2} we also have
$$
|\Sigma_t(x,z)|\leq\sqrt{\lambda}|G(t,z)|\|u\|_{\mL^\infty(T)}.
$$
Hence, for any $\alpha\in[1,2]$,
\begin{align}\label{k5}
\sL^G_{\nu}u\leq \left(2\nu(B^c_R)+\lambda^{\frac{\alpha}{2}}\int_{|z|<R}|G(t,z)|^{\alpha}\nu(\dif z)\right)\|u\|_{\mL^\infty(T)}.
\end{align}
For $n>|u(0,0)|$, if we define the stopping time
$$
\tau_n:=\tau\wedge\inf\{t\geq 0: |Z_t|\geq n\},
$$
then $t\mapsto Z_{t\wedge\tau_n}$ is a bounded martingale.
Thus, by  the definition of $Z_t$, \eqref{k4} with $\beta=1$, \eqref{k2} and \eqref{k5},  we have
\begin{align*}
\mE u(t\wedge\tau_n,X_{t\wedge\tau_n})-u(0,0)&\leq -\mE\left(\int_0^{t\wedge\tau_n}\!\big(\det a(s)\big)^{\frac{1}{d+1}}f(s,X_s)\dif s\right)\\
&+\mE\bigg(\sqrt{\lambda}\!\int_0^{t\wedge\tau_n}\!\dif |V_s|+\lambda\!\int_0^{t\wedge\tau_n}\!(\tr a(s)+1)\dif s\\
&+2\nu(B^c_R)+\lambda^{\frac{\alpha}{2}}\!\!\int_0^{t\wedge\tau_n}\!\!\!\int_{|z|<R}|G(t,z)|^{\alpha}\nu(\dif z)\dif s\bigg)\|u\|_{\mL^\infty(T)}\\
&\leq -\mE\left(\int_0^{t\wedge\tau_n}\!\big(\det a(s)\big)^{\frac{1}{d+1}}f(s,X_s)\dif s\right)\\
&\quad+\Big(\sqrt{\lambda}\mV+\lambda(\mA+t)+2\nu(B^c_R)+\lambda^{\frac{\alpha}{2}}\mG_\alpha\Big)\|u\|_{\mL^\infty(T)}.
\end{align*}
Taking into account (\ref{k2}),  we get
\begin{align*}
&\mE\left(\int_0^{t\wedge\tau_n}\!\big(\det a(s)\big)^{\frac{1}{d+1}}f(s,X_s)\dif s\right)\\
&\quad\lesssim\Big(\sqrt{\lambda}\mV+\lambda(\mA+1)+\lambda^{\frac{\alpha}{2}}\mG_\alpha+1\Big)\lambda^{-d/(2(d+1))}\|f\|_{\mL^{d+1}(T)},
\end{align*}
which, by taking $\lambda^{-1}=\mV^2\vee\mA\vee\mG_\alpha^{\frac{2}{\alpha}}\vee 1$ and letting $n\to\infty$,  implies  \eqref{Kr1} for $p=d+1$.
Finally, for $p>d+1$, by H\"older's inequality, we have
\begin{align*}
\mE\left(\int_0^{T\wedge\tau}\!\!\big(\det a(t)\big)^{\frac{1}{p}}f(t,X_t)\dif t\right)&\lesssim\left(\mE\!\int_0^{T\wedge\tau}\!\!\big(\det a(t)\big)^{\frac{1}{d+1}}
|f(t,X_t)|^{\frac{p}{d+1}}\dif t\right)^{\frac{d+1}{p}}\\
&\lesssim\Big(1+\mV^2+\mA+\mG_\alpha^{\frac{2}{\alpha}}\Big)^{\frac{d}{2p}}\|f\|_{\mL^p(T)}.
\end{align*}
The proof is finished.
\end{proof}

\subsection{Non-degenerate diffusion SDEs with jumps}

Below, for the moment we suppose that $X_t$ satisfies the following equation:
\begin{align}\label{form}
\begin{split}
X_t&=X_0+\int_0^t\!\sigma_s(X_s)\dif W_s+\int_0^t\!\!\!\int_{|z|\leq  R}g_s(X_{s-},z)\tilde N(\dif s, \dif z)\\
&\qquad+\int_0^t\!\!\!\int_{|z|> R}g_s(X_{s-},z)N(\dif s, \dif z)+\int_0^t\xi(s)\dif s,
\end{split}
\end{align}
where $\xi(t)$ is a measurable $\sF_t$-adapted process.
The following lemma is an easy consequence of Lemma \ref{kk}.
\bl
Let $X_t$ be of the form (\ref{form}). Suppose that $\sigma\sigma^*$ is bounded and uniformly positive definite,
and for some $\alpha\in[1,2]$ and $q\geq d+1$,
$\Gamma^{0,\alpha}_{0, R}(g)\in\mL^q(T)$, where $\Gamma^{0,\alpha}_{0, R}(g)$ is defined by \eqref{LK1}.
Then for any $p\geq d+1$ and $\delta>0$, there is a constant $c_\delta>0$ such that for any stopping time $\tau$ and $f\in \mL^p(T)$,
\begin{align}\label{kry01}
\mE\left(\int_0^{T\wedge\tau}\!\!f(s,X_s)\dif s\right)\leq \left(c_\delta+\delta\mE\left(\int_{0}^{T\wedge\tau}|\xi(s)|\dif s\right)\right)\|f\|_{\mL^p(T)}.
\end{align}
\el
\begin{proof}
Without loss of generality, we assume $\mE\left(\int_{0}^{T\wedge\tau}|\xi(s)|\dif s\right)<\infty$. In order to use Lemma \ref{kk}, we take
$$
m(t):=\int^t_0\sigma_s(X_s)\dif W_s,\ \ V(t):=\int^t_0|\xi(s)|\dif s,\ \ G(t,z):=g_t(X_{t-},z).
$$
Thus, by the assumption on $\sigma$, for any $p\geq d+1$, by \eqref{Kr1}, there is a constant $c>0$ such that for all $f\in\mL^p(T)$,
\begin{align}\label{Kr11}
\mE\left(\int_0^{T\wedge\tau}f(t,X_t)\dif t\right)\leq
c\Big(1+\mV^2+\mG_\alpha^{\frac{2}{\alpha}}\Big)^{\frac{d}{2p}}\|f\|_{\mL^{p}(T)}.
\end{align}
Here, $\mV:=\mE\left(\int^{T\wedge\tau}_0|\xi(s)|\dif s\right)$ and
$$
\mG_\alpha:=\mE\left(\int_0^{T\wedge\tau}\!\!\!\int_{|z|<R}|g_t(X_{t},z)|^{\alpha}\dif t\nu(\dif z)\right)
=\mE\left(\int_0^{T\wedge\tau}\Gamma^{0,\alpha}_{0, R}(g_t)(X_{t})\dif t\right).
$$
By \eqref{Kr11} with $f=\Gamma^{0,\alpha}_{0, R}(g)$ and the assumption, we have
$$
\mG_\alpha\leq c\Big(1+\mV^2+\mG_\alpha^{\frac{2}{\alpha}}\Big)^{\frac{d}{2q}}\|\Gamma^{0,\alpha}_{0, R}(g)\|_{\mL^q(T)}
\leq c(1+\mV^{d/q})+\tfrac{1}{2}\mG_\alpha,
$$
which implies $\mG_\alpha\leq c(1+\mV^{d/q})$.
Thus, we get \eqref{kry01} by \eqref{Kr11} and Young's inequality.
\end{proof}
In the above estimate, it is required $p\geq d+1$, which is too strong for our purpose.
Below we use Theorem \ref{pde} to obtain better integrability index $p$. The price  we
have to pay is to strengthen the assumption on $\Gamma^{0,\alpha}_{0, R}(g)$ since we need to use Theorem \ref{pde}.
\bl\label{lem}
Let $X$ be of the form (\ref{form}) and $\Gamma^{0,2}_{0,R}(g)$ be defined by \eqref{LK1}. Let $T>0$. Suppose  that {\bf (H$^\sigma$)} holds and
$$
\Gamma^{0,2}_{0,R}(g)\in\mL^\infty(T)\mbox{ and }\lim_{\eps\to 0}\|\Gamma^{0,2}_{0,\eps}(g)\|_{\mL^\infty(T)}=0.
$$
Then for any $p,q\in(1,\infty)$ with $\frac{d}{p}+\frac{2}{q}<1$ and each $\delta>0$, there is a constant $c_\delta>0$ such that for any stopping time $\tau$ and
$0\leq t_0\leq t_1\leq T$ and $f\in \mL^q_p(t_0, t_1)$,
\begin{align}\label{UY8}
\mE\left(\int_{t_0\wedge\tau}^{t_1\wedge\tau}\!\!f(s,X_s)\dif s\Big|\sF_{t_0\wedge\tau}\right)\leq
\|f\|_{\mL^q_p(t_0,t_1)}\left[c_\delta+\delta\mE\left(\int_{t_0\wedge\tau}^{t_1\wedge\tau}\!\!|\xi(s)|\dif s\Big|\sF_{t_0\wedge\tau}\right)\right].
\end{align}
Moreover, if $\xi\equiv 0$, then we can relax $p,q$ to satisfy $\frac{d}{p}+\frac{2}{q}<2$.
\el
\begin{proof}
We may assume without loss of generality that $f\in C^{\infty}_0(\mR^{d+1})$ and
$$
\mE\left(\int_0^{T\wedge\tau}|\xi(s)|\dif s\right)<+\infty.
$$
Let $r$ be large enough so that
$$
\tfrac{d}{r}+\tfrac{2}{r}\leq\tfrac{d}{p}+\tfrac{2}{q}<1.
$$
Let $\lambda_0$ be the constant in Theorem \ref{pde}. For $\lambda\geq\lambda_0$ and $t_1\in(0,T]$, since $f\in\mL^q_p(t_1)\cap\mL^r(t_1)$,
by Theorem \ref{pde}, there exists a unique solution $u\in \mH^{2,q}_{p}(t_1)\cap\mH^{2,r}_r(t_1)$ with $\p_tu\in \mL^{r}(t_1)$
to the following backward equation:
$$
\p_tu+(\sL^a_2-\lambda)u+\sL^g_{\nu, R} u=f,\quad u(t_1)=0,
$$
where $a=\sigma\sigma^*/2$.
Let $\phi$ be a non-negative smooth function on $\mR^{d+1}$ with support in $\{x\in\mR^{d+1}: |x|\leq 1\}$ and
$\int_{\mR^{d+1}}\phi(t,x)\dif t\dif x=1$. Set
$$
\phi_n(t,x):=n^{d+1}\phi(nt,nx)
$$
and extend $u(t,x)$ to $\mR$
by setting $u(t,x)=0$ for $t\geq t_1$ and $u(t,x)=u(0,x)$ for $t\leq 0$. Define
\begin{align}\label{CON}
u_n(t,x):=u*\phi_n(t,x):=\int_{\mR^{d+1}}u(s,y)\phi_n(t-s,x-y)\dif s\dif y
\end{align}
and
\begin{align}  \label{un}
f_n:=\p_tu_n+(\sL^a_2-\lambda)u_n+\sL^g_{\nu, R}u_n,
\end{align}
where $\sL^g_{\nu, R}$ is defined by \eqref{UY7}.
Since $\tfrac{d}{p}+\tfrac{2}{q}<1$,
by the property of convolution and using \eqref{es1} with $\gamma=1$ and $p'=q'=\infty$, there is a constant
$c>0$ independent of $n$ such that for all $\lambda\geq 1$ and $t_0\in[0,t_1]$,
\begin{align}\label{BG3}
\|u_n\|_{\mH^{1,\infty}_\infty(t_0, t_1)}\leq \|u\|_{\mH^{1,\infty}_\infty(t_0, t_1)}
\leq c\lambda^{\frac{1}{2}(\frac{d}{p}+\frac{2}{q}-1)}\|f\|_{\mL^q_p(t_0,t_1)},
\end{align}
and
\begin{align*}
\|f_n-f\|_{\mL^r_r(t_1)}&\leq \lambda\|u_n-u\|_{\mL^r_r(t_1)}+\|\p_t (u_n-u)\|_{\mL^r_r(t_1)}\\
&\quad+c\|\nabla^2(u_n-u)\|_{\mL^r_r(t_1)}+\|\sL^g_{\nu,R}(u_n-u)\|_{\mL^r_r(t_1)}\\
&\leq \|\p_t (u_n-u)\|_{\mL^r_r(t_1)}+c\|u_n-u\|_{\mH^{2,r}_{r}(t_1)}\rightarrow0,\quad n\rightarrow\infty,
\end{align*}
where we have used the same estimate as in \eqref{BB2}.
Therefore, by the Krylov estimate (\ref{kry01}), we have
\begin{align}\label{Lim0}
\lim_{n\rightarrow\infty}\mE\left(\int_0^{T\wedge\tau}\!\big|f_n(s,X_s)-f(s,X_s)\big|\dif s\right)\leq c\lim_{n\rightarrow\infty}\|f_n-f\|_{\mL^r_r(T)} =0.
\end{align}
Now, applying It\^o's formula to $u_n(t,x)$ , we have
\begin{align*}
u_n(t,X_t)&=u_n(0,X_0)+\int_0^{t}\!\!\Big(\p_s u_n+\sL^a_{2}u_n+\sL^g_{\nu} u_n\Big)(s,X_s)\dif s\\
&\quad+\int_0^{t}\!\!\xi(s)\cdot\nabla u_n(s,X_s)\dif s+\mbox{a martingale}.
\end{align*}
Thus, by Doob's optional stopping theorem and \eqref{un}, we obtain
\begin{align*}
&\quad\mE\Big( u_n(t_1\wedge\tau,X_{t_1\wedge\tau})|{\sF_{t_0\wedge\tau}}\Big)-u_n(t_0\wedge\tau,X_{t_0\wedge\tau})\\
&=\mE\left(\int_{t_0\wedge\tau}^{t_1\wedge\tau}\!\!\Big(\p_s u_n+\sL^a_{2}u_n+\sL^g_{\nu} u_n\Big)(s,X_s)\dif s\Big|{\sF_{t_0\wedge\tau}}\right)\\
&\quad+\mE\left(\int_{t_0\wedge\tau}^{t_1\wedge\tau}\xi(s)\cdot\nabla u_n(s,X_s)\dif s\Big|{\sF_{t_0\wedge\tau}}\right)\\
&\geq \mE\left(\int_{t_0\wedge\tau}^{t_1\wedge\tau}\!\!\Big(\lambda u_n(s,X_s)+f_n(s,X_s)\Big)\dif s\Big|{\sF_{t_0\wedge\tau}}\right)\\
&\quad-2\|u_n\|_{\mL^\infty(t_0,t_1)}\nu(B^c_R) t_1
-\|\nabla u_n\|_{\mL^\infty(t_0,t_1)}\mE\left(\int_{t_0\wedge\tau}^{t_1\wedge\tau}|\xi(s)|\dif s\Big|{\sF_{t_0\wedge\tau}}\right),
\end{align*}
which implies that by \eqref{BG3},
\begin{align*}
&\mE\left(\int_{t_0\wedge\tau}^{t_1\wedge\tau}f_n(s,X_s)\dif s\Big|{\sF_{t_0\wedge\tau}}\right)
\leq \Big(2+\lambda T+2\nu(B^c_R) T\Big)\|u_n\|_{\mL^\infty(t_0,t_1)}\\
&\qquad\qquad+\|\nabla u_n\|_{\mL^\infty(t_0,t_1)}\mE\left(\int_{t_0\wedge\tau}^{t_1\wedge\tau}|\xi(s)|\dif s\Big|{\sF_{t_0\wedge\tau}}\right)\\
&\qquad\leq \left[c_\lambda+c\lambda^{\frac{1}{2}(\frac{d}{p}+\frac{2}{q}-1)}
\mE\left(\int_{t_0\wedge\tau}^{t_1\wedge\tau}|\xi(s)|\dif s\Big|{\sF_{t_0\wedge\tau}}\right)\right]\|f\|_{\mL^q_p(t_0,t_1)}.
\end{align*}
Letting $n\to\infty$ and $\lambda$ be large enough, by \eqref{Lim0} we get \eqref{UY8}.
If $\xi\equiv 0$, then we only need to control $\|u\|_{\mL^\infty(t_0,t_1)}$, which follows by \eqref{es1} with $\vartheta=0$ and $p'=q'=\infty$.
\end{proof}
\br
Lemma \ref{lem} will be used to derive the Krylov estimate for SDE with polynomial growth drift in the proof of ergodicity for SDEs with singular drifts.
\er

%As a direct consequence, we have the following Krvlov estimate for the solutions of SDE (\ref{sde1}).
We also have the following Krvlov estimate for the solutions of SDE (\ref{sde1}).

\bt\label{enou}
Let $T>0$. Assume that {\bf (H$^\sigma$)} holds and for some $p_1,q_1\in(2,\infty]$ with
$\frac{d}{p_1}+\frac{2}{q_1}<1$,
$$
b\in\mL^{q_1}_{p_1}(T), \ \Gamma^{0,2}_{0,R}(g)\in\mL^\infty(T), \ \  \lim_{\eps\to 0}\|\Gamma^{0,2}_{0,\eps}(g)\|_{\mL^\infty(T)}=0.
$$
Then for any $p,q\in(1,\infty)$ with $\tfrac{d}{p}+\tfrac{2}{q}<2$,  the solution $X$ of SDE (\ref{sde1}) satisfies  Krylov's estimate with index $p,q$.
\et
\begin{proof}
(i) First of all, we show that $X$ satisfies  Krylov's estimate for all $p,q\in(1,\infty)$ with $\tfrac{d}{p}+\tfrac{2}{q}<1$.
By Lemma \ref{lem}, it suffices to show that for all $0\leq t_0\leq t_1\leq T$,
\begin{align}\label{HD1}
\mE\left(\int_{t_0}^{t_1}|b_s(X_s)|\dif s\Big|{\sF_{t_0}}\right)\leq c\|b\|_{\mL^{q_1}_{p_1}(t_0,t_1)}.
\end{align}
For $n\in\mN$, define a stopping time
$$
\tau_n:=\inf\left\{t>0: \int^t_0 |b_s(X_s)|\dif s\geq n\right\}.
$$
Taking $\xi(s)=b_s(X_s)$ and $f=|b|$ in Lemma \ref{lem}, we get that for every $\delta>0$
and $0\leq t_0\leq t_1\leq T$,
\begin{align*}
&\mE\left(\int_{t_0\wedge\tau_n}^{t_1\wedge\tau_n}\!\!|b_s(X_s)|\dif s\Big|{\sF_{t_0\wedge\tau_n}}\right)
\leq \left[c_\delta+\delta\mE\left(\int_{t_0\wedge\tau_n}^{t_1\wedge\tau_n}\!\!|b_s(X_s)|\dif s\Big|{\sF_{t_0\wedge\tau_n}}\right)
\right]\|b\|_{\mL^{q_1}_{p_1}(t_0, t_1)}.
\end{align*}
Choosing $\delta$ be small enough such that
$$
\delta\|b\|_{\mL^{q_1}_{p_1}(T)}<\tfrac{1}{2},
$$
we obtain that for all $0\leq t_0\leq t_1\leq T$,
$$
\mE\left(\int_{t_0\wedge\tau_n}^{t_1\wedge\tau_n}|b_s(X_s)|\dif s\Big|{\sF_{t_0\wedge\tau_n}}\right)\leq c\|b\|_{\mL^{q_1}_{p_1}(t_0,t_1)},
$$
where $c$ is independent of $n$.
Letting $n\rightarrow\infty$, we get \eqref{HD1}.

(ii) In this step we show that $X$ satisfies the Krylov estimate for $p=p_1/2$ and $q=q_1/2$.
Without loss of generality, we assume $p_1,q_1\in(2,\infty)$.
Assume $f\in C^{\infty}_0(\mR^{d+1})$. Let $\lambda_0$ be the constant in Theorem \ref{pde}. For $\lambda\geq\lambda_0$ and $t_1\in(0,T]$,
since $f\in\mL^q_p(t_1)\cap\mL^{q_1}_{p_1}(t_1)$,
by Theorem \ref{pde}, there exists a unique solution $u\in \mH^{2,q}_{p}(t_1)\cap\mH^{2,q_1}_{p_1}(t_1)$ with $\p_tu\in \mL^{q_1}_{p_1}(t_1)$
to the following backward equation:
$$
\p_tu+(\sL^a_2-\lambda)u+\sL^g_{\nu, R} u+\sL^{b}_1u=f,\quad u(t_1)=0.
$$
Let $u_n:=u*\phi_n$ be defined as in \eqref{CON}, and
$$
f_n:=\p_tu_n+(\sL^a_2-\lambda)u_n+\sL^g_{\nu, R}u_n+\sL^{b}_1u_n.
$$ 
As in the proof of Lemma \ref{lem} we have
$$
\mE\Big( u_n(t_1,X_{t_1})|{\sF_{t_0}}\Big)-u_n(t_0,X_{t_0})\\
=\mE\left(\int_{t_0}^{t_1}(f_n+\lambda u_n)(s,X_s)\dif s\Big|{\sF_{t_0}}\right),
$$
which implies by \eqref{es1} with $\vartheta=0$ and $p'=q'=\infty$ that
\begin{align}\label{HD5}
\mE\left(\int_{t_0}^{t_1}f_n(s,X_s)\dif s\Big|{\sF_{t_0}}\right)
\leq (\lambda T+2)\|u_n\|_{\mL^\infty(t_1)}\leq c\|f\|_{\mL^q_p(t_1)}.
\end{align}
Noticing that 
$$
\lim_{n\to\infty}\|f_n-f\|_{\mL^{q_1}_{p_1}(t_1)}=0,
$$ 
by step (i) and taking limits $n\to\infty$ for \eqref{HD5}, we get
\begin{align}\label{HD4}
\mE\left(\int_{t_0}^{t_1}f(s,X_s)\dif s\Big|{\sF_{t_0}}\right)
\leq c\|f\|_{\mL^q_p(t_1)}=c\|f\|_{\mL^{q_1/2}_{p_1/2}(t_1)}.
\end{align}

(iii) By \eqref{HD4}, we have for all $0\leq t_0\leq t_1\leq T$,
$$
\mE\left(\int_{t_0}^{t_1}|b_s(X_s)|^2\dif s\Big|{\sF_{t_0}}\right)\leq c\|b\|_{\mL^{q_1}_{p_1}(t_0,t_1)}. 
$$
By Lemma \ref{Cor1}, for any $\lambda>0$, there is a constant $c>0$ such that for all $0\leq t_0<t_1\leq T$,
\begin{align}\label{BG1}
\mE\left(\exp\left\{\lambda\int_{t_0}^{t_1}|b_s(X_s)|^2\dif s\right\}\Big|{\sF_{t_0}}\right)\leq c.
\end{align}
Define for $\gamma\in\mR$,
$$
\cE^{(\gamma)}_{t_0,t_1}:=\exp\left\{\gamma\int^{t_1}_{t_0} (\sigma^{-1}_sb_s)(X_s)\dif W_s-\frac{\gamma^2}{2}\int^{t_1}_{t_0}|\sigma^{-1}_sb_s|^2(X_s)\dif s\right\}.
$$
By Novikov's criterion, $t\mapsto\cE^{(\gamma)}_{0,t}$ is an exponential martingale. Hence, by \eqref{BG1} and H\"older's inequality,
\begin{align}\label{KJ2}
\mE\Big((\cE^{(1)}_{t_0,t_1})^\gamma|\sF_{t_0}\Big)\leq
\left(\mE\left(\exp\left\{(2\gamma^2-\gamma)\int_{t_0}^{t_1}|\sigma^{-1}_sb_s|^2(X_s)\dif s\right\}\Big|{\sF_{t_0}}\right)\right)^{1/2}\leq c.
\end{align}
Define a new probability $\mQ_{t_0, t_1}:=\cE^{(1)}_{t_0,t_1}\mP$. 
By Girsanov's theorem, under the probability measure $\mQ_{t_0,t_1}$, after time $t_0$,
$\tilde W_t:=W_t+\int^t_{t_0}(\sigma^{-1}_sb_s)(X_s)\dif s$ is still a Brownian motion and
$N(\dif t,\dif z)$ is still a Poisson random measure with the same compensator $\dif t\nu(\dif z)$. Moreover,
$X_t$ satisfies
\begin{align*}
X_t&=X_{t_0}+\int^t_{t_0}\sigma_s(X_s)\dif \tilde W_s+\int^t_{t_0}\!\!\int_{|z|<R}g_s(X_{s-},z)\tilde N(\dif s,\dif z)
+\int^t_{t_0}\!\!\int_{|z|> R}g_s(X_{s-},z)N(\dif s,\dif z).
\end{align*}
Hence, by Lemma \ref{lem} with $\xi\equiv 0$, for any $p,q\in(1,\infty)$ with $\tfrac{d}{p}+\tfrac{2}{q}<2$,
\begin{align}\label{HD8}
\mE^{\mQ_{t_0,t_1}}\left(\int_{t_0}^{t_1}\!\!f\big(s,X_s\big)\dif s\bigg|\sF_{t_0}\right)\leq c\|f\|_{\mL^{q}_{p}(t_0,t_1)}.
\end{align}
Noticing that for any nonnegative random variable $\zeta$,
$$
\mE\Big(\zeta\cE^{(1)}_{t_0,t_1}|\sF_{t_0}\Big)=\mE^{\mQ_{t_0,t_1}}\Big(\zeta|\sF_{t_0}\Big)\mE\Big(\cE^{(1)}_{t_0,t_1}|\sF_{t_0}\Big),
$$
by \eqref{HD8},  \eqref{KJ2} and suitable H\"older's inequality, we get the desired Krylov estimate.
\end{proof}

\subsection{SDEs driven by pure jump L\'evy noises}
In this subsection we assume $\nu(\dif z)=\dif z/|z|^{d+\alpha}$ for some $\alpha\in(1,2)$,
and show the Krylov estimate for pure jump cases.
First of all, we have
\bl\label{Le57}
Under {\bf (H$^g$)}, it holds that
\begin{align}\label{499}
\sL^g_\nu u(x)=\sL^\kappa_\alpha u(x)+\bar b^g_t(x)\cdot\nabla u(x),
\end{align}
where $\sL^\kappa_\alpha$ is defined by \eqref{49} with $\kappa$ satisfying \eqref{Con2},
$\bar b^g_t(x)$ is bounded and H\"older continuous in $x$ uniformly with respect to $t$.
\el
\begin{proof}
By \eqref{HF3}, one sees that \eqref{499} holds.
We now check that $\kappa$ has the desired property \eqref{Con2}.
By \eqref{CN1} and $g_t(x,0)=0$, the map $z\mapsto g_t(x,z)$ admits an inverse $g^{-1}_t(x,z)$ so that
$$
g^{-1}_t(x,0)=0,\ \ c_1^{-1}|z-z'|\leq|g^{-1}_t(x,z)-g^{-1}_t(x,z')|\leq c_1|z-z'|.
$$
In particular,
\begin{align}\label{HF1}
c^{-1}_1|z|\leq |g^{-1}_t(x,z)|\leq c_1|z|,\quad\|\nabla_z g^{-1}\|_\infty\leq c_1.
\end{align}
Moreover, for $x,y\in\mR^d$, letting $\tilde z:=g^{-1}_t(x,z)$, we have
\begin{align}\label{HF2}
\begin{split}
|g^{-1}_t(x,z)-g^{-1}_t(y,z)|&=|g^{-1}_t(y,g_t(y,\tilde z))-g^{-1}_t(y,g_t(x,\tilde z))|\\
&\leq c_1 |g_t(y,\tilde z)-g_t(x,\tilde z)|\leq c^2_1 |x-y|^\beta|\tilde z|\leq c^3_1|x-y|^\beta|z|.
\end{split}
\end{align}
Noticing that
$$
\nabla_zg^{-1}_t(x,z)=[\nabla_zg_t]^{-1}(x,g^{-1}_t(x,z)).
$$
by \eqref{HF1} and \eqref{HF2}, it is easy to see that the $\kappa$ defined by \eqref{ka} satisfies \eqref{Con2}.
\end{proof}
The following lemma is similar to Lemma \ref{lem}.
\bl\label{Le63}
Suppose that {\bf (H$^g$)} holds and $X_t$ satisfies
$$
X_t=X_0+\int^t_0\!\!\!\int_{|z|<R}g_s(X_{s-},z)\tilde N(\dif s,\dif z)+\int^t_0\!\!\!\int_{|z|\geq R}\eta_s(z)N(\dif s,\dif z)+\int^t_0\xi(s)\dif s,
$$
where $\eta:\mR_+\times\Omega\times\mR^d\to\mR^d$ is a predictable process.
For any $p,q\in[1,\infty]$ with $\frac{d}{p}+\frac{\alpha}{q}<\alpha-1$ and each $\delta>0$,
there is a constant $c_\delta>0$ such that for all $0\leq t_0<t_1\leq T$, any stopping time $\tau$ and $f\in\mL^q_p(T)$,
\begin{align}\label{UU4}
\begin{split}
\mE\left(\int_{t_0\wedge\tau}^{t_1\wedge\tau}f(s,X_s)\dif s\Big|{\sF_{t_0\wedge\tau}}\right)
\leq \left[c_\delta+\delta\mE\left(\int_{t_0\wedge\tau}^{t_1\wedge\tau}|\xi(s)|\dif s\Big|{\sF_{t_0\wedge\tau}}\right)\right]\|f\|_{\mL^q_p(T)}.
\end{split}
\end{align}
\el
\begin{proof}
Without loss of generality, we assume $f\in C^\infty_c(\mR_+\times\mR^d)$.
By Lemma \ref{Le57}, Theorem \ref{Th63} and Remark \ref{Rem44},  for some $\eps\in(0,2-\alpha)$, there is a unique
$u\in\mH^{\alpha+\eps,\infty}_\infty(T)$ solving the following equation
$$
\p_t u+\sL^g_\nu u-\lambda u=f,\ u(t_1)=0.
$$
By It\^o's formula and Doob's optional stopping theorem, we have
\begin{align*}
&\quad\mE\Big( u(t_1\wedge\tau,X_{t_1\wedge\tau})|{\sF_{t_0\wedge\tau}}\Big)-u(t_0\wedge\tau,X_{t_0\wedge\tau})\\
&=\mE\left(\int_{t_0\wedge\tau}^{t_1\wedge\tau}\!\!\Big(\p_s u_n+\sL^g_{\nu,R}u_n+\bar\sL^\eta_{\nu, R}u_n+\sL^b_1 u_n\Big)(s,X_s)\dif s\Big|{\sF_{t_0\wedge\tau}}\right)\\
&=\mE\left(\int_{t_0\wedge\tau}^{t_1\wedge\tau}\Big(\big(\lambda u+f-\bar\sL^g_{\nu,R}u+\bar\sL^\eta_{\nu,R}u\big)(s,X_s)
+\xi(s)\cdot\nabla u(s,X_s)\Big)\dif s\Big|{\sF_{t_0\wedge\tau}}\right),
\end{align*}
where $\bar\sL^g_{\nu,R}$ and $\bar\sL^\eta_{\nu,R}$ are defined as in \eqref{UY7}.
Hence, by Theorem \ref{Th62} with $\vartheta=1$ and $q'=p'=\infty$, we get for $\lambda\geq 1$,
\begin{align*}
&\mE\left(\int_{t_0\wedge\tau}^{t_1\wedge\tau}f(s,X_s)\dif s\Big|{\sF_{t_0\wedge\tau}}\right)
\leq(t_1-t_0)(\lambda+2+4\nu(B^c_R))\|u\|_{\mL^\infty(T)}\\
&\qquad\qquad+\|\nabla u\|_{\mL^\infty(T)}
\mE\left(\int_{t_0\wedge\tau}^{t_1\wedge\tau}|\xi(s)|\dif s\Big|{\sF_{t_0\wedge\tau}}\right)\\
&\quad\leq c\lambda^{\frac{1}{q}+\frac{d}{\alpha p}+\frac{1}{\alpha}-1}\|f\|_{\mL^q_p(T)}\left(T(\lambda+2+4\nu(B^c_R))
+\mE\left(\int_{t_0\wedge\tau}^{t_1\wedge\tau}|\xi(s)|\dif s\Big|{\sF_{t_0\wedge\tau}}\right)\right),
\end{align*}
which yields the estimate by letting $\lambda$ be large enough since $\frac{d}{p}+\frac{\alpha}{q}<\alpha-1$.
\end{proof}

That the above Krylov estimate required $\frac{d}{p}+\frac{\alpha}{q}<\alpha-1$ is too strong for later use. Below we relax it to
$\frac{d}{p}+\frac{\alpha}{q}<\alpha$, which is similar to Theorem \ref{enou}.
\bt\label{Th65}
Let $T>0$ and $p_1,q_1\in(1,\infty)$ with $\frac{d}{p_1}+\frac{\alpha}{q_1}<\alpha-1$
and $b\in\mL^{q_1}_{p_1}(T)$. Suppose that {\bf (H$^g$)} holds and $X_t$ satisfies
$$
X_t=X_0+\int^t_0\!\!\!\int_{|z|<R}g_s(X_{s-},z)\tilde N(\dif s,\dif z)+\int^t_0\!\!\!\int_{|z|\geq R}\eta_s(z)N(\dif s,\dif z)+\int^t_0b_s(X_s)\dif s,
$$
where $\eta:\mR_+\times\Omega\times\mR^d\to\mR^d$ is a predictable process.
Then for any $p,q\in(1,\infty)$ with $\frac{d}{p}+\frac{\alpha}{q}<\alpha$, the Krylov estimate hold for $X$ with index $p,q$.
\et
\begin{proof}
First of all, we show that for all $p,q\in[1,\infty]$ with $\frac{d}{p}+\frac{\alpha}{q}<\alpha-1$,
\begin{align}\label{UU3}
\mE\left(\int_{t_0}^{t_1}f(s,X_s)\dif s\Big|{\sF_{t_0}}\right)\leq c\|f\|_{\mL^q_p(T)},\ 0\leq t_0<t_1\leq T.
\end{align}
For $n>0$, define
$$
\tau_n:=\inf\left\{t\geq 0: \int^t_0|b_s|(X_s)\dif s\geq n\right\}.
$$
In \eqref{UU4}, if we take $f=|b|$, $\xi(s)=|b_s|(X_s)$ and $\delta=\frac{1}{2}\|b\|_{\mL^{q_1}_{p_1}(T)}$, then
$$
\mE\left(\int_{t_0\wedge\tau_n}^{t_1\wedge\tau_n}|b_s|(X_s)\dif s\Big|{\sF_{t_0\wedge\tau_n}}\right)\leq c\|b\|_{\mL^{q_1}_{p_1}(T)}.
$$
Letting $n\to\infty$, we further have
$$
\mE\left(\int_{t_0}^{t_1}|b_s|(X_s)\dif s\Big|{\sF_{t_0}}\right)\leq c\|b\|_{\mL^{q_1}_{p_1}(T)}.
$$
Substituting this into \eqref{UU4} with $\tau=T$, we get \eqref{UU3}.

Below, without loss of generality, we assume $f\in C^\infty_0(\mR_+\times\mR^d)$.
Let $b^n:=b*\phi_n$ be defined as in \eqref{CON}.
Since $b_n\in\mH^{1,\infty}_\infty(T)$, by Lemma \ref{Le57}, Theorem \ref{Th62} and Remark \ref{Rem44}, for  $\eps$ small enough,
there exists a unique $u_n\in\mH^{\alpha+\eps,\infty}_\infty(T)$ solve the following equation
\begin{align}\label{UU2}
\p_t u_n+\sL^g_\nu u_n+\sL^{b^n}_1 u_n=f.
\end{align}
By It\^o's formula and \eqref{UU2}, we have
\begin{align*}
&\mE\Big( u_n(t_1,X_{t_1})|{\sF_{t_0}}\Big)-u_n(t_0,X_{t_0})\\
&=\mE\left(\int_{t_0}^{t_1}\!\!\Big(\p_s u_n+\sL^g_{\nu,R}u_n+\bar\sL^\eta_{\nu, R}u_n+\sL^b_1 u_n\Big)(s,X_s)\dif s\Big|{\sF_{t_0}}\right)\\
&=\mE\left(\int_{t_0}^{t_1}\!\!\Big(f-\bar\sL^g_{\nu, R}u_n+\bar\sL^\eta_{\nu, R}u_n+(b-b^n)\cdot\nabla u_n\Big)(s,X_s)\dif s\Big|{\sF_{t_0}}\right),
\end{align*}
where $\bar\sL^\eta_{\nu, R}$ is defined as in \eqref{UY7} in terms of $\eta$.
Hence, by \eqref{UU3} and \eqref{NR1} with $\vartheta=0,1$ and $p'=q'=\infty$, we have
\begin{align*}
\mE\left(\int_{t_0}^{t_1}f(s,X_s)\dif s\Big|{\sF_{t_0}}\right)
&\leq(t_1-t_0)(2+4\nu(B^c_R))\|u_n\|_{\mL^\infty(t_0, t_1)}+c\|\nabla u_n\|_{\mL^\infty(T)}\|b-b^n\|_{\mL^{q_1}_{p_1}(T)}\\
&\leq c\|f\|_{\mL^q_p(t_0,t_1)}+c\|f\|_{\mL^{q_1}_{p_1}(t_0,t_1)}\|b-b^n\|_{\mL^{q_1}_{p_1}(T)},
\end{align*}
where $c$ is independent of $n$ due to $\|b_n\|_{\mL^{q_1}_{p_1}(T)}\leq \|b\|_{\mL^{q_1}_{p_1}(T)}$.
Letting $n\to\infty$, we obtain the desired estimate.
\end{proof}

\section{Strong well-posedness of SDEs with jumps}

\subsection{Proof of Theorem \ref{main1}}

Below we fix $T>0$ and assume that {\bf (H$^\sigma$)} holds
and for some $p,q\in(2,\infty)$ with $\frac{d}{p}+\frac{2}{q}<1$,
$$
|\nabla\sigma|,\ b,\ (\Gamma^{1,2}_{0,R}(g))^{1/2}\in \mL^{q}_p(T),
$$
and
$$
\Gamma^{0,2}_{0,R}(g)\in\mL^\infty(T),\ \lim_{\eps\to 0}\|\Gamma^{0,2}_{0,\eps}(g)\|_{\mL^\infty(T)}=0,
$$
where $\Gamma^{j,2}_{0,R}(g)$ is defined by \eqref{LK1}.

Consider the following backward second order partial integral differential equation:
\begin{align}
\p_tu+(\sL^a_2-\lambda)u+\sL^b_1u+\sL^g_{\nu,R}u+b=0,\quad u(T)=0. \label{pde2}
\end{align}
Since $\frac{d}{p}+\frac{2}{q}<1$, by Theorem \ref{pde}, for $\lambda$ large enough, there is a unique solution $u\in\mH^{2,q}_{p}(T)$ to the  above equation with
$$
\|u\|_{\mL^\infty(T)}+\|\nabla u\|_{\mL^\infty(T)}\leq \tfrac{1}{2}.
$$
Let $u_\infty(t,x):=u(t,x)$ and $u_n$ be defined as in \eqref{CON}. Define for $n\in\mN\cup\{\infty\}$,
\begin{align*}
\Phi_n(t,x):=x+u_n(t,x).
\end{align*}
Since for each $t\in[0,T]$,
$$
\tfrac{1}{2}|x-y|\leq\big|\Phi_n(t,x)-\Phi_n(t,y)\big|\leq \tfrac{3}{2}|x-y|,
$$
the map $x\rightarrow\Phi_n(t,x)$ forms a $C^1$-diffeomorphism and
\begin{align}
1/2\leq \|\nabla\Phi_n\|_{\mL^\infty(T)},\ \|\nabla\Phi^{-1}_n\|_{\mL^\infty(T)}\leq 2,   \label{upd}
\end{align}
where $\Phi^{-1}_n(t,\cdot)$ is the inverse of $\Phi_n(t,\cdot)$ and
$$
\Phi^{-1}_n(t,y)=y-u_n(t,\Phi^{-1}_n(t,y)).
$$
The following limits are easily verified by the definition, $u\in\mH^{2,q}_{p}(T)$ and \eqref{upd}:
\begin{align}\label{UY9}
\begin{split}
&\lim_{n\to\infty}\|\nabla^j\Phi_n-\nabla^j\Phi_\infty\|_{\mL^\infty(T)}=0,\ \ \lim_{n\to\infty}\|\nabla^j\Phi^{-1}_n-\nabla^j\Phi^{-1}_\infty\|_{\mL^\infty(T)}=0,\ j=0,1,\\
&\qquad\lim_{n\to\infty}\|(\Phi_n-\Phi_\infty)\chi_m\|_{\mH^{2,q}_{p}(T)}=0,\ \ \lim_{n\to\infty}\|(\Phi^{-1}_n-\Phi^{-1}_\infty)\chi_m\|_{\mH^{2,q}_{p}(T)}=0,
\end{split}
\end{align}
where $\chi_m$ is defined by \eqref{Cut}.
Now let us define $\Phi_t(x):=\Phi_\infty(t,x)$ and
\begin{align}\label{c1}
\begin{split}
&\tilde\sigma_t(y):=\big(\nabla\Phi_t\cdot \sigma_t\big)\circ\Phi_t^{-1}(y),\quad \tilde b_t(y):=\lambda u\big(t,\Phi_t^{-1}(y)\big),  \\
&\qquad\tilde g_t(y,z):=\Phi_t\Big(\Phi_t^{-1}(y)+g_t\big(\Phi_t^{-1}(y),z\big)\Big)-y.
\end{split}
\end{align}
We have
\bp\label{BP}
\begin{enumerate}[(i)]
\item $\tilde\sigma$ satisfies {\bf (H$^\sigma$)} and $|\nabla\tilde\sigma|\in \mL^{q}_p(T)$, $\tilde b\in\mH^{1,\infty}_\infty(T)$ and
$$
(\Gamma^{1,2}_{0,R}(\tilde g))^{1/2}\in \mL^{q}_p(T),\ \
\Gamma^{0,2}_{0,R}(\tilde g)\in\mL^\infty(T),\ \lim_{\eps\to 0}\|\Gamma^{0,2}_{0,\eps}(\tilde g)\|_{\mL^\infty(T)}=0.
$$
\item $\lim_{n\to\infty}\big\|(\p_s+\sL^\sigma_2+\sL^b_1+\sL^g_{\nu, R})\Phi_n-\lambda u\big\|_{\mL^q_p(T)}=0$.

\item $\lim_{n\to\infty}\big\|\big((\p_s+\sL^{\tilde\sigma}_2+\sL^{\tilde b}_1+\sL^{\tilde g}_{\nu, R})\Phi^{-1}_n-b\circ\Phi^{-1}\big)\chi_m\big\|_{\mL^q_p(T)}=0$,
where $\chi_m$ is defined by \eqref{Cut}.
\end{enumerate}
\ep
\begin{proof}
(i) We only show $(\Gamma^{1,2}_{0,R}(\tilde g))^{1/2}\in \mL^{q}_p(T)$. The others are direct by definition.
Let $\bar g_t(x,z):=\Phi_t\big(x+g_t(x,z)\big)-\Phi_t(x)$. By \eqref{upd}, it suffices to show $(\Gamma^{1,2}_{0,R}(\bar g))^{1/2}\in \mL^{q}_p(T)$. Noticing that
\begin{align*}
|\nabla_x\bar g_t(x,z)|&=\big|(\nabla\Phi_t)\big(x+g_t(x,z)\big)\cdot\big(\mI+\nabla_xg_t(x,z)\big)-\nabla\Phi_t(x)\big|\\
&\leq\sup_y|y|^{-1}|(\nabla\Phi_t)\big(x+y\big)-\nabla\Phi_t(x)|\cdot|g_t(x,z)|+2|\nabla_xg_t(x,z)|,
\end{align*}
in view of $p>d$, by Lemma \ref{impo}, we have
\begin{align*}
\big\|(\Gamma^{1,2}_{0,R}(\bar g))^{1/2}\big\|_{\mL^q_p(T)}&\leq\Big\|\sup_y|y|^{-1}|(\nabla\Phi)\big(\cdot+y\big)-\nabla\Phi(\cdot)|\Big\|_{\mL^q_p(T)}\\
&\times\big\|(\Gamma^{0,2}_{0,R}(g))^{1/2}\big\|_{\mL^\infty(T)}+2\big\|(\Gamma^{1,2}_{0,R}(g))^{1/2}\big\|_{\mL^q_p(T)}\\
&\lesssim\|\nabla u\|_{\mH^{1,q}_p(T)}\big\|\Gamma^{0,2}_{0,R}(g)\big\|^{1/2}_{\mL^\infty(T)}+2\big\|(\Gamma^{1,2}_{0,R}(g))^{1/2}\big\|_{\mL^q_p(T)}.
\end{align*}
(ii) By \eqref{pde2} and the same calculations as in \eqref{BB2}, one sees that
\begin{align*}
&\|(\p_s+\sL^\sigma_2+\sL^b_1+\sL^g_{\nu, R})\Phi_n-\lambda u\|_{\mL^q_p(T)}\\
&=\|(\p_s+\sL^\sigma_2+\sL^b_1+\sL^g_{\nu, R})(u_n-u)\|_{\mL^q_p(T)}\\
&\lesssim\|u_n-u\|_{\mH^{2,q}_p(T)}\to 0,\ \ n\to\infty.
\end{align*}
(iii) For simplicity, we drop the time variable and set
\begin{align*}
&\tilde\sigma_n(y):=\big(\nabla\Phi_n\cdot \sigma\big)\circ\Phi^{-1}_n(y),\quad \tilde b_n(y):=\lambda u\big(\Phi^{-1}_n(y)\big),  \\
&\qquad\tilde g_n(y,z):=\Phi_n\Big(\Phi^{-1}_n(y)+g_t\big(\Phi^{-1}_n(y),z\big)\Big)-y.
\end{align*}
By cumbersome calculations, we have
$$
(\p_s+\sL^{\tilde\sigma_n}_2+\sL^{\tilde b_n}_1+\sL^{\tilde g_n}_{\nu, R})\Phi^{-1}_n=b\circ\Phi^{-1}_n.
$$
The limit in (iii) now follows by \eqref{UY9}.
\end{proof}

Now, by Propositions \ref{Pr313} and \ref{BP}, we have

\bl\label{zvon}
Let $\Phi_t(x)$ be defined as above. Then $X_t$ solves SDE \eqref{sde1}
if and only if $Y_t:=\Phi_t(X_t)$ solves the following SDE:
\begin{align}\label{sde2}
\dif Y_t=\tilde\sigma_t(Y_t)\dif W_t+\tilde b_t(Y_t)\dif t+\int_{|z|<R}\tilde g_t(Y_{t-},z)\tilde{N}(\dif t, \dif z)
+\int_{|z|\geq R}\tilde g_t(Y_{t-},z){N}(\dif t, \dif z),
\end{align}
where $\tilde\sigma,\tilde b$ and $\tilde g$ are defined by \eqref{c1}.
\el

Now we can give

\begin{proof}[Proof of Theorem \ref{main1}]

By Lemma \ref{zvon},
it suffices to prove Theorem \ref{main1} for SDE \eqref{sde2}. For the sake of simplicity, we shall drop the tilde over $\tilde\sigma,\tilde b$ and $\tilde g$.

(i) Define
$$
\sigma^{(n)}_t(y):=\sigma_t*\phi_n(y),\ \  g^{(n)}_t(y,z):= g_t(\cdot,z)*\phi_n(y),
$$
where $\phi_n$ is the mollifiers in $\mR^d$.
Since $\sigma$ satisfies {\bf (H$^\sigma$)},
there is a $n_0$ large enough such that for all $n\geq n_0$,
$$
\mbox{$\sigma^{(n)}$ satisfies {\bf (H$^\sigma$)} uniformly with respect to $n$},
$$
and
\begin{align*}
&\|\nabla\sigma^{(n)}\|_{\mL^q_p(T)}\leq\|\nabla\sigma\|_{\mL^q_p(T)},\
\|(\Gamma^{1,2}_{0,R}( g^{(n)}))^{1/2}\|_{\mL^q_p(T)}\leq \|(\Gamma^{1,2}_{0,R}( g))^{1/2}\|_{\mL^q_p(T)},\\
&\|\Gamma^{0,2}_{0,R}( g^{(n)})\|_{\mL^\infty(T)}\leq \|\Gamma^{0,2}_{0,R}( g)\|_{\mL^\infty(T)},
\ \ \lim_{\eps\to 0}\sup_n\|\Gamma^{0,2}_{0,\eps}( g^{(n)})\|_{\mL^\infty(T)}=0.
\end{align*}
Let $Y^{(n)}$ solve the following SDE with no big jumps:
\begin{align}\label{HF7}
Y^{(n)}_t=y+\int_0^t\!\!\sigma^{(n)}_s(Y^{(n)}_s)\dif W_s+\int_0^t b_s(Y^{(n)}_s)\dif s
+\int_0^t\!\!\!\int_{|z|<R} g^{(n)}_s(Y^{(n)}_{s-},z)\tilde{N}(\dif s, \dif z).
\end{align}
Since $\sigma^{(n)},  g^{(n)}$ satisfy the assumptions of Theorem \ref{enou} uniformly with respect to $n$, by Theorem \ref{enou},
$Y^{(n)}$ satisfies the Krylov estimate for all $p',q'$ with $\frac{d}{p'}+\frac{2}{q'}<2$ and the Krylov constant is independent of $n$.
Thus, by Theorem \ref{Sta} with $r=1$, we have for any $\theta\in(0,1)$,
$$
\mE\left(\sup_{t\in[0,T]}|Y^{(n)}_t-Y^{(m)}_t|^{2\theta}\right)
\lesssim\|\sigma^{(n)}-\sigma^{(m)}\|^{2\theta}_{\mL^2_\infty(T)}
+\big\|\Gamma^{0,2}_{0,R}( g^{(n)}- g^{(m)})\big\|^\theta_{\mL^1_\infty(T)}.
$$
Since $p>d$, by \eqref{w12} with $\mB=L^2(B_R;\nu)$ and $\mB=\mR^d\otimes\mR^d$ respectively, we have
\begin{align}
\Gamma^{0,2}_{0,R}( g^{(n)}_t- g_t)(y)&=\int_{|z|<R}\left|\int_{\mR^d}(g_t(y-y',z)-g_t(y,z))\phi_n(y')\dif y'\right|^2\nu(\dif z)\no\\
&\leq\left(\int_{\mR^d}\|g(y-y',\cdot)-g(y,\cdot)\|_{L^2(B_R;\nu)}\phi_n(y')\dif y'\right)^2\no\\
&\lesssim\left(\int_{\mR^d}|y'|^{1-d/p}\phi_n(y')\dif y'\right)^2\|(\Gamma^{1,2}_{0,R}(g_t))^{1/2}\|^2_p\no\\
&\leq n^{-2+2d/p}\|(\Gamma^{1,2}_{0,R}(g_t))^{1/2}\|^2_p,\label{EG3}
\end{align}
and
$$
\|\sigma^{(n)}_t(y)- \sigma_t(y)\|\lesssim n^{-1+d/p}\|\nabla\sigma_t\|_p.
$$
Therefore,
$$
\mE\left(\sup_{t\in[0,T]}|Y^{(n)}_t-Y^{(m)}_t|^{2\theta}\right)\lesssim (n^{-1+d/p}+m^{-1+d/p})^{2\theta}\to 0 \mbox{ as $n,m\to\infty$},
$$
and there exists a c\`adl\`ag $\sF_t$-adapted process $Y$ such that
$$
\lim_{n\to\infty}\mE\left(\sup_{t\in[0,T]}|Y^{(n)}_t-Y_t|^{2\theta}\right)=0.
$$
By Remark \ref{Rem23}, $Y_t$ also satisfies the Krylov estimate with index $p,q$.
By taking limits for \eqref{HF7}, one finds that $Y_t$ solves
\begin{align}\label{SDE1}
Y_t=y+\int_0^t\!\!\sigma_s(Y_s)\dif W_s+\int_0^t\!\! b_s(Y_s)\dif s
+\int_0^t\!\!\!\int_{|z|<R}\!\!\! g_s(Y_{s-},z)\tilde{N}(\dif s, \dif z).
\end{align}
For example,  letting $Y^{\infty}_t:=Y_t$, by \eqref{EG3}, we have
\begin{align*}
&\sup_{n\in\mN\cup\{\infty\}}\mE\left|\int_0^t\!\!\!\int_{|z|<R}(g^{(m)}_s(Y^{(n)}_{s-},z)-g_s(Y^{(n)}_{s-},z))\tilde{N}(\dif s, \dif z)\right|^2\\
&=\sup_{n\in\mN\cup\{\infty\}}\mE\int_0^t\!\!\!\int_{|z|<R}| g^{(m)}_s(Y^{(n)}_{s-},z)-
 g_s(Y^{(n)}_{s-},z)|^2\nu(\dif z)\dif s\\
&\lesssim \|\Gamma^{0,2}_{0,R}(g^{(m)}- g)\|_{\mL^1_\infty(T)}\to 0,\ m\to\infty,
\end{align*}
and for each $m\in\mN$,
$$
\lim_{n\to\infty}\mE\left|\int_0^t\!\!\!\int_{|z|<R}\!\!\! g^{(m)}_s(Y^{(n)}_{s-},z)\tilde{N}(\dif s, \dif z)
-\int_0^t\!\!\!\int_{|z|<R}\!\!\! g^{(m)}_s(Y_{s-},z)\tilde{N}(\dif s, \dif z)\right|^2=0.
$$
Combining the above two estimates, we obtain
$$
\lim_{n\to\infty}\mE\left|\int_0^t\!\!\!\int_{|z|<R}\!\!\! g^{(n)}_s(Y^{(n)}_{s-},z)\tilde{N}(\dif s, \dif z)
-\int_0^t\!\!\!\int_{|z|<R}\!\!\! g_s(Y_{s-},z)\tilde{N}(\dif s, \dif z)\right|^2=0.
$$

(ii) To show \eqref{feller}, we first consider SDE \eqref{SDE1}. By the classical Bismut-Elworthy-Li's formula (see \cite{W-X-Zh}),  we have for any $h\in\mR^d$,
\begin{align}
\nabla_h\mE\varphi\big(Y^{(n)}_t(y)\big)=\frac{1}{t}\mE\Bigg[\varphi\big(Y^{(n)}_t(y)\big)
\int_0^t\big[\sigma\big(Y^{(n)}_s(y)\big)\big]^{-1}\nabla_h Y^{(n)}_s(y)\dif W_s\Bigg],  \label{grad}
\end{align}
where $\nabla_h Y_t^{(n)}(y):=\lim_{\eps\to 0}[Y_t^{(n)}(y+\eps h)-Y_t^{(n)}(y)]/\eps$ is the derivative flow of $Y^{(n)}_t(y)$ with respect to the initial value $y$.

Now by Theorem \ref{Sta}, we have for any $\theta\in(0,1)$,
$$
\mE|Y_t^{(n)}(y)-Y_t^{(n)}(y')|^{2\theta}\leq c|y-y'|^{2\theta},
$$
where $c$ is independent of $n$. Let $\theta\in(1/2,1)$. By Theorem \ref{th1} with $p=2\theta$ and $q=r=\infty$, we get
$$
\sup_n\sup_y\mE\left(\sup_{t\in[0,T]}|\nabla Y_t^{(n)}(y)|^{2\theta}\right)\leq c.
$$
Hence, by \eqref{grad} and Burkholder's inequality,
\begin{align*}
\sup_y|\nabla \mE\varphi\big(Y^{(n)}_t(y)\big)|&\leq \frac{\|\varphi\|_\infty\|\sigma^{-1}\|_\infty}{t}\sup_y\mE\Bigg[
\int_0^t|\nabla Y^{(n)}_s(y)|^2\dif s\Bigg]^{1/2}\\
&\leq \frac{\|\varphi\|_\infty\|\sigma^{-1}\|_\infty}{\sqrt{t}}\sup_y\mE\Bigg[
\sup_{s\in[0,t]}|\nabla Y^{(n)}_s(y)|\Bigg]\leq c\|\varphi\|_\infty t^{-1/2},
\end{align*}
which means that
$$
|\mE\varphi\big(Y^{(n)}_t(y)\big)-\mE\varphi\big(Y^{(n)}_t(y')\big)|\leq  c_T\|\varphi\|_\infty t^{-1/2}|y-y'|.
$$
By taking limits $n\to\infty$ we get
$$
\mbox{Var}(P_t(y,\cdot)-P_t(y',\cdot))=\sup_{\varphi\in C_b(\mR^d), \|\varphi\|_\infty\leq 1}|\mE\varphi\big(Y_t(y)\big)-\mE\varphi\big(Y_t(y')\big)|\leq  c_T t^{-1/2}|y-y'|,
$$
where $P_t(y,\cdot)$ denotes the law of $Y_t(y)$.

(iii) To allow the large jump in the equation, we shall use the interlacing technique. More precisely, let ${\rm p}_s$ be a point function
on $\mR_+$ with values in $B^c_R$, $\mu$ the associated counting measure, i.e.,
$$
\mu([0,t],A):=\#\{{\rm p}_s\in A: s\in[0,t]\},\ \ A\in\sB(B^c_R).
$$
Let $\tau^{\rm p}_n:=\inf\{t>0: \mu([0,t]; B^c_R)=n\}$ be the $n$-th jump time of $t\mapsto\mu([0,t]; B^c_R)$. Let $Y_{s,t}(y)$
solve SDE \eqref{sde2} with initial value $Y_{s,s}(y)=y$. Define $Y^{\rm p}_t(y)$ recursively by
$$
Y^{\rm p}_t:=Y^{\rm p}_t(y):=
\left\{
\begin{aligned}
&Y_{\tau^{\rm p}_{n-1},t}(Y^{\rm p}_{\tau^{\rm p}_{n-1}}(y)),\qquad \qquad t\in[\tau^{\rm p}_{n-1},\tau^{\rm p}_n),\\
&Y^{\rm p}_{\tau^{\rm p}_{n}-}(y)+g_{\tau^{\rm p}_{n}}(Y^{\rm p}_{\tau^{\rm p}_n-}(y),{\rm p}_{\tau^{\rm p}_n}),\ \ t=\tau^{\rm p}_n.
\end{aligned}
\right.
$$
It is easy to see that $Y^{\rm p}_t$ solves the following SDE with starting point $Y^{\rm p}_0=y$:
\begin{align*}
\dif Y^{\rm p}_t=\sigma_t(Y^{\rm p}_t)\dif W_t+b_t(Y^{\rm p}_t)\dif t+\!\!\int_{|z|<R}g_t(Y^{\rm p}_{t-},z)\tilde N(\dif t,\dif z)
+\!\!\int_{|z|\geq R}g_t(Y^{\rm p}_{t-},z)\mu(\dif t,\dif z).
\end{align*}
In particular, if we let ${\rm p}^N_s$ be the Poisson point process with values in $B^c_R$ associated to the Poisson random measure $N(\dif t,\dif z)$, i.e.,
$$
N([0,t], A)=\#\{{\rm p}^N_s\in A: s\in[0,t]\},\ \ A\in\sB(B^c_R),
$$
then ${\rm p}^N$ is independent with $X$. Therefore, $\tilde Y_t:=Y^{{\rm p}^N}_t$ solves SDE
\eqref{sde1} with $Y_0=y$.

Next we show \eqref{feller} for $\tilde Y_t(y)$. We adopt the same argument as used in \cite{W-X-Zh}.
We first look at it for $Y^{\rm p}_t(y)$.
Observing that
$$
Y^{\rm p}_t(y)=
\left\{
\begin{aligned}
&Y_t(y),\ \ t<\tau^{\rm p}_1,\\
&Y_{\tau^{\rm p}_1-}(y)+g_{\tau^{\rm p}_1}(Y^{\rm p}_{\tau^{\rm p}_1-}(y),{\rm p}_{\tau^{\rm p}_1}),\ \ t=\tau^{\rm p}_1,\\
&Y_{\tau^{\rm p}_{1},t}(Y^{\rm p}_{\tau^{\rm p}_{1}}(y)),\ \ t\in[\tau^{\rm p}_1,\tau^{\rm p}_2),\\
&\cdots\cdots,
\end{aligned}
\right.
$$
by what we have proved in step (ii), and since $Y_{s,t}(\cdot)$ and $Y_{0,s}(\cdot)$ are independent, one sees that
$$
|\mE\varphi\big(Y^{\rm p}_t(y)\big)-\mE\varphi\big(Y^{\rm p}_t(y')\big)|\leq  c_T\|\varphi\|_\infty (t\wedge\tau^{\rm p}_1)^{-1/2}|y-y'|,\ \ t\in[0,T].
$$
Hence,
\begin{align*}
|\mE\varphi\big(\tilde Y_t(y)\big)-\mE\varphi\big(\tilde Y_t(y')\big)|
\leq c_T\|\varphi\|_\infty \mE \Big(t\wedge\tau^{{\rm p}^N}_1\Big)^{-1/2}|y-y'|,\ \ t\in[0,T].
\end{align*}
Since the random variable $\tau^{{\rm p}^N}_1=\inf\{t>0: N([0,t]; B^c_R)=1\}$ obeys the exponential distribution with parameter $\nu(B^c_R)$,
by easy calculations, we have
$$
\mE \Big(t\wedge\tau^{{\rm p}^N}_1\Big)^{-1/2}\leq ct^{-1/2}.
$$
Thus, we get \eqref{feller} for $\tilde Y_t(y)$. The proof is complete.
\end{proof}

\subsection{Proof of Theorem \ref{main3}} Let $T>0$
and $\nu(\dif z)=|z|^{-d-\alpha}\dif z$ for some $\alpha\in(1,2)$. Below, we assume {\bf (H$^g$)}, and
for some $\theta\in(1-\frac{\alpha}{2},1)$,
$p\in(\frac{2d}{\alpha}\vee 2,\infty)$ and $q\in(\frac{2\alpha}{\alpha+2(\theta-1)},\infty)$ with $\frac{d}{p}+\frac{\alpha}{q}<\frac{\alpha}{2}$,
$$
\big(\Gamma^{1,2}_{0,R}(g)\big)^{1/2}\in \mL^{q}_p(T),\ b\in \mH^{\theta,q}_p(T).
$$
We also fix
$$
\gamma\in((d/p+1-\theta)\vee1\vee (1+\alpha/2-\theta), \alpha-\alpha/q).
$$
Consider the following backward nonlocal equation:
\begin{align*}
\p_tu+(\sL^g_\nu-\lambda)u+\sL^b_1u+b=0,\quad u(T)=0.
\end{align*}
For any $\lambda\geq 0$, by Theorem \ref{Th63}, there is a unique solution $u\in\mH^{\gamma+\theta,q}_{p}(T)$ to the  above equation.  Moreover,
by Sobolev's embedding \eqref{Sob},
\begin{align*}
b\in \mL^{q}_{p_1}(T), \mbox{ where } p_1:=dp/(d-\theta p)\mbox{ satisfies } \tfrac{d}{p_1}+\tfrac{\alpha}{q}<\tfrac{\alpha}{2}-\theta<\alpha-1,
\end{align*}
and using \eqref{NR1} with $p'=q'=\infty$ and $\vartheta=\gamma$, for $\lambda$ large enough, we have
\begin{align}\label{EW1}
\sup_{t\in[0,T]}\|u(t)\|_{C^\gamma_b}\leq \|u\|_{\mH^{\gamma,\infty}_\infty(T)}\leq \tfrac{1}{2}.
\end{align}
Let $u_\infty=u$ and $u_n=u*\phi_n$ for $n\in\mN$.
Define for $n\in\mN\cup\{\infty\}$,
$$
\Phi_n(t,x):=x+u_n(t,x).
$$
Since for each $t\in[0,T]$, \begin{align}\label{KL1}
\tfrac{1}{2}|x-y|\leq\big|\Phi_n(t,x)-\Phi_n(t,y)\big|\leq \tfrac{3}{2}|x-y|,
\end{align}
the map $x\rightarrow\Phi_n(t,x)$ forms a $C^1$-diffeomorphism and
\begin{align}
1/2\leq \|\nabla\Phi_n\|_{\mL^\infty(T)},\ \|\nabla\Phi^{-1}_n\|_{\mL^\infty(T)}\leq 2,   \label{upd0}
\end{align}
where $\Phi^{-1}_n(t,\cdot)$ is the inverse of $\Phi_n(t,\cdot)$.  Moreover, by \eqref{EW1}, we also have
\begin{align}\label{EW2}
\nabla\Phi_n,\nabla\Phi^{-1}_n\mbox{ are H\"older continuous uniformly with respect to $t, n$.}
\end{align}
As above, define $\Phi_t(x):=\Phi_\infty(t,x)$ and
\begin{align}\label{c11}
\begin{split}
\tilde b_t(y):=\lambda u\big(t,\Phi_t^{-1}(y)\big), \quad\tilde g_t(y,z):=\Phi_t\Big(\Phi_t^{-1}(y)+g_t\big(\Phi_t^{-1}(y),z\big)\Big)-y.
\end{split}
\end{align}
\bp\label{Pr63}
\begin{enumerate}[(i)]
\item $\tilde b\in\mH^{1,\infty}_\infty(T)$ and $\tilde g$ satisfies {\bf (H$^{ g}$)}, $(\Gamma^{1,2}_{0,R}(\tilde g))^{1/2}\in\mL^q_p(T).$

\item $\lim_{n\to\infty}\big\|(\p_s+\sL^b_1+\sL^g_{\nu})\Phi_n-\lambda u\big\|_{\mL^q_p(T)}=0$.

\item $\lim_{n\to\infty}\big\|\big((\p_s+\sL^{\tilde b}_1+\sL^{\tilde g}_{\nu})\Phi^{-1}_n-b\circ\Phi^{-1}\big)\chi_m\big\|_{\mL^q_p(T)}=0$,
where $\chi_m$ is defined by \eqref{Cut}.
\end{enumerate}
\ep
\begin{proof}
(i) It is clear that $\tilde b\in\mH^{1,\infty}_\infty(T)$ by definition.
Let $\bar g_t(x,z):=\Phi_t\big(x+g_t(x,z)\big)-\Phi_t(x)$. We show that
$\bar g$ satisfies {\bf (H$^{ g}$)} and $(\Gamma^{1,2}_{0,R}(\bar g))^{1/2}\in \mL^{q}_p(T)$.
Clearly, $\bar g_t(x,0)=0$ by $g_t(x,0)=0$, and by \eqref{KL1},
$$
(2c_1)^{-1}|z-z'|\leq|\bar g_t(x,z)-\bar g_t(x,z')|\leq 2c_1|z-z'|.
$$
Moreover, notice that
$$
\nabla_z \bar g_t(x,z)=(\nabla\Phi_t)\big(x+g_t(x,z)\big)\cdot\nabla_z g_t(x,z).
$$
Since $g$ satisfies {\bf (H$^{ g}$)}, by \eqref{upd0} and \eqref{EW2}, it is easy to see that $\bar g$ also satisfies {\bf (H$^{ g}$)},
and so does $\tilde g$.

On the other hand, define
$$
U_t(x):=\sup_{y}|y|^{1-\theta-\gamma}|\nabla u_t(x+y)-\nabla u_t(x)|.
$$
Notice that
\begin{align*}
|\nabla_x \bar g_t(x,z)|&=|\nabla\Phi_t(x+g_t(x,z))(\mI+\nabla g_t(x,z))-\nabla\Phi_t(x)|\\
&\leq|\nabla u_t(x+g(x,z))-\nabla u_t(x)|+2|\nabla_x g_t(x,z)|\\
&\leq U_t(x) |g(x,z)|^{\theta+\gamma-1} +2|\nabla_x g_t(x,z)|.
\end{align*}
Hence, by {\bf (H$^g$)}, we have
\begin{align*}
\Gamma^{1,2}_{0,R}(\bar g)(x)=\int_{|z|<R}|\nabla\bar g_t(x,z)|^2\nu(\dif z)\lesssim
U^2(x)\int_{|z|<R}|z|^{2(\theta+\gamma-1)-d-\alpha}\dif z+\Gamma^{1,2}_{0,R}(g)(x).
\end{align*}
By Lemma \ref{impo}, since $p(\theta+\gamma-1)>d$, we have
\begin{align*}
\|U\|_{\mL^q_p(T)}\lesssim\|\nabla u\|_{\mH^{\theta+\gamma-1,q}_p(T)}\leq\|u\|_{\mH^{\theta+\gamma,q}_p(T)}.
\end{align*}
Since $\theta+\gamma-1>\frac{\alpha}{2}$, we get $(\Gamma^{1,2}_{0,1}(\bar g))^{1/2}\in\mL^q_p(T)$ and so $(\Gamma^{1,2}_{0,1}(\tilde g))^{1/2}\in\mL^q_p(T)$.

(ii) and (iii) can be proven by the same calculations as in Proposition \ref{BP}.
\end{proof}

Since $\frac{d}{p_1}+\frac{\alpha}{q}<\alpha-1$ and $b\in\mL^{q}_{p_1}(T)\subset\mK^\alpha_d$, by Theorem \ref{Th65},
any solution $X$  of SDE \eqref{sde1} with $\sigma\equiv0$
satisfies the Krylov estimate for all $p',q'\in(1,\infty)$ with $\frac{d}{p'}+\frac{\alpha}{q'}<\alpha$.
As in Lemma \ref{zvon}, we have

\bl\label{zvon1}
Let $\Phi_t(x)$ be defined as above. Then $X_t$ solves SDE
$$
\dif X_t=b_t(X_t)\dif t
+\int_{|z|<R} g_t(X_{t-},z)\tilde{N}(\dif t, \dif z)+\int_{|z|\geq R} g_t(X_{t-},z){N}(\dif t, \dif z)
$$
if and only if $Y_t:=\Phi_t(X_t)$ solves the following SDE:
\begin{align}\label{Sde02}
\dif Y_t=\tilde b_t(Y_t)\dif t+\int_{|z|<R}\tilde g_t(Y_{t-},z)\tilde{N}(\dif t, \dif z)+\int_{|z|\geq R}\tilde g_t(Y_{t-},z){N}(\dif t, \dif z),
\end{align}
where $\tilde b$ and $\tilde g$ are defined by \eqref{c11}.
\el

Now we can give

\begin{proof}[Proof of Theorem \ref{main3}]
By Lemma \ref{zvon1},
it suffices to prove the theorem for SDE \eqref{Sde02}.

(i) Let $\tilde g^{(n)}_t(y,z):=\tilde g_t(\cdot,z)*\phi_n(y)$. By (i) of Proposition \ref{Pr63}, there is a $n_0$ large enough such that for all $n\geq n_0$,
$$
\mbox{$\tilde g^{(n)}$ satisfies {\bf (H$^g$)} with constant $c_1$ independent of $n$},
$$
and
\begin{align*}
&\|(\Gamma^{1,2}_{0,R}( \tilde g^{(n)}))^{1/2}\|_{\mL^q_p(T)}\leq \|(\Gamma^{1,2}_{0,R}( \tilde g))^{1/2}\|_{\mL^q_p(T)}.
\end{align*}
Let $Y^{(n)}$ satisfy
\begin{align}\label{HF6}
Y^{(n)}_t=y+\int_0^t\tilde b_s(Y^{(n)}_s)\dif s
+\int_0^t\!\!\!\int_{|z|<R} \tilde g^{(n)}_s(Y^{(n)}_{s-},z)\tilde{N}(\dif s, \dif z).
\end{align}
By Theorem \ref{Th65}, for any $p', q'$ with $\frac{d}{p'}+\frac{\alpha}{q'}<\alpha$,
$Y^{(n)}$ satisfies Krylov's estimate with index $p',q'$ and Krylov's constant independent of $n$.
Thus, by Theorem \ref{Sta} with $r=1$, for any $\theta\in(0,1)$, we have
$$
\mE\left(\sup_{t\in[0,T]}|Y^{(n)}_t-Y^{(m)}_t|^\theta\right)
\lesssim\big\|\Gamma^{0,2}_{0,R}(\tilde g^{(n)}-\tilde g^{(m)})\big\|^\theta_{\mL^1_\infty(T)},
$$
which converges to zero as $n,m\to\infty$ by \eqref{EG3}.
Therefore, there exists a c\`adl\`ag $\sF_t$-adapted process $Y$ such that
$$
\lim_{n\to\infty}\mE\left(\sup_{t\in[0,T]}|Y^{(n)}_t-Y_t|^\theta\right)=0,
$$
and by Remark \ref{Rem23}, $Y_t$ also satisfies the Krylov estimate with index $p,q$.
By taking limits for \eqref{HF6}, one finds that $Y_t$ solves
$$
Y_t=y+\int_0^t\tilde b_s(Y_s)\dif s+\int_0^t\!\!\!\int_{|z|<R}\tilde g_s(Y_{s-},z)\tilde{N}(\dif s, \dif z).
$$
The uniqueness follows by Theorem \ref{Sta}.
For the large jump, we use the same technique as used in the proof of Theorem \ref{main1}.

(ii) To show the existence and estimates of the distribution density of $Y_t(y)$, we use Theorem \ref{heat}.
Consider the operator $\sL:=\sL^{\tilde g}_\nu+\sL^{\tilde b}_1=\sL^{\tilde \kappa}_\alpha+\sL^{\bar b^{\tilde g}}_1
+\sL^{\tilde b}_1$, where we have used
\eqref{499}. Since $\tilde b$ and $\bar b^{\tilde g}$ are bounded and H\"older continuous in $x$,
by Theorem \ref{heat}, the operator $\sL$ admits a fundamental solution
$\tilde\rho(s,y;t,y')$ so that the conclusions of Theorem \ref{heat} still holds. In particular,
$\tilde\rho(s,y;t,y')$ is a family of transition probability density functions.
It determines a Feller process
$$
  \Big(\Omega,\sF, (\mP_{s,y})_{(s,y)\in\mR_+\times\mR^d}; (Y_t)_{t\geq 0}\Big),
$$
with the property that
$$
  \mP_{s,y}\big(Y_t=y,\,0\leq t\leq s\big)=1,
$$
and for $r\in[s,t]$ and $E\in\cB(\mR^d)$,
\begin{align*}
  \mE_{s,y}\big(Y_{t}\in E \,|\, Y_{r}\big)=\int_{E}\tilde\rho(r,Y_{r};s,y')\dif y'.
\end{align*}
Moreover, for any $f\in C^2_b(\mR^d)$, it follows from the Markov property of $Y$ that under $\mP_{s,y}$,
with respect to the filtration $\sF_t:=\sigma\{Y_r, r\leq t\}$,
\begin{align*}
  M^f_s:=f(Y_s)-f(X_t)-\int^s_t\sL_r f(Y_r)\dif r\ \mbox{ is a martingale}.
\end{align*}
In other words, $\mP_{s,y}$ solves the martingale problem for $(\sL, C^2_b (\mR^d))$.
On the other hand, by \cite{Ab-Ka} or \cite{Ch-Zh1}, we know that the martingale problem for $\sL$ is well-posed, and by It\^o's formula, any solution of SDE \eqref{sde1}
is a martingale solution of $\sL$. The desired estimates \eqref{GD5} and \eqref{GD4} for $\tilde\rho$ follows by Theorem \ref{heat}.
\end{proof}

\section{Ergodicity of SDEs with jumps}

This section is devoted to the study of  the existence and uniqueness of  invariant probability measures associated with the time-independent SDE \eqref{sde4}.
\subsection{SDEs with dissipative drifts}
Below we assume that $\sigma$ is continuous and for some $r>-1$ and $\kappa_1,\kappa_2,\kappa_3>0$,
\begin{align}\label{EQ3}
2\<x,b(x)\>+\|\sigma(x)\|^2\leq-\kappa_1|x|^{2+r}+\kappa_2,\ \ |b(x)|\leq \kappa_3(1+|x|^{1+r}),
\end{align}
and for any $\eps>0$ and $\lambda\geq R$,
\begin{align}\label{EQ4}
\Gamma^{0,2}_{0,\lambda}(g)(x)+\Gamma^{0,1}_{\lambda,\infty}(g)(x)\leq\eps|x|^{1+r}+c_{\eps,\lambda}.
\end{align}
We first show the non-explosion and some moment estimates of the unique strong solution to SDE \eqref{sde4}.
\bl\label{Le71}
Under \eqref{EQ3} and \eqref{EQ4}, there is no explosion to SDE \eqref{sde4}. Moreover,
for any $\vartheta\in(0,1)$, there is a constant $c>0$ such that for all $t>0$ and $x\in\mR^d$,
\begin{align}\label{BB7}
\int^t_0\mE|X_s(x)|^{1+r}\dif s+\left[\mE\left(\sup_{s\in[0,t]}|X_s(x)|^{\vartheta}\right)\right]^{1/\vartheta}
\leq c(|x|+t+1),
\end{align}
and
\begin{align}\label{BB6}
\mE|X_t(x)|\leq
\left\{
\begin{aligned}
&c \e^{-t/c}|x|+c,\ \ r=0,\\
& c(1+t^{-1/2}),\  \ r>0.
\end{aligned}
\right.
\end{align}
\el

\begin{proof}
Let $h(x):=\sqrt{1+|x|^2}$. By It\^o's formula, we have
$$
\dif h(X_t)=[\sL^\sigma_2h+\sL^b_1 h+\sL^g_{\nu} h](X_t)\dif t+\dif M_t,
$$
where $M_t$ is a local martingale.
Noticing that
$$
\p_i h(x)=x_i(1+|x|^2)^{-1/2}/2
$$
and
$$
\p_i\p_j h(x)=(1+|x|^2)^{-1/2}\delta_{ij}/2-3x_ix_j(1+|x|^2)^{-3/2}/4,
$$
we have
\begin{align}
\sL^\sigma_2 h(x)
+\sL^b_1h(x)&\leq\big(\|\sigma(x)\|^2+2\<x,b(x)\>\big)(1+|x|^2)^{-1/2}/4.\label{EQ1}
\end{align}
On the other hand, observing that
\begin{align*}
&|h(x+y)-h(x)|\leq|y|\int^1_0|\nabla h(x+sy)|\dif s\leq |y|/2,\\
&\qquad h(x+y)-h(x)-y\cdot \nabla h(x)\leq |y|^2/2,
\end{align*}
we have
\begin{align}\label{EQ2}
\begin{split}
\sL^g_\nu h(x)&=\int_{\mR^d}\Big[h\big(x+g(x,z)\big)-h(x)-1_{|z|<R}g(x,z)\cdot\nabla h(x)\Big]\nu(\dif z)\\
&\leq\frac{1}{2}\int_{|z|<R}|g(x,z)|^2\nu(\dif z)+\frac{1}{2}\int_{|z|\geq R}|g(x,z)|\nu(\dif z)\\
&=\Big(\Gamma^{0,2}_{0,R}(g)(x)+\Gamma^{0,1}_{R,\infty}(g)(x)\Big)/2.
\end{split}
\end{align}
By \eqref{EQ1}, \eqref{EQ2} and \eqref{EQ3}, \eqref{EQ4}, there are $c_1,c_2>0$ only depending on $\kappa_i$ such that
\begin{align*}
[\sL^\sigma_2 h+\sL^b_1 h+\sL^g_{\nu} h](x)\leq-c_1(1+|x|^2)^{(1+r)/2}+c_2.
\end{align*}
Hence,
\begin{align*}
\dif h(X_t)\leq -c_1h(X_t)^{1+r}\dif t+c_2\dif t+\dif M_t.
\end{align*}
Letting $\tau_n:=\inf\{t>0: |X_t|\geq n\}$, we have
$$
c_1\mE\left(\int^{t\wedge\tau_n}_0h(X_s)^{1+r}\dif s\right)\leq h(x)+c_2 t,
$$
and by Lemma \ref{gron}, for any $\vartheta\in(0,1)$,
$$
\mE\left(\sup_{s\in[0,t\wedge\tau_n]}h(X_s)^\vartheta\right)\leq c_\vartheta(h(x)+c_2t)^\vartheta,
$$
which yields that $\tau_n\to\infty$ as $n\to\infty$. By taking limits $n\to\infty$ , we then obtain \eqref{BB7}.
Moreover, we also have
\begin{align*}
\dif \mE h(X_{t})/\dif t
\leq
\left\{
\begin{aligned}
&-c_1\mE\left(h(X_t)\right)+c_2,\ \ r=0,\\
&-c_1\left(\mE h(X_t)\right)^{1+r}+c_2,\ \ r>0.
\end{aligned}
\right.
\end{align*}
Solving this differential inequality, we get \eqref{BB6}.
\end{proof}
The following lemma is useful for showing the irreducibility in the non-degenerate diffusion case.
\bl\label{Le72}
For given $x_0\not=y_0\in\mR^d$ and $m\geq 1$, let $Z_t$ solve the following SDE:
\begin{align}
\dif Z_t&=[b(Z_t)-m(Z_t-y_0)/2]\dif t +  \sigma(Z_t)\dif W_t\no\\
&+\int_{|z|<R}  g(Z_{t-},z)\tilde{N}(\dif t, \dif z)+\int_{|z|\geq R}  g(Z_{t-},z)N(\dif t, \dif z),\quad Z_0=x_0.\label{Eq5}
\end{align}
Under \eqref{EQ3} and \eqref{EQ4}, for any $0<a<|x_0-y_0|$ and $T>0$, there exists an $m$ large enough such that
\begin{align}\label{EB4}
\mP(|Z_T(x_0)-y_0|>a)<1/2,
\end{align}
and for any $\vartheta\in(0,1)$,
\begin{align}
\mE\left(\sup_{t\in[0,T]}|Z_t|^{\vartheta}\right)<\infty.   \label{EB5}
\end{align}
\el
\begin{proof}
First of all, by using the same argument as in estimating \eqref{BB7}, we have (\ref{EB5}).
Let us show \eqref{EB4}. For $\lambda>0$, define
$$
\tau_\lambda:=\inf\{t\geq 0: N([0,t],B^c_\lambda)=1\}.
$$
Let $T>0$ be fixed. Since $\tau_\lambda$ obeys the exponential distribution with parameter $\nu(B^c_\lambda)$,
one can choose $\lambda\geq R$ large enough so that
\begin{align}\label{RR}
\mP(\tau_\lambda\leq T)=1-\e^{-T\nu(B^c_\lambda)}\leq 1/4.
\end{align}
For this $\lambda$, let $Z^\lambda_t$ solve the following SDE with starting point $Z^\lambda_0=x_0$,
\begin{align*}
\dif Z^\lambda_t&=[b_\lambda(Z^\lambda_t)-m(Z^\lambda_t-y_0)/2]\dif t +  \sigma(Z^\lambda_t)\dif W_t
+\int_{|z|<\lambda}  g(Z^\lambda_{t-},z)\tilde{N}(\dif t, \dif z),
\end{align*}
where $b_\lambda(x):=b(x)+\int_{R\leq|z|<\lambda}  g(x,z)\nu(\dif z)$.
 Clearly,
\begin{align}\label{RR0}
Z_t=Z^\lambda_t,\ \ t\in[0,\tau_\lambda).
\end{align}
By It\^o's formula and \eqref{EQ3}, \eqref{EQ4}, we have
\begin{align*}
\e^{mt}\mE|Z^\lambda_t-y_0|^2&=|x_0-y_0|^2+\mE\int^t_0\e^{ms}\Big(2\<Z^\lambda_s-y_0,b_\lambda(Z^\lambda_s)\>+\|\sigma(Z^\lambda_s)\|^2\Big)\dif s\\
&\quad+\mE\int^t_0\e^{ms}\int_{|z|\leq\lambda}|g(Z^\lambda_s,z)|^2\nu(\dif z)\dif s\\
&=|x_0-y_0|^2+\mE\int^t_0\e^{ms}\Big(2\<Z^\lambda_s,b(Z^\lambda_s)\>+\|\sigma(Z^\lambda_s)\|^2-2\<y_0,b_\lambda(Z^\lambda_s)\>\Big)\dif s\\
&\quad+\mE\int^t_0\e^{ms}\left(2\Big\<Z^\lambda_s,\int_{R\leq|z|<\lambda}  g(Z^\lambda_s,z)\nu(\dif z)\Big\>
+\Gamma^{0,2}_{0,\lambda}(g)(Z^\lambda_s)\right)\dif s\\
&\leq|x_0-y_0|^2+\mE\int^t_0\e^{ms}\Big(-\kappa_1|Z^\lambda_s|^{2+r}+\kappa_2\Big)\dif s\\
&\quad+2|y_0|\mE\int^t_0\e^{ms}\Big(\kappa_3(|Z^\lambda_s|^{1+r}+1)+\Gamma^{0,1}_{R,\lambda}(g)(Z^\lambda_s)\Big)\dif s\\
&\quad+\mE\int^t_0\e^{ms}\Big(2|Z^\lambda_s|\,\Gamma^{0,1}_{R,\lambda}(g)(Z^\lambda_s)+\Gamma^{0,2}_{0,\lambda}(g)(Z^\lambda_s)\Big)\dif s\\
&\leq |x_0-y_0|^2+c(\e^{mt}-1)/m,
\end{align*}
where $c>0$ is independent of $m$.
From this we derive that for $m$ large enough,
$$
\mP(|Z^\lambda_T(x_0)-y_0|>a)\leq \frac{\mE|Z^\lambda_T(x_0)-y_0|^2}{a^2}\leq\frac{\e^{-mT}|x_0-y_0|^2}{a^2}+\frac{c(1-\e^{-mT})}{ma^2}\leq 1/4,
$$
which together with \eqref{RR} and \eqref{RR0} yields that
$$
\mP(|Z_T(x_0)-y_0|>a)\leq\mP(|Z_T(x_0)-y_0|>a, T<\tau_\lambda)+\mP(T\geq\tau_\lambda)\leq 1/2.
$$
The proof is complete.
\end{proof}

For each $m\in\mN$, let $\chi_m(x)$ be the cutoff function in \eqref{Cut}. Let
$$
\sigma_m(x):=\sigma(x\chi_m(x)),\ b_m(x):=\chi_m(x)b(x),\ g_m(x,z):=g(x\chi_m(x),z).
$$
In the following we assume that one of the following conditions holds:
\begin{enumerate}[{\bf (C1)}]
\item for each $m\in\mN$, $(\sigma^m, b^m,g^m)$ satisfies the assumptions of Theorem \ref{main1}.
\item for each $m\in\mN$, $(0, b^m,g^m)$ satisfies the assumptions of Theorem \ref{main3}.
\end{enumerate}
Let $P_t\varphi(x):=\mE \varphi(X_t(x))$. We have

\bl\label{strong}
Under {\bf (C1)} or {\bf (C2)}, and \eqref{EQ3} and \eqref{EQ4},
the semigroup $P_t$ has the strong Feller property and irreducibility.
\el
\begin{proof}
(i) Let $X^m_t(x)$ be the solution of SDE (\ref{sde4}) corresponding to $(\sigma_m,b_m,g_m)$.
In the case of {\bf (C1)}, by \eqref{feller}, for any bounded measurable function $f$ and $t>0$,
\begin{align}
x\mapsto \mE f(X^m_t(x))\mbox{ is continuous.}\label{EB6}
\end{align}
In the case of {\bf (C2)}, by the gradient estimate \eqref{GD4}, we still have \eqref{EB6}.

Now fix $K>0$. For $x\in\mR^d$ and $m>K$, define a stopping time
$$
\tau^x_m:=\left\{t\geq 0: |X_t(x)|\geq m\right\}.
$$
By Chebyshev's inequality and \eqref{BB7}, we have
\begin{align}
\lim_{m\to\infty}\sup_{|x|\leq K}\mP(t\geq\tau^x_m)\leq\lim_{m\to\infty}\sup_{|x|\leq K}\mE\left(\sup_{s\in[0,t]}|X_s(x)|^{\vartheta}\right)/m^{\vartheta}=0.\label{EB7}
\end{align}
Moreover, by the local uniqueness of solutions to SDE (\ref{sde4}), we have
$$
X_t(x)=X^m_t(x),\ \ |x|\leq K,\ \ t\in[0,\tau^x_m).
$$
Let $f$ be a bounded measurable function. For any $x,y\in B_K$, we have
\begin{align*}
&\quad|\mE(f(X_t(x))-f(X_t(y)))|\\
&\leq\big|\mE\big(f(X_t(x))-f(X_t(y))1_{t<\tau^x_m\wedge\tau^y_m}\big)\big|+2\|f\|_\infty \mP(t\geq\tau^x_m\wedge\tau^y_m)\\
&=\big|\mE\big(f(X^m_t(x))-f(X^m_t(y))1_{t<\tau^x_m\wedge\tau^y_m}\big)\big|+2\|f\|_\infty \mP(t\geq\tau^x_m\wedge\tau^y_m)\\
&\leq\big|\mE\big(f(X^m_t(x))-f(X^m_t(y))\big)\big|+4\|f\|_\infty \mP(t\geq\tau^x_m\wedge\tau^y_m)\\
&\leq\big|\mE\big(f(X^m_t(x))-f(X^m_t(y))\big)\big|+4\|f\|_\infty \big(\mP(t\geq\tau^x_m)+\mP(t\geq\tau^y_m)\big),
\end{align*}
which together with (\ref{EB6}), (\ref{EB7}) yields the continuity of $x\mapsto \mE(f(X_t(x)))$.

(ii) For the irreducibility, it suffices to prove that for any $T$, $a>0$ and $x_0, y_0\in\mR^d$,
\begin{align}
\mP\big(|  X_T(x_0)-y_0|\leq a)>0.    \label{01}
\end{align}
In the case of {\bf (C1)}, we use Lemma \ref{Le72} and Girsanov's transformation to show \eqref{01}, see \cite{Re-Wu-Zh}.
Let $Z_t(x_0)$ solve SDE (\ref{Eq5}) and
set for $K>0$,
$$
\tau_K:=\inf\{t: |Z_t(x_0)|\geq K\}.
$$
By (\ref{EB4}) and (\ref{EB5}), we may fix $N$ and $m$ large enough so that
\begin{align}
\label{qandp}
\mP(\tau_K\leq T)+\mP(|Z_T(x_0)-y_0|>a)<1.
\end{align}
Define
$$
U_t:=-m~ \sigma(Z_t)^{-1}(Z_t-y_0),\ \ \tilde W_t:=W_t+\int^{t\wedge\tau_K}_0 U_s\dif s,
$$
and
$$
\cE_T:=\exp\left(\int_0^{T\wedge\tau_K}U_s\dif W_s-\frac{1}{2}\int_0^{T\wedge\tau_K}|U_s|^2\dif s\right).
$$
Since $|U_{t\wedge \tau_K}|^2$ is bounded, we have $\mE[\cE_T]=1$.
By Girsanov's theorem (see \cite[Theorem 132]{R}), under the new probability measure $\mQ:=\cE_T\mP$, $\tilde W_t$ is still a Brownian motion,
and $N(\dif t, \dif z)$ is a Poisson random measure with the same compensator $\dif t\nu(\dif z)$. In view of (\ref{qandp}), we also have
\begin{align*}
\mQ(\{\tau_K\leq T\}\cup\{|Z_T(x_0)-y_0|>a\})<1.
\end{align*}
Note that the solution $Z_t$ of (\ref{Eq5}) also solves the following SDE:
\begin{align*}
Z_{t\wedge\tau_K}&=x_0+\int_0^{t\wedge \tau_K}  b(Z_s)\dif s+\int_0^{t\wedge\tau_K} \sigma(Z_s)\dif \tilde W_s\\
&\quad+\int_0^{t\wedge \tau_K}\!\!\!\int_{|z|<R}  g(Z_{s-},z)\tilde{N}(\dif s, \dif z)+\int_0^{t\wedge \tau_K}\!\!\!\int_{|z|\geq R}  g(Z_{s-},z)N(\dif s, \dif z).
\end{align*}
Set
$$
\theta_K:=\inf\{t: |X_t|\geq K\}.
$$
Then the law uniqueness for (\ref{sde4}) yields that the law of $\{(X_t1_{\{t\leq\theta_K\}})_{t\in [0,T]},\theta_K\}$ under $\mP$ is the same as that of
$\{(Z_t1_{\{t\leq \tau_K\}})_{t\in [0,T]},\tau_K\}$ under $\mQ$. Hence
\begin{align*}
\mP(|X_T(x_0)-y_0|>a)&\leq\mP(\{\theta_K\leq T\}\cup\{\theta_K\geq T, |X_T(x_0)-y_0|>a\})\\
&=\mQ(\{\tau_K\leq T\}\cup\{\tau_K\geq T, |Z_T(x_0)-y_0|>a\})\\
&\leq \mQ(\{\tau_K\leq T\}\cup\{|Z_T(x_0)-y_0|>a\})<1,
\end{align*}
which implies (\ref{01}).

In the case of {\bf (C2)}, let $D_m:=\{x: |x|<m\}$ be a ball containing $x_0$ and $B_a(y_0)$. We have
\begin{align*}
\mP\big(| X_T(x_0)-y_0|\leq a)&\geq\mP\big(X_T(x_0)\in B_a(y_0); T<\tau_{D_m})\\
&=\mP\big(X^m_T(x_0)\in B_a(y_0); T<\tau_{D_m}),
\end{align*}
where $X^m_T(x_0)$ is the solution of SDE (\ref{sde4}) corresponding to $(0,b_m,g_m)$.
By \eqref{GD5} and Theorem \ref{Th79}, since $\varrho_i(t,r)=c_it(t^{1/\alpha}+r)^{-d-\alpha}$ with $c_1<c_2$ satisfy {\bf (H$^\varrho$)},
we get \eqref{01}. The proof is complete.
\end{proof}

Now we can give

\begin{proof}[Proof of Theorem \ref{main22}]
By \eqref{BB6}, the existence of invariant probability measures for $P_t$ follows by the standard Bogoliov-Krylov's argument.
The uniqueness is a direct consequence of the strong Feller property and irreducibility. Moreover, still by the strong Feller property and irreducibility,
we can derive easily that for any $y\in\mR^d$ and $r,t>0$,
$$
\inf_{x\in B_r} \mP\big(X_t(x)\in B_r(y)\big)>0.
$$
Combing this with \eqref{BB6} and \cite[Theorem 2.5]{Go-Ma}, we get the desired results.
\end{proof}
\subsection{SDEs with singular and dissipative drifts}

In this subsection we study the ergodicity of SDE \eqref{sde4} with  singular and dissipative drift.
The main idea is to use Zvonkin's transformation to kill the singular part.
First of all, we consider the  case of non-degenerate diffusion,
and show the following non explosion and Krylov's estimate.
\bl\label{Le74}
Under  {\bf (H$^\sigma$)},  {\bf (H$^b$)} and \eqref{JF},
any solution $X_t(x)$ to SDE (\ref{sde4}) does not explode.
Moreover, for any $T>0$ and $f\in L^{p'}(\mR^d)$ with $p'>d$,
\begin{align}
\mE\left(\int_0^T\!f\big(X_s(x)\big)\dif s\right)\leq c(|x|+1)\|f\|_{p'},    \label{kry6}
\end{align}
where $c>0$ is independent of $x$.
\el
\begin{proof}
For $n>0$, let $\tau_n:=\inf\{t\geq 0: |X_t|\geq n\}$.
By Lemma \ref{lem}, for any $T>0$, $p'>d$ and $\delta>0$, there exists a constant $c_\delta>0$ such that
for any $f\in L^q(\mR^d)$,
\begin{align}\label{B09}
\mE\left(\int_0^{T\wedge\tau_n}\!\!f(X_s)\dif s\right)\leq
\left(c_\delta+\delta\mE\left(\int_0^{T\wedge\tau_n}|b_1+b_2|(X_s)\dif s\right)\right)\|f\|_{p'}.
\end{align}
Since $b_1\in L^{p}(\mR^d)$ with $p>d$, for every $\delta_0>0$,
we can take $f=|b_1|$ and choose $\delta$ small enough such that $\delta\|b_1\|_{p}<\delta_0$ in the above inequality to get
\begin{align}\label{BB4}
\begin{split}
&\mE\left(\int_0^{T\wedge\tau_n}\!\!b_1(X_s)\dif s\right)\leq
c_{\delta_0}+\delta_0\mE\left(\int_0^{T\wedge\tau_n}|b_2(X_s)|\dif s\right)\\
&\qquad\qquad\leq c_{\delta_0}+\kappa_3\delta_0\mE\left(\int_0^{T\wedge\tau_n}(1+|X_s|^2)^{(1+r)/2}\dif s\right).
\end{split}
\end{align}
On the other hand,
let $h(x):=\sqrt{1+|x|^2}$. By It\^o's formula, we have
\begin{align}\label{BB44}
\mE h(X_{t\wedge\tau_n})=h(x)+\mE\int^{t\wedge\tau_n}_0[\sL^\sigma_2 h+\sL^b_1 h+\sL^g_{\nu} h](X_s)\dif s.
\end{align}
As the calculations in Lemma \ref{Le71}, by the assumptions, we have
\begin{align*}
\sL^\sigma_2 h(x)&\leq\tfrac{1}{2}(\sigma^{ik}\sigma^{ik})(x)(1+|x|^2)^{-1/2}\leq c,\\
\sL^b_1h(x)&\leq(-\kappa_1|x|^{2+r}+\kappa_2)(1+|x|^2)^{-1/2}+|b_1(x)|\\
&\leq-\kappa_1(1+|x|^2)^{(1+r)/2}/2+c+|b_1(x)|,
\end{align*}
and
\begin{align*}
\sL^g_\nu h(x)&=\int_{\mR^d}\Big[h\big(x+g(x,z)\big)-h(x)-1_{|z|\leq R}g(x,z)\cdot\nabla h(x)\Big]\nu(\dif z)\\
&\leq\int_{|z|<R}|g(x,z)|^2\nu(\dif z)+\int_{|z|\geq R}|g(x,z)|\nu(\dif z)\leq c.
\end{align*}
Hence,  by \eqref{BB44} and \eqref{BB4} with $\delta_0$ small enough, we obtain
\begin{align*}
\mE (1+|X_{T\wedge\tau_n}|^2)^{1/2}
&\leq (1+|x|^2)^{1/2}-\frac{\kappa_1}{2}\mE\int^{T\wedge\tau_n}_0 (1+|X_s|^2)^{(1+r)/2}\dif s\\
&\quad+\mE\int^{T\wedge\tau_n}_0 |b_1(X_s)|\dif s+ct\\
&\leq (1+|x|^2)^{1/2}-\frac{\kappa_1}{4}\mE\int^{T\wedge\tau_n}_0 (1+|X_s|^2)^{(1+r)/2}\dif s+c_T,
\end{align*}
which implies that $\lim_{n\to\infty}\tau_n=\infty$ and
$$
\mE (1+|X_T|^2)^{1/2}+\frac{\kappa_1}{4}\mE\int^T_0 (1+|X_s|^2)^{(1+r)/2}\dif s\leq(1+|x|^2)^{1/2}+c_T.
$$
Substituting this into \eqref{B09} and \eqref{BB4}, we obtain \eqref{kry6}.
\end{proof}
To perform Zvonkin's transformation, we need to solve a related elliptic equation,
which is a consequence of Theorem \ref{pde}.
\bt\label{Th51}
Suppose that {\bf (H$^\sigma$)} holds and $b\in L^p(\mR^d)$ for some $p>d$, and
\begin{align*}
\Gamma^{0,2}_{0,R}(g)\in L^\infty(\mR^d),\ \ \lim_{\eps\to 0}\|\Gamma^{0,2}_{0,\eps}(g)\|_\infty=0.
\end{align*}
Then for some $\lambda_1\geq 1$ large enough and for all $\lambda\geq\lambda_1$ and $f\in L^p(\mR^d)$,
there exists a unique solution $u\in H^2_p$
to the following elliptic equation:
\begin{align}\label{ELL}
(\sL^\sigma_2-\lambda) u+\sL^g_{\nu,R}u+\sL_1^{b}u =f,
\end{align}
%where $\sL^g_{\nu,R}u$ is defined by \eqref{UY7},
and for any $p'\in[p,\infty]$ and $\vartheta\in(0,2)$ with $\frac{d}{p}<2-\vartheta+\frac{d}{p'}$,
\begin{align}
\lambda^{\frac{1}{2}(2-\vartheta+\frac{d}{p'}-\frac{d}{p})}\|u\|_{\vartheta,p'}+\|\nabla^2u\|_{p}\leq c\|f\|_p.\label{es11}
\end{align}
%where $c=c(d,p,\vartheta, p', \lambda_1,\|b\|_{p})>0$.
\et
\begin{proof}
As usual, it suffices to show the apriori estimate \eqref{es11}. Let $u\in H^2_p$ solve \eqref{ELL}.
Let $T>0$ and $\phi(t)$ be a nonnegative and nonzero smooth function with support in $(0,T)$. Let $\bar u(t,x):=u(x)\phi(t)$. It is easy to see that
$\bar u$ satisfies the following parabolic equation:
$$
\p_t\bar u+(\sL^\sigma_2-\lambda) \bar u+\sL^g_{\nu,R}\bar u+\sL_1^{b}\bar u=u\phi'+f\phi.
$$
Thus, by Theorem \ref{pde}, there is a $\lambda_0\geq 1$ depending on $\|b\|_p$ and $\|\Gamma^{0,2}_{0,R}(g)\|_\infty$
such that for all $\lambda\geq\lambda_0$, $p'\in[p,\infty]$ and $\vartheta\in(0,2)$ with $\frac{d}{p}<2-\vartheta+\frac{d}{p'}$,
$$
\lambda^{\frac{1}{2}(2-\vartheta+\frac{d}{p'}-\frac{d}{p})}\|\bar u\|_{\mH^{\vartheta,\infty}_{p'}(T)}+\|\nabla^2\bar u\|_{\mL^\infty_p(T)}\leq c\|u\phi'+f\phi\|_{\mL^\infty_p},
$$
which implies that
\begin{align}\label{KG1}
\lambda^{\frac{1}{2}(2-\vartheta+\frac{d}{p'}-\frac{d}{p})}\|u\|_{\vartheta,p'}+\|\nabla^2 u\|_{p}\leq
c\|\phi\|^{-1}_{\infty}\Big(\|u\|_p\|\phi'\|_\infty+\|f\|_p\|\phi\|_{\infty}\Big).
\end{align}
Letting $p'=p$ and $\vartheta=0$ in \eqref{KG1} and choosing $\lambda_1\geq \lambda_0$ large enough, we get
$$
\|u\|_p\leq c\|f\|_p.
$$
Finally, substituting this into \eqref{KG1}, we obtain the desired estimate \eqref{es11}.
\end{proof}

Below we assume that {\bf (H$^\sigma$)} holds and for some $p>d$,
$$
b_1,\ |\nabla\sigma|,\ (\Gamma^{1,2}_{0,R}(g))^{1/2}\in L^p(\mR^d),\ \
\Gamma^{0,2}_{0,R}(g)\in L^\infty(\mR^d),\ \lim_{\eps\to 0}\|\Gamma^{0,2}_{0,\eps}(g)\|_{\infty}=0.
$$
Now consider the following elliptic equation system:
\begin{align*}
(\sL^\sigma_2-\lambda) u+\sL^g_{\nu,R}u+\sL_1^{b_1}u =b_1.
\end{align*}
By \eqref{es11}, there are $c, \lambda_1\geq 1$ such that for all $\lambda\geq\lambda_1$,
\begin{align}\label{BB8}
\|u\|_\infty+\|\nabla u\|_\infty\leq c\lambda^{\frac{1}{2}(\frac{d}{p}-1)}.
\end{align}
Define
$$
\Phi(x):=x+ u(x).
$$
By \eqref{BB8} with $\lambda$ large enough, the map $x\rightarrow\Phi(x)$ forms a $C^1$-diffeomorphism and
\begin{align*}
1/2\leq \|\nabla \Phi\|_{\infty},\|\nabla\Phi^{-1}\|_{\infty}\leq 2,
\end{align*}
where $\Phi^{-1}$ is the inverse of $\Phi$.

By Lemma \ref{Le74} and Theorem \ref{Th51},
the following result can be shown in the same way as in Lemma \ref{zvon}. We omit the details.
\bl\label{Le53}
$X_t$ solves SDE \eqref{sde4} if and only if $Y_t:=\Phi(X_t)$ solves
\begin{align*}
\begin{split}
\dif Y_t=\tilde\sigma(Y_t)\dif W_t+\tilde b(Y_t)\dif t+\int_{|z|<R}\tilde g(Y_{t-},z)\tilde N(\dif t,\dif z)+\int_{|z|\geq R}\tilde g(Y_{t-},z)N(\dif t,\dif z),
\end{split}
\end{align*}
where $y:=\Phi(x)$ and
\begin{align*}
\begin{split}
&\tilde \sigma(y):=\big(\nabla\Phi\cdot \sigma\big)\circ\Phi^{-1}(y),\quad \tilde b(y):=(\lambda u+\nabla\Phi\cdot b_2)\circ\Phi^{-1}(y),\\
&\qquad\quad\tilde g(y,z):=\Phi\big(\Phi^{-1}(y)+g\big(\Phi^{-1}(y),z\big)\big)-y.
\end{split}
\end{align*}
\el
The following proposition shows that the dissipativity \eqref{diss} is preserved under Zvonkin's transformation.
\bp\label{Pr77}
Under \eqref{diss}, for $\lambda$ large enough, there are $\tilde\kappa_1,\tilde\kappa_2,\tilde\kappa_3>0$ such that for all $y\in\mR^d$,
\begin{align*}
\<y, \tilde b(y)\>\leq -\tilde\kappa_1|y|^{2+r}+\tilde\kappa_2\quad\text{and}\quad|\tilde b(y)|\leq \tilde\kappa_3(1+|y|^{1+r}).
\end{align*}
\ep
\begin{proof}
Noticing that
$$
y=\Phi^{-1}(y)+ u\big(\Phi^{-1}(y)\big),\ \ \nabla\Phi(x)=\mI+\nabla u(x),
$$
by the definition of $\tilde b$ and \eqref{diss},
we have
\begin{align*}
\<y, \tilde b(y)\>&=\lambda\<y, u(\Phi^{-1}(y))\>+\<y,b_2(\Phi^{-1}(y))\>+\<y,(b_2\nabla u)(\Phi^{-1}(y))\>\\
&\leq \lambda\| u\|_{\infty}|y| +\big\langle\Phi^{-1}(y),  b_2\big(\Phi^{-1}(y)\big)\big\rangle+\|u\|_{\infty}\cdot|b_2\big(\Phi^{-1}(y)\big)|\\
&\qquad\qquad\qquad+\|\nabla u\|_\infty|y|\cdot|b_2\big(\Phi^{-1}(y)\big)|\\
&\leq \lambda\|u\|_\infty |y|-\kappa_1|\Phi^{-1}(y)|^{2+r}+\kappa_2\\
&\quad+\kappa_3(1+|\Phi^{-1}(y)|^2)^{(1+r)/2}(\| u\|_{\infty}+\|\nabla u\|_\infty|y|)\\
&\leq \lambda\|u\|_\infty |y|-\kappa_1(|y|-\|u\|_\infty)^{2+r}+\kappa_2\\
&\quad+\kappa_3(1+(|y|+\|u\|_\infty)^2)^{(1+r)/2}(\| u\|_{\infty}+\|\nabla u\|_\infty|y|)\\
&\leq (c_1\|\nabla u\|_\infty-\tfrac{\kappa_1}{2})|y|^{2+r}+c_\lambda,
\end{align*}
where $c_1$ only depends on $\kappa_3$ and $r$. By \eqref{BB8} with $\lambda$ large enough so that
$c_1\|\nabla u\|_\infty\leq\tfrac{\kappa_1}{4}$, we get the first estimate. The second estimate is easy.
\end{proof}

Now we can give

\begin{proof}[Proof of Theorem \ref{main2}]

By  Lemma \ref{Le53}, Proposition \ref{Pr4} and Proposition \ref{Pr77}, one can use Theorem \ref{main22} to conclude the same conclusions.
It remains to show that the invariant probability measure $\mu$ has a density $\rho\in L^q(\mR^d)$ with $q<d/(d-1)$. By Zvonkin's transformation Lemma \ref{Le53},
we may assume $b_1=0$. Let $f\in C^\infty_0(\mR^d)$, and for $p>d$ let $u\in H^2_p$ solve the following elliptic equation:
$$
(\sL^\sigma_2-\lambda) u+\sL^g_{\nu,R}u=f.
$$
Let $u_n=u*\phi_n$ be the mollifying approximation of $u$ and define
$$
f_n:=(\sL^\sigma_2-\lambda) u_n+\sL^g_{\nu,R}u_n.
$$
By It\^o's formula, we have
$$
\mE u_n(X_T)=u_n(x)+\mE\left(\int^T_0(f_n+\lambda u_n+b\cdot\nabla u_n)(X_t)\dif t\right).
$$
Noticing that
$$
\|f_n-f\|_p\leq c\|u_n-u\|_{2,p},
$$
by Krylov's estimate \eqref{kry6} and \eqref{diss} we have
\begin{align*}
&\mE\left(\int^T_0f(X_t)\dif t\right)=\lim_{n\to\infty}\mE\left(\int^T_0f_n(X_t)\dif t\right)\\
&\quad\leq(\lambda+2)\|u_n\|_\infty+\|\nabla u_n\|_\infty\mE\left(\int^T_0|b|(X_t)\dif t\right)\\
&\quad\leq(\lambda+2)\|u\|_\infty+\kappa_3\|\nabla u\|_\infty\mE\left(\int^T_0(1+|X_t|^{1+r})\dif t\right),
\end{align*}
which yields by \eqref{BB7} and \eqref{es11} that
$$
\mE\left(\int^T_0f(X_t)\dif t\right)\leq c(1+|x|+T)\|f\|_p,
$$
where $c$ is independent of $T$ and $x$. By \eqref{Erg} we get for any $p>d$,
$$
\mu(f)\leq c\|f\|_p,\ \ f\in C^\infty_0(\mR^d),
$$
which in turn implies that $\mu$ has a density $\rho\in L^{p/(p-1)}(\mR^d)$. The proof is complete.
\end{proof}

The proof of  Theorem \ref{main4} is similar. We sketch the proof below. As in Lemma \ref{Le74}, the following lemma can be proven by Lemma \ref{Le63}.
\bl\label{Le78}
Under  {\bf (H$^g$)},  {\bf (H$^b$)} and \eqref{JF},
any solution $X_t(x)$ to SDE (\ref{sde4}) does not explode.
Moreover, for any $T>0$ and $f\in L^{p'}(\mR^d)$ with $p'>d/(\alpha-1)$,
$$
\mE\left(\int_0^T\!f\big(X_s(x)\big)\dif s\right)\leq c(|x|+1)\|f\|_{p'},
$$
where $c>0$ is independent of $x$.
\el
We also have the solvability of the  following  non-local elliptic equation.
\bt\label{Th78}
Let $\alpha\in(1,2)$ and $\sL^\kappa_\alpha$ be defined by \eqref{49}, where $\kappa$ satisfies \eqref{Con2}.
Let $\theta\in(0,1)$ and $\gamma\in(0,\alpha)$ with $\gamma+\theta<\alpha+\eps$, where $\eps$ is the same as in (vii) of Theorem \ref{heat}.
Suppose that $b\in H^\theta_p(\mR^d)$ for some $p>d/(\alpha+\theta-1)$.
Then for some $\lambda_1\geq 1$ large enough and for all $\lambda\geq\lambda_1$ and $f\in H^\theta_p(\mR^d)$,
there exists a unique solution $u\in H^{\gamma+\theta}_p$ to the following nonlocal elliptic equation:
\begin{align}\label{eq00}
(\sL^\kappa_\alpha-\lambda) u+\sL_1^{b}u =f,
\end{align}
so that
\begin{align}
\lambda^{1-\frac{\gamma}{\alpha}}\|u\|_{\gamma+\theta, p}\leq c\|f\|_{\theta,p}.\label{es1101}
\end{align}
Moreover, for any $\vartheta\in(0,\alpha)$ and $p'\in[p,\infty]$ with $\frac{d}{p}<\alpha-\vartheta+\frac{d}{p'}$,
\begin{align}
\lambda^{\frac{1}{\alpha}(\alpha-\vartheta+\frac{d}{p'}-\frac{d}{p})}\|u\|_{\vartheta, p'}\leq c\|f\|_p.\label{es101}
\end{align}
\et
\begin{proof}
We show the apriori estimate \eqref{es1101} and \eqref{es101}. Suppose $u\in H^{\gamma+\theta}_p$ satisfies \eqref{eq00}.
Let $T>0$ and $\phi(t)$ be a nonnegative and nonzero smooth function with support in $(0,T)$. Let $\bar u(t,x):=u(x)\phi(t)$. Then
$$
\p_t \bar u+(\sL^g_\nu-\lambda+\sL^b_1)\bar u=u\phi'+f\phi.
$$
By Theorem \ref{Th63}, we have
$$
\lambda^{1-\frac{\gamma}{\alpha}}\|\bar u\|_{\mH^{\gamma+\theta, \infty}_p(T)}\leq c\|u\phi'+f\phi\|_{\mH^{\theta,\infty}_p(T)},
$$
which implies that
$$
\lambda^{1-\frac{\gamma}{\alpha}}\|u\|_{\gamma+\theta, p}\|\phi\|_\infty\leq c\Big(\|u\|_{\theta,p}\|\phi'\|_\infty+\|f\|_{\theta,p}\|\phi\|_\infty\Big).
$$
Letting $\lambda$ be large enough, we get \eqref{es1101}.
On the other hand, by \eqref{NR1} we also have
$$
\lambda^{\frac{1}{\alpha}(\alpha-\vartheta+\frac{d}{p'}-\frac{d}{p})}\|\bar u\|_{\mH^{\gamma', \infty}_{p'}(T)}\leq c\|u\phi'+f\phi\|_{\mL^p(T)},
$$
which also implies \eqref{es101} as above.
\end{proof}

Below we assume that $\sigma\equiv0$, $\nu(\dif z)=|z|^{-d-\alpha}\dif z$, {\bf (H$^g$)} holds, and
for some $\theta\in(1-\alpha/2,1)$ and $p>2d/\alpha$,
$$
(I-\Delta)^{\theta/2}b_1, (\Gamma^{1,2}_{0,R}(g))^{1/2}\in L^p(\mR^d),\ \
\Gamma^{0,2}_{0,R}(g)\in L^\infty(\mR^d).
$$
Consider the following nonlocal elliptic equation system:
\begin{align*}
(\sL^g_{\nu}-\lambda)u+\sL_1^{b_1}u =b_1.
\end{align*}
By \eqref{es101}, there are $c, \lambda_1\geq 1$ such that for all $\lambda\geq\lambda_1$,
\begin{align}\label{B0B8}
\|u\|_\infty+\|\nabla u\|_\infty\leq c\lambda^{\frac{1}{\alpha}(\frac{d}{p}+1-\alpha)}.
\end{align}
Define
$$
\Phi(x):=x+ u(x).
$$
By \eqref{B0B8} with $\lambda$ large enough, the map $x\rightarrow\Phi(x)$ forms a $C^1$-diffeomorphism and
\begin{align*}
1/2\leq \|\nabla \Phi\|_{\infty},\|\nabla\Phi^{-1}\|_{\infty}\leq 2,
\end{align*}
where $\Phi^{-1}$ is the inverse of $\Phi$.

By Lemma \ref{Le78} and Theorem \ref{Th78},
the following result can be shown in the same way as in Lemma \ref{zvon}. We omit the details.
\bl\label{Le710}
$X_t$ solves SDE \eqref{sde4} with $\sigma=0$ if and only if $Y_t:=\Phi(X_t)$ solves
\begin{align*}
\begin{split}
\dif Y_t=\tilde b(Y_t)\dif t+\int_{|z|<R}\tilde g(Y_{t-},z)\tilde N(\dif t,\dif z)+\int_{|z|\geq R}\tilde g(Y_{t-},z)N(\dif t,\dif z),
\end{split}
\end{align*}
where $y:=\Phi(x)$ and
$$
\tilde b(y):=(\lambda u+\nabla\Phi\cdot b_2)\circ\Phi^{-1}(y),\quad\tilde g(y,z):=\Phi\big(\Phi^{-1}(y)+g\big(\Phi^{-1}(y),z\big)\big)-y.
$$
\el
\begin{proof}[Proof of Theorem \ref{main4}]

By  Lemma \ref{Le710}, Proposition \ref{Pr4} and Proposition \ref{Pr77},  the result follows by Theorem \ref{main22}.
As for the conclusion that $\mu$ has a density $\rho\in L^q(\mR^d)$ with $q<d/(d-\alpha+1)$, it follows by Theorem \ref{Th78}
and the same argument as
used in the proof of Theorem \ref{main2}.
\end{proof}

\subsection{Positivity of Dirichlet heat kernel}
Let $\rho(t,x,y)$ be a family of jointly continuous transition probability density functions in $\mR^d$. Let $(X,\mP_x)_{x\in\mR^d}$
be the associated Markov processes, that is,
$\mP_x(X_0=x)=1$ and for any $t>0$,
$$
\int_A\rho(t,x,y)\dif y=\mP_x(X_t\in A),\ \ A\in\sB(\mR^d).
$$
Let $D$ be a domain (bounded open subset of $\mR^d$), and $\tau_D:=\{t>0: X_t\notin D\}$ be the exit time of $X$ from $D$.
Let $X^D$ be the killed Markov process outside $D$, and $P^D_t$ the transition probability of $X^D$, that is,
$$
P^D_t(x,A):=\mP_x(t<\tau_D; X_t\in A),\ \ A\in\sB(D).
$$
Define
$$
\rho^D(t,x,y):=\rho(t,x,y)-r^D(t,x,y), \
$$
where
$$
r^D(t,x,y):=\mE^x[\tau_D<t; \rho(t-\tau_D;X(\tau_D),y)].
$$
Let $\varrho_i(t,r):\mR_+\times\mR_+\to\mR_+$, $i=1,2$ be two continuous functions and satisfy that
\begin{enumerate}[{\bf (H$^\varrho$)}]
\item For each $t>0$, the map $r\mapsto\varrho_i(t,r)$ is decreasing, and for each $\delta>0$,
$$
\sup_{t>0, r>\delta}\varrho_2(t,r)<\infty,
$$
and there are $t_0=t_0(\delta)$ and $R=R(\delta)>0$
such that $t\mapsto\varrho_i(t,\delta)$ is increasing on $(0,t_0)$ and
$$
\varrho_1(t,\delta/R)>\varrho_2(t,\delta),\ \ t\in(0,t_0).
$$
\end{enumerate}
The following result is essentially due to Hunt (cf. \cite[Theorem 2.4]{Ch-Zh2}).
\bt\label{Th79}
Let $\varrho_1$ and $\varrho_2$ satisfy {\bf (H$^\varrho$)}. Suppose that
\begin{align}\label{He}
\varrho_1(t,|x-y|)\leq \rho(t,x,y)\leq \varrho_2(t,|x-y|).
\end{align}
Then $\rho^D$ is the transition probability density function of $X^D$, i.e., for any $t>0$,
$$
P^D_t(x,A)=\int_A\rho^D(t,x,y)\dif y,\ x\in\mR^d,\ A\in\sB(D).
$$
Moreover, $\rho^D$ is continuous and strictly positive on $\mR_+\times D\times D$ and for $0<s<t<\infty$ and $x,y\in\mR^d$,
\begin{align}\label{CK}
\rho^D(t,x,y)=\int_D \rho^D(s,x,z)\rho^D(t-s,z,y)\dif z.
\end{align}
\et
\begin{proof}
We only show the strict positivity of $\rho^D(t,x,y)$. The others are completely same as in \cite[Theorem 2.4]{Ch-Zh2}.
Fix $x,y\in D$ and let $d(y,\p D)$ be the distance of $y$ to the boundary $\p D$. Let $\delta\in(0,d(y,\p D))$ be given. By the assumption on $\varrho_i$,
there are $t_0=t_0(\delta)>0$ and $R>0$ such that
$$
\varrho_1(t,\delta/R)>\varrho_2(t,\delta),\ \ t\in(0,t_0).
$$
Hence, by the definition of $r^D$ , \eqref{He} and the assumptions of $\varrho_2$, we have for $t\in(0,t_0)$,
$$
r^D(t,x,y)\leq \mE^x[\tau_D<t; \varrho_2(t-\tau_D,|X(\tau_D)-y|)]\leq\varrho_2(t,\delta).
$$
Consequently, if $|x-y|\leq \delta/R\leq\delta<\rho(y,\p D)$, then
\begin{align}\label{GD7}
\rho^D(t,x,y)\geq\varrho_1(t,|x-y|)-\varrho_2(t,\delta)\geq \varrho_1(t,\delta/R)-\varrho_2(t,\delta)>0.
\end{align}
Now for any $t>0$ and $x,y\in D$. Let $\Gamma$ be a curve in $D$ connecting $x$ and $y$. Let $\delta:=\rho(\Gamma,\p D)$.
Let $n$ be large enough such that $t\leq n t_0(\delta)$ and there are points $a_0,a_1,\cdots, a_{n+1}$ on $\Gamma$ with $a_0=x$, $a_{n+1}=y$
and $a_i\in B(a_{i-1},\delta/(3R))$.
Notice that for $x_{i-1}\in B(a_{i-1},\delta/(3R))$ and $x_{i}\in B(a_{i},\delta/(3R))$,
$$
|x_{i}-x_{i-1}|\leq|x_{i}-a_i|+|a_{i-1}-x_{i-1}|+|a_i-a_{i-1}|\leq \delta/R.
$$
By C-K equation \eqref{CK} and \eqref{GD7}, we have
\begin{align*}
&\rho^D(t,x,y)=\int_D\cdots\int_D \rho^D(\tfrac{t}{n}, x,x_1)\cdots \rho^D(\tfrac{t}{n}, x_n,y)\dif x_1\cdots\dif x_n\\
&\geq \int_{B(a_1,\delta/(3R))}\cdots\int_{B(a_n,\delta/(3R))} \rho^D(\tfrac{t}{n}, x,x_1)\cdots \rho^D(\tfrac{t}{n}, x_n,y)\dif x_1\cdots\dif x_n>0.
\end{align*}
The proof is complete.
\end{proof}

\begin{center}
\bf Acknowledgement
\end{center}

The authors would like to thank Zhen-Qing Chen, Zenghu Li, Renming Song, Feng-Yu Wang, Yinchao Xie
and Guohuan Zhao for their quite useful conversations.
%\bigskip

\end{document}